\documentclass[11pt]{article}
\def\diff{\mathsf{d}}

\usepackage{amscd}      % to be able to use "CD" environment
\usepackage{amssymb}
\usepackage[intlimits]{amsmath}
\usepackage{theorem}

\newtheorem{prop}{Proposition}[section]
\newtheorem{lem}[prop]{Lemma}

\newtheorem{cor}[prop]{Corollary}
\newtheorem{them}[prop]{Theorem}

\theorembodyfont{\upshape}

\newtheorem{defn}[prop]{Definition}

\newtheorem{numrmk}[prop]{Remark}

\newtheorem{numex}[prop]{Example}

\newtheorem{rmk}{Remark}

\newenvironment{pf}{\begin{trivlist}\item[]{\sc Proof.}}%
           {\nolinebreak $\Box$ \end{trivlist}}

\newcommand{\noprint}[1]{}

\renewcommand{\tilde}{\widetilde}

\newcommand{\toto}{\rightrightarrows}

\newcommand{\rr}{{\mathbb R}}

\newcommand{\cala}{{\cal A}}
\def\lcf{\lbrack\! \lbrack}
\def\rcf{\rbrack\! \rbrack}
\newcommand{\alp }{\alpha }
\newcommand{\bet }{\beta }
\newcommand{\parrr}[1]{\frac{\partial }{\partial #1}}
\def\gpd{\rightrightarrows}

\newcommand{\ldiag}[1]%
       {\makebox[0cm]{${\scriptstyle#1}\downarrow\phantom{\scriptstyle#1}$}}
\newcommand{\ldiagup}[1]%
       {\makebox[0cm]{${\scriptstyle#1}\uparrow\phantom{\scriptstyle#1}$}}
\newcommand{\rdiag}[1]%
       {\makebox[0cm]{$\phantom{\scriptstyle#1}\downarrow{\scriptstyle#1}$}}
\newcommand{\sediagr}[1]%
       {\makebox[0cm]{$\phantom{\scriptstyle#1}\searrow{\scriptstyle#1}$}}
\newcommand{\nediagr}[1]%
       {\makebox[0cm]{$\phantom{\scriptstyle#1}\nearrow{\scriptstyle#1}$}}
\newcommand{\rdiagup}[1]%
       {\makebox[0cm]{$\phantom{\scriptstyle#1}\uparrow{\scriptstyle#1}$}}
\newcommand{\swdiag}[1]%
       {\makebox[0cm]{$\phantom{\scriptstyle#1}\swarrow{\scriptstyle#1}$}}
\newcommand{\sediag}[1]%
       {\makebox[0cm]{${\scriptstyle#1}\searrow\phantom{\scriptstyle#1}$}}
\newcommand{\nediag}[1]%
       {\makebox[0cm]{${\scriptstyle#1}\nearrow\phantom{\scriptstyle#1}$}}

\newcommand{\doublearrowstack}[2]%
                      {{{{\scriptstyle#1}\atop{\textstyle\longrightarrow}}\atop{{\textstyle\longrightarrow}\atop{\scriptstyle#2}}}}
\newcommand{\rightleftarrowstack}[2]%
                      {{{{\scriptstyle#1}\atop{\textstyle\longrightarrow}}\atop{{\textstyle\longleftarrow}\atop{\scriptstyle#2}}}}
\newcommand{\leftrightarrowstack}[2]%
                      {{{{\scriptstyle#1}\atop{\textstyle\longleftarrow}}\atop{{\textstyle\longrightarrow}\atop{\scriptstyle#2}}}}

\newcommand{\overtoparrow}%
{\makebox[0cm]{\beginpicture \setcoordinatesystem units
<.8cm,.4cm> point at 0 0 \setplotarea x from -3 to 3, y from 0 to
1 \setquadratic \plot -3 0 0 1 3 0 / \put{\vector(3,-1){0}}[Bl] at
3 0
\endpicture}}

\newcommand{\underbottomarrow}%
{\makebox[0cm]{\beginpicture \setcoordinatesystem units
<.8cm,.4cm> point at 0 0 \setplotarea x from -3 to 3, y from 0 to
1 \setquadratic \plot -3 1 0 0 3 1 / \put{\vector(3,1){0}}[Bl] at
3 1
\endpicture}}

\newcommand{\ses}[5]%
{0\longrightarrow#1\stackrel{#2}{ \longrightarrow}#3\stackrel{#4}{
\longrightarrow}#5\longrightarrow0}

\newcommand{\dt}[6]%
{#1\stackrel{#2}{longrightarrow}#3
\stackrel{#4}{\longrightarrow}#5 \stackrel{#6}{\longrightarrow}
#1[1]}

\newcommand{\cat}[1]%
{(\mbox{\rm #1})}

\newcommand{\gm}{\Gamma }
\newcommand{\lon }{\longrightarrow }
\newcommand{\backl}{\mathbin{\vrule width1.5ex height.4pt\vrule height1.5ex}}
\newcommand{\per}{\backl}
\newcommand{\be}{\begin{eqnarray*}}
\newcommand{\ee}{\end{eqnarray*}}
\newcommand{\smalcirc}{\mbox{\tiny{$\circ $}}}

\newcommand{\half}{\frac{1}{2}}

\let\Vec=\overrightarrow
\let\ceV=\overleftarrow

\begin{document}
\title{{\bf Universal lifting theorem and quasi-Poisson groupoids}}
\author{David Iglesias Ponte$^1$\thanks{ Research partially supported by
MCYT grant BFM2003-01319. }
, Camille Laurent-Gengoux$^2$, Ping Xu$^1$
 \thanks{ Research partially supported by NSF
       grant DMS03-06665 and NSA grant 03G-142. }
\\[5pt] {\small\it $^1$ Department of Mathematics, Pennsylvania State University}\\[5pt]
{\small\it $^2$ D\'epartement de Math\'ematiques,  Universit\'e de
Poitiers}\\[5pt]
{\small\it e-mail: iglesias@math.psu.edu,
Camille.Laurent@math.univ-poitiers.fr, ping@math.psu.edu} }
\date{}

\sloppy \maketitle

\begin{abstract}
We prove the universal lifting theorem:  for an $\alpha$-simply
connected and $\alpha$-connected Lie groupoid $\gm$ with Lie
algebroid $A$,  the graded Lie algebra of multi-differentials on $A$  is
isomorphic to that  of multiplicative multi-vector fields on
$\gm$. As a consequence, we obtain the
integration theorem for a quasi-Lie bialgebroid, which generalizes
various  integration theorems in the literature in  special cases.

The  second goal of the paper is  the  study of basic properties
of quasi-Poisson groupoids. In particular, we prove that a group
pair $(D, G)$ associated to a Manin quasi-triple $(\mathfrak d,
\mathfrak g, \mathfrak h)$ induces a quasi-Poisson groupoid on the
transformation groupoid $G\times D/G\toto D/G$. Its momentum map
corresponds exactly with the $D/G$-momentum map of Alekseev and
Kosmann-Schwarzbach.
\end{abstract}

\tableofcontents

\section{Introduction}
The notion of quasi-Lie bialgebroids, first introduced by
Roytenberg \cite{Roy},  is a  natural generalization of Lie
bialgebroids \cite{MackenzieX:1994}. It also generalizes
Drinfeld's quasi-Lie bialgebras \cite{Drinfeld}, the classical
limit of quasi-Hopf algebras. Roytenberg defines   quasi-Lie
bialgebroids using Lie algebroid analogues of Manin quasi-triples.
Equivalently, one may define a quasi-Lie bialgebroid as follows.
Given a Lie algebroid $A$, it is known that the anchor map
together with the bracket on $\gm (A)$ extends to a graded Lie
bracket on $\oplus _k\gm (\wedge^k A)$, which makes it into a
Gerstenhaber algebra  $ (\oplus _k \gm (\wedge^k A), \lcf \cdot,
\cdot \rcf, \wedge )$ \cite{Xu:1999}. Then a quasi-Lie
bialgebroid is a Lie algebroid $A$ equipped with a degree
$1$-derivation $\delta :\gm (\wedge ^\bullet A)\to \gm (\wedge
^{\bullet +1}A)$ with respect to the Gerstenhaber algebra
structure on $\oplus _k\gm (\wedge ^k A)$ such that $\delta^2
=\lcf \Omega , \cdot \rcf $ for some $\Omega \in \gm (\wedge^3 A
)$ satisfying $\delta \Omega =0$. When $ \Omega=0$, one obtains a
differential Gerstenhaber algebra on $\oplus _k\gm (\wedge ^k A)$,
which is known to be equivalent to the  usual definition  of Lie
bialgebroids \cite{Kosmann,Xu:1999}.

Quasi-Lie bialgebroids arise naturally in the study of the
so-called {\em twisted Poisson structures}. Recall \cite{CX,SW}
that a twisted Poisson structure consists of  a pair $(\pi,\phi)$,
where $\pi$ is a bivector field on $M$ and $\phi$ is a closed $3$-form,
which  satisfies the equation
\begin{equation}\label{phipoisson}
\frac12 [\pi, \pi]= (\wedge^3\pi^\sharp)(\phi),
\end{equation}
where $\pi^\sharp$ is the vector bundle homomorphism $T^*M\to TM$
induced by $\pi$ (i.e., $\pi^\sharp
(x)(\sigma):=\pi(x)(\sigma,\cdot )$, with $x\in M$, $\sigma\in
T^*_xM$). As explained in \cite{SW}, a twisted Poisson structure
induces a Lie algebroid structure on $T^*M$ with anchor map
$\pi^\sharp$ and Lie bracket of sections $\sigma$ and $\tau$
defined by
\begin{equation}\label{tLie}
[\sigma, \tau] := {\mathcal L}_{\pi^\sharp (\sigma)}\tau-{\mathcal
L}_{\pi^\sharp (\tau)}\sigma -d\pi(\sigma,\tau)+\phi(\pi^\sharp
(\sigma),\pi^\sharp (\tau),\cdot ).
\end{equation}
We will denote this Lie algebroid by  ${T^*M_{(\pi,\phi)}}$.
Sections of its exterior algebra are ordinary differential forms.
There is  a derivation $\delta: \Omega^\bullet
(M)\to\Omega^{\bullet +1}(M)$ deforming the de~Rham differential
$d$ by $\phi$, which is defined as follows. For $f\in
C^\infty(M)$,  $\delta f=d f$, and
$\delta\sigma = d\sigma -\pi^\sharp (\sigma )\per \phi$,
if $\sigma\in\Omega^1(M)$. It turns out that
$\delta[\sigma, \tau]=[{\delta\sigma}, \tau]+[\sigma,
{\delta\tau}]$, $\forall\sigma,\tau\in\Omega^1(M)$,
and that $\delta^2=[\phi, \cdot ]$. Thus, one obtains a quasi-Lie
bialgebroid $({T^*M_{(\pi,\phi)}}, \delta )$.

As it is known \cite{MackenzieX:2000}, a Lie bialgebroid
integrates to a Poisson groupoid. It is natural to expect that a
quasi-Lie bialgebroid should integrate to a quasi-Poisson
groupoid. For  the quasi-Lie bialgebroid corresponding to a
twisted Poisson manifold, this is true as shown by Cattaneo-Xu
\cite{CX} (it was proved in \cite{CX} that it integrates to  a
twisted symplectic groupoid, which is indeed a special example of
quasi-Poisson groupoids). However, the proof in \cite{CX} relies on
a twisted version of Poisson sigma-model, which does not apply for
a general quasi-Lie bialgebroid. The  method used by Mackenzie-Xu
\cite{MackenzieX:2000} to integrate a Lie bialgebroid, on the
other hand, does not admit a  generalization to the quasi case
either. Therefore, we must seek some new method to tackle this
integration problem.

For this purpose, we introduce the notion of multi-differentials
on a Lie algebroid. By a $k$-differential on a Lie algebroid $A$,
we mean  a linear operator $\delta :\gm (\wedge ^\bullet A)\to \gm
(\wedge ^{\bullet +k-1}A)$ satisfying
\begin{equation}
\begin{array}{l}
\delta (P\wedge Q)=(\delta P)\wedge Q+(-1)^{p(k-1)}P\wedge \delta
Q, \\[7pt] \delta \lcf P, Q\rcf =\lcf \delta P, Q\rcf +(-1)^{(p-1)(k-1)}
\lcf P, \delta Q\rcf ,
\end{array}
\end{equation}
for all $P\in \gm (\wedge ^pA)$ and $Q\in\gm (\wedge ^q A)$ (See
Definition \ref{k-differential} for an equivalent definition). The
space of all multi-differentials  $\cala = \oplus_k  \cala_k$
becomes a graded Lie algebra under the graded commutator.  Thus,
  a  quasi-Lie bialgebroid \cite{Roy} exactly corresponds to a
2-differential  whose square is a coboundary, i.e., $\delta :\gm
(\wedge ^\bullet A)\to \gm (\wedge ^{\bullet +1}A)$ such that
$\delta^2  =\lcf \Omega , \cdot \rcf $ for some $\Omega \in \gm
(\wedge^3 A )$ satisfying $\delta \Omega =0$.

On the other hand, using the generalized coisotropic calculus of
Weinstein \cite{Weinstein:1988}, we introduce the notion  of {\em
multiplicative multivector fields} on a Lie groupoid. The space of
multiplicative $k$-vector fields on $\gm$ is denoted by $\mathfrak
{X} _{mult}^k(\gm )$. One  proves that $\oplus_k {\mathfrak
{X}}_{mult}^k(\gm )$ is closed under the Schouten bracket, and
therefore is a Gerstenhaber subalgebra of $\oplus_k
{\mathfrak{X}}^k (\gm )$.

The main theorem is the following

\begin{quote}
{\bf Universal lifting theorem} Assume that $\gm \toto M$ is an
$\alpha$-simply connected and $\alpha$-connected Lie groupoid with
Lie algebroid $A$. Then  $\oplus_k  \cala_k$  is isomorphic to
$\oplus_k {\mathfrak{X}}^k_{mult}(\gm)$ as graded Lie algebras.
\end{quote}

As an immediate consequence, one concludes that there is a
bijection between quasi-Lie bialgebroids $(A, \delta, \Omega)$ and
quasi-Poisson groupoids $(\gm, \pi, \Omega )$, where
 $\gm \toto M$ is an $\alpha$-simply connected and
$\alpha$-connected Lie groupoid integrating the Lie algebroid $A$.
In particular, when $\Omega=0$, one obtains a simpler proof of the
Lie bialgebroid  integration theorem  of Mackenzie-Xu
\cite{MackenzieX:2000}. When $(A, \delta, \Omega)$ is the
quasi-Lie bialgebroid corresponding to a twisted Poisson manifold,
one recovers the integration theorem  of Cattaneo-Xu \cite{CX}.

%To prove this theorem, one direction is pretty simple. Given a
%multiplicative $k$-vector field $\Pi\in \mathfrak X
%^k_{mult}(\gm)$, for any $f \in C^{\infty} (M)$ and $X\in \gm (A)$,
%one proves that $[\Pi , \alp ^\ast f ] $ and $[\Pi, \Vec{X} ]$ are
%right invariant, where $\Vec{X}$ denotes the right invariant
%vector field on $\gm$ corresponding to $X\in \gm (A)$. Therefore,
%there exists $\delta_\Pi f \in \gm (\wedge^{k-1}A)$ and
%$\delta_\Pi X \in \gm (\wedge^{k}A)$ such that $\forall f\in
%C^\infty (M)$ and $X\in \gm (A)$,
%\[
%\Vec{\delta _\Pi f} =[\Pi , \alp ^\ast f ] , \ \ \Vec{\delta _\Pi
%X}=[\Pi, \Vec{X} ].
%\]
To prove this theorem, one direction is pretty simple. Given a
multiplicative $k$-vector field $\Pi\in \mathfrak X
^k_{mult}(\gm)$, for any $P\in \gm (\wedge ^p A)$, one proves that
$[\Pi, \Vec{P} ]$ is right invariant, where $\Vec{P}$ denotes the
right invariant $p$-vector field on $\gm$ corresponding to $P\in
\gm (\wedge ^p A)$. Therefore, there exists $\delta_\Pi P \in \gm
(\wedge^{k+p-1}A)$ such that for any $P\in \gm (\wedge ^p A)$,
\[
\Vec{\delta _\Pi P}=[\Pi, \Vec{P} ].
\]
Thus, $\delta_\Pi$ is indeed a $k$-differential. When $\gm$ is a
Lie group, the above construction is called the {\em inner
derivative} \cite{Lu:1990}. To prove the other direction, we
realize  the groupoid  $\gm$ as the moduli space of the space of
all $A$-pathes $P(A)$ moduli the gauge transformations. Such a
characterization was first obtained by Cattaneo-Felder motivated
by the Poisson sigma model when the Lie algebroid is the cotangent
Lie algebroid associated to a Poisson manifold \cite{CF}. The
general case was due to Crainic-Fernandes \cite{Crainic}.
Heuristically, the idea can be described as follows. Let $\delta$
be  a $k$-differential. Then $\delta$ naturally induces a
$k$-vector field $\pi_\delta$ on $A$, which is linear along the
fibers. It in turns gives rise to a $k$-vector field
$\tilde{\pi}_{\delta}$ on the path space $\tilde{P}(A)$. We then
prove that $\tilde{\pi}_{\delta}$ induces a $k$-vector field on
the moduli space of the space of all $A$-pathes $P(A)$.

The  second part of the paper is devoted to the  study of basic
properties of quasi-Poisson groupoids. Similar to Poisson
groupoids,  a quasi-Poisson groupoid $\gm$ also defines a momentum
theory via the so-called Hamiltonian $\gm$-spaces. Important
properties of such spaces are also studied.

A fundamental example, which is also a driving force for our
study, is the quasi-Poisson  groupoid induced by a Manin
quasi-triple $(\mathfrak d, \mathfrak g, \mathfrak h)$. Given such
a quasi-triple $(\mathfrak d, \mathfrak g, \mathfrak h)$, there  is
an associated  quasi-Lie bialgebra $(\mathfrak g, \delta, \phi )$,
where $\delta: \wedge ^\bullet \mathfrak g \to   \wedge ^{\bullet
+1}\mathfrak g$ is a derivation of the Gerstenhaber algebra
$\wedge \mathfrak g$  such that $\delta^2 =[\phi , \cdot ]$. It is
simple to see that  $\delta$ extends to a 2-differential of the
transformation Lie algebroid $\mathfrak g\times D/G\to D/G$, where $D$
and $G$ are  connected and simply connected Lie groups with Lie
algebras $\mathfrak d$ and $ \mathfrak g$ respectively, and
$\mathfrak g$ acts on $D/G$ as the infinitesimal action of the
left $G$-multiplication on $D/G$.  We explicitly  describe the
corresponding quasi-Poisson groupoid $G\times D/G \toto D/G$, and
prove that in this case the corresponding Hamiltonian $\gm$-spaces
are equivalent to the quasi-Poisson spaces with $D/G$-momentum
maps in the sense of Alekseev-Kosmann-Schwarzbach \cite{AK}.
However, our approach does not require $\mathfrak h$ to be
admissible.

When $(\mathfrak d, \mathfrak g, \mathfrak h)$ is a Manin triple
and $G$ is a complete Poisson group, we recover the Poisson
groupoid $G\times G^*\toto G^*$  \cite{LW} (which is symplectic in
this case). When $G$ is not necessary complete, the Poisson
groupoid $G\times D/G \toto D/G$, which is in fact a symplectic
groupoid integrating the Poisson structure on $D/G$, can be
considered as a replacement of Lu-Weinstein's symplectic groupoid
$G\times G^*\toto G^*$.

A particularly   interesting case is the Manin quasi-triple
$(\mathfrak g\oplus \mathfrak g, \Delta ({\mathfrak g}), \half
\Delta_{-} ({\mathfrak g}))$ corresponding to a Lie algebra
$\mathfrak g$ equipped with a non-degenerate symmetric pairing. In
this case, one obtains a quasi-Poisson groupoid structure on
$G\times G\toto G$, where $G$ acts on $G$ by conjugation. The
discussion on this part is done in Section 4.

Finally, some remarks are in order about  notations. For a Lie
groupoid $\gm\gpd M$, by  $\alp ,\bet :\gm\to M$, we denote  the
source and target maps. Two elements $g,h\in \gm$ are composable
if  $\bet (g)=\alp (h)$. We denote by $\gm^{(2)}\subset \gm\times
\gm$ the subset of composable pairs in $\gm \times \gm$. By
$i:\gm \to \gm$, $g\mapsto i(g)=g^{-1}$, we denote   the
inversion, and  $\epsilon :M\to \gm$, $x\mapsto \epsilon
(x)=\tilde{x}$, the unit map.

The cotangent Lie groupoid is denoted by $T^\ast \gm\gpd A^\ast$,
$A^\ast$ being the dual of the Lie algebroid $A$, and the
structural functions of this groupoid are
$\tilde{\alp},\tilde{\bet}$ for the projections, $\cdot$ for the
multiplication, $\tilde{\mbox{{\i}}}$ for the inversion and
$\tilde{\epsilon}$ for the unit map.

\medskip
{\bf Acknowledgments.}
We  would like to thank several institutions
for their hospitality while work on this project was being done:
University of La Laguna (Iglesias Ponte), Penn State University
(Laurent-Gengoux), and Universit\'e Lyon 1 and Universit\'e Pierre et Marie
Curie (Xu).  Iglesias wishes to thank Fulbright/Spanish ministry
of Education and Culture for a postdoctoral fellowship.
We would like to thank Alberto Cattaneo,  Marius Crainic, Rui Loja Fernandes, 
Jiang-hua Lu, and Alan Weinstein for  useful discussions.

\section{Multiplicative $k$-vector fields on Lie groupoids}

\subsection{Coisotropic submanifolds}

In this section, we generalize Weinstein's coisotropic calculus
\cite{Weinstein:1988} to multivector fields.

\begin{defn}
Let $V$ be a vector space and  $\Pi \in \wedge ^k V$. We say that
a subspace $W$ of $V$ is  {\em coisotropic} with respect to $\Pi$
if
\[
\Pi  (\xi ^1,\ldots ,\xi ^k)=0
\]
for all $\xi ^1,\ldots ,\xi ^k\in W^{\smalcirc}$, where
$W^{\smalcirc}$ is the annihilator of $W$, that is
$W^{\smalcirc}=\{ \xi \in V^\ast \, | \, \xi _{|W}=0\}$.
\end{defn}
A slight modification of coisotropic calculus for bivector fields
developed in \cite{Weinstein:1988} will allow us to show several
properties in the sequel. We give now the following important
lemma.
\begin{lem}\label{coisotropic-lemma}
Let $V_1,V_2$ be vector spaces and $\Pi _i\in \wedge^k V_i$ for
$i=1, 2$.  If $R\subseteq V_1\times V_2$ is coisotropic with
respect to $\Pi _1\oplus \Pi _2$ and $C\subseteq V_2$ is
coisotropic with respect to $\Pi _2$ then
\[ R(C) =\{ u\in V_1\, |\, \exists\, v \in C \mbox{ with }
(u,v)\in R  \} \] is coisotropic with respect to $\Pi _1$.
\end{lem}
\begin{pf} Let $R^\ast\subseteq V^\ast _1\times V_2^\ast$ be the subspace
given by $$R^\ast = \{ (\xi ^1,\xi ^2)\in V^\ast_1\times
V^\ast_2\, |\, \langle \xi ^1, v_1\rangle =\langle \xi ^2,
v_2\rangle, \,\, \forall\, (v_1,v_2)\in R\}.$$ Then we have
$R(C)^{\smalcirc}=R^\ast (C^{\smalcirc} )$.

Now, let $\xi ^1,\ldots ,\xi ^k \in R(C)^{\smalcirc}=R^\ast
(C^{\smalcirc} )$. Then, there exists $\varphi ^1,\ldots ,\varphi
^k \in C^{\smalcirc} $ such that $(\xi ^i,\varphi ^i)\in R^\ast$,
i.e., $(\xi ^i,-\varphi ^i)\in R^{\smalcirc}$. Thus,
\[ \begin{array}{lll}
0&=& (\Pi _1\oplus \Pi _2)((\xi ^1,-\varphi ^1),\ldots ,(\xi
^k,-\varphi ^k))=\Pi _1 (\xi ^1,\ldots ,\xi ^k )+(-1)^k \Pi _2
(\varphi ^1,\ldots
,\varphi ^k)\\[5pt] &=&\Pi _1 (\xi ^1,\ldots ,\xi ^k ).
\end{array}
\]
That is, $R(C)$ is coisotropic. \end{pf} A generalization of the
notion of coisotropy to manifolds is the following.
\begin{defn}
Let $M$ be an arbitrary manifold and $\Pi \in \mathfrak X ^k(M)$.
A submanifold $S$ of $M$ is said to be {\em coisotropic} with
respect to $\Pi$ if $T_x S$ is coisotropic with respect to
$\Pi(x)$ for all $x\in S$.
\end{defn}

\begin{numrmk}\label{coisotrop-charact}
{\rm
%By $N^* S$ we denote the conormal bundle of $S$, i.e.,
%$N^\ast _x S=(T_xS)^{\smalcirc}$ for any $x\in S$.
 It is simple  to see that
$S$ is coisotropic with respect to the  multivector field  $\Pi$
if and only if $\Pi (df^1,\ldots ,df^k)_{|S}=0$ for any
$f^1,\ldots ,f^k \in C^\infty (M)$ such that $f^i_{|S}=0$ (see
\cite{Weinstein:1988} for the case of bivector fields).}
\end{numrmk}

\begin{prop}\label{Schouten}
If $S$ is coisotropic with respect to the multivector fields $\Pi$
and $\Pi '$, so is with respect to the multivector field $[\Pi
,\Pi ']$.
\end{prop}
\begin{pf}
If $A$ is any subset of $\{ 1,2,\ldots ,(k+k'-1)\}$, let $A'$
denote its complement and $|A|$ the number of elements in $A$. If
$|A|=l$ and the elements in $A$ are $\{ i_1,\ldots ,i_l\}$ in
increasing order, let us write $f_A$ for the ordered $k$-tuple
$(f^{i_1},\ldots ,f^{i_l})$. Furthermore, we write $\varepsilon
_A$ for the sign of the permutation which rearranges the elements
of the ordered $(k+k'-1)$-tuple (A',A), in the original order.
Then, the Schouten bracket of $\Pi \in \mathfrak X^k(M)$ and $\Pi
'\in \mathfrak X^{k'}(M)$ is given \cite{BhaskaraV} by
\begin{equation}\label{eq:Schoutenbracket}
\begin{array}{rl}
[ \Pi  ,\Pi ' ] (df^1,\ldots ,df^{k+k'-1})=(-1)^{k+1}\Big (
&\displaystyle \sum_{|A|=k'} \varepsilon _A \Pi (d(\Pi '
(f_A)),df_{A'})\\& +(-1)^{kk'} \displaystyle \sum_{|B|=k}
\varepsilon _B \Pi ' (d(\Pi (f_B)),df_{B'})\Big ) .
\end{array}
\end{equation}

Our result  thus follows from the characterization of a
 coisotropic submanifold
as described in Remark \ref{coisotrop-charact}.
\end{pf}

\subsection{Definition and examples}
In this section, we will introduce the notion of  multiplicative
multivector fields and  show that it generalizes several
concepts which have previously appeared in the literature.

Throughout this section, we fix a Lie groupoid $\gm\gpd M$ and
denote its Lie algebroid by $A$.  Moreover, we  denote by
$\Lambda$ the graph of the groupoid multiplication, that is,
\[
\Lambda =\{ (g,h,gh)\, | \, \bet (g)=\alp (h) \}.
\]

\begin{defn}\label{k-multiplicative}
Let $\gm\gpd M$ be a Lie groupoid and $\Pi \in \mathfrak X ^k
(\gm)$
 a $k$-vector field on $\gm$. We say that $\Pi$ is {\em
multiplicative}, denoted as $\Pi \in \mathfrak X ^k _{mult}
(\gm)$, if $\Lambda$ is coisotropic with respect to $\Pi \oplus
\Pi \oplus (-1)^{k+1}\Pi$.
\end{defn}

An interesting characterization of multiplicative multivector
fields  is the following:
\begin{prop}\label{charact-multiplicative}
Let $\gm\gpd M$ be a Lie groupoid and $\Pi \in \mathfrak X ^k
(\gm)$
 a $k$-vector field on $\gm$. Then, the following are equivalent:
\begin{itemize}
\item[{\it i)}] $\Pi$ is multiplicative;
\item[{\it ii)}] for any $\mu _g^i,\nu _h^i\in {T^\ast \gm}$, such that
$\tilde{\bet}(\mu ^i_g)=\tilde{\alp}(\nu ^i_h)$
\begin{equation}\label{multiplicative}
\Pi (gh)(\mu _g^1\cdot \nu _h^1,\ldots ,\mu _g^k\cdot \nu _h^k)=
\Pi (g)(\mu _g^1,\ldots ,\mu _g^k)+ \Pi (h)(\nu _h^1,\ldots ,\nu
_h^k);
\end{equation}
\item[{\it iii)}]
the linear skew-symmetric function $F_\Pi$
 on $T^\ast \gm\times _\gm \stackrel{(k)}{\ldots }\times _\gm T^\ast \gm$
induced by $\Pi$,
$$F_\Pi (\mu^1, \ldots , \mu^k)=\Pi (\mu^1, \ldots , \mu^k),$$
  is a  1-cocycle with respect to the  Lie groupoid
\begin{equation}\label{k-groupoid}
T^\ast \gm\times _\gm \stackrel{(k)}{\ldots }\times _\gm T^\ast
\gm \gpd A^\ast \times _M \stackrel{(k)}{\ldots }\times _M A^\ast
.
\end{equation}
Here $T^\ast \gm\times _\gm\stackrel{(k)}{\ldots}\times _\gm
T^\ast \gm = \{ (\mu ^1,\ldots ,\mu ^k) \in T^\ast \gm\times
\ldots \times T^\ast \gm | \tau (\mu ^1)=\ldots =\tau (\mu ^k)\}$,
$\tau :T^\ast \gm\to \gm$ is the bundle projection, $A^\ast \times
_M \stackrel{(k)}{\ldots }\times _M A^\ast = \{ (\eta ^1,\ldots
,\eta ^k) \in A^\ast \times \ldots \times A^\ast | p _\ast(\eta
^1)=\ldots =p_\ast (\eta ^k)\}$, $p_\ast :A^\ast \to M$, and
 $T^\ast \gm\times _\gm \stackrel{(k)}{\ldots }
\times _\gm T^\ast \gm\gpd A^\ast \times _M
\stackrel{(k)}{\ldots }\times _M A^\ast$ is considered as
 a Lie subgroupoid of the
direct product  groupoid $(T^\ast \gm)^k\gpd (A^\ast )^k$.
\end{itemize}
\end{prop}
\begin{pf}
It follows from the fact that $(\mu _g,\nu _h,\gamma _{gh})\in
N^\ast _{(g,h,gh)}\Lambda$ if and only if $\gamma _{gh}=-\mu
_g\cdot \nu _h$, for $\mu _g,\nu _h\in T^\ast \gm$.
\end{pf}
\begin{numrmk}\label{multiplicative-bivectors}
{\rm From Proposition \ref{charact-multiplicative} one can deduce
that a bivector field $\Pi$ on a Lie groupoid $\Gamma$ is
multiplicative if and only if $\Pi ^\sharp :T^\ast \Gamma \to
T\Gamma$ is a Lie groupoid morphism, i.e.,
\[
\Pi ^\sharp (\mu _g\cdot \nu _h)=\Pi ^\sharp (\mu _g)\cdot \Pi
^\sharp (\nu _h),
\]
for $\mu _g,\nu_h\in T^\ast \Gamma$ such that $\tilde{\bet}(\mu
_g)=\tilde{\alp}(\nu _h)$.}
\end{numrmk}

A direct consequence of Proposition \ref{Schouten} and Definition
\ref{k-multiplicative} is the following.
\begin{prop}
The Schouten bracket of multiplicative multivector fields is still
multiplicative. That is, $(\oplus _k \mathfrak X _{mult}^k
(\gm),[\cdot ,\cdot ])$
 is a graded Lie subalgebra of $(\oplus _k \mathfrak X ^k
(\gm),[\cdot ,\cdot ])$.
\end{prop}

Below we  list some well-known examples.
\begin{numex}
{\rm Let $G$ be a Lie group and $\Pi \in \mathfrak X^k(G)$ a
$k$-vector field. From $\tilde{\bet}(\mu ^i_g)=\tilde{\alp}(\nu
^i_h)$ we deduce that
\[
\mu ^i_g=(L_{g^{-1}})^\ast (R_h)^\ast \nu ^i_h .
\]
Moreover, we have that
\[
\mu _g^i\cdot \nu _h^i=(L_{g^{-1}})^\ast \nu ^i_h.
\]
Therefore, we see that Eq. (\ref{multiplicative}) is equivalent to
\[ \Big ((L_{g^{-1}})_\ast \Pi (gh) -(L_{g^{-1}})_\ast (R_h)_\ast
\Pi (g)-\Pi (h)\Big ) (\nu ^1_h,\ldots ,\nu ^k_h)=0.
\]
That is, $\Pi (gh) =(R_h)_\ast \Pi (g)+(L_g)_\ast\Pi (h)$. The
converse is obvious. Therefore, our definition of multiplicative
$k$-vector fields is indeed a generalization of the usual notion
for Lie groups (see \cite{LuW:1990}).}
\end{numex}
\begin{numex}
{\rm We say that a function $\sigma$ on $\gm$ is multiplicative if
$\sigma (gh)=\sigma (g)+\sigma (h)$ for $(g,h)\in \gm^{(2)}$.
Therefore, multiplicative functions are multiplicative $0$-vector
fields.}
\end{numex}
\begin{numex}
{\rm  A vector field $X\in \mathfrak X(\gm)$ is said to be
multiplicative if it is a Lie groupoid morphism $X:\gm\to T\gm$
from a Lie groupoid $\gm\gpd M$ to its corresponding tangent Lie
groupoid $T\gm\gpd TM$ (see \cite{MackenzieX:1998}). Therefore, we
have $X(gh)=X(g)\cdot X(h)$ for $\bet (g)=\alp (h)$. Using this
fact and that
\[
(\mu _g\cdot \nu _h)(u_g\cdot v_h)=\mu _g (u_g)+\nu _h(v_h),\mbox{
for all }(u_g ,v_h)\in T\gm^{(2)},
\]
we deduce that $X$ is a multiplicative $1$-vector field in the
sense of Definition \ref{k-multiplicative}.}
\end{numex}
\begin{numex}
{\rm  From Definition \ref{k-multiplicative}, we see that the
Poisson tensor on a Poisson groupoid is a multiplicative bivector
field \cite{Weinstein:1988}.}
\end{numex}
\begin{numex}
{\rm  Given a Lie groupoid $\gm$ with Lie algebroid $A$, if  $P
\in \gm (\wedge ^k A)$, then $\Vec{P}-\ceV{P}$ is a multiplicative
$k$-vector field on $\gm\gpd M$. }
\end{numex}

\subsection{Multiplicative and affine multivector fields}
Let $\Omega$ be the submanifold  of $\gm\times \gm\times \gm\times \gm$
consisting of elements $(l,h,g,w)$ such that $w=hl^{-1}g$.
$\Omega$ is called the {\em affinoid diagram} corresponding to the
groupoid $\gm$ by Weinstein \cite{Weinstein:1990}. A
characterization of multiplicative $k$-vector fields in terms of
the affinoid diagram is the following:
\begin{prop}\label{charac}
Let $\Pi$ be a $k$-vector field on a Lie groupoid. Then $\Pi$ is
multiplicative if and only if $\Omega$ is coisotropic with respect
to $\Pi\oplus (-1)^{k+1}\Pi \oplus (-1)^{k+1}\Pi \oplus \Pi$ and
$M$ is coisotropic with respect to $\Pi$.
\end{prop}
\begin{pf} Suppose that $\Pi$ is a multiplicative $k$-vector field on
$\gm$. Using Eq. (\ref{multiplicative}) and the identity
$\tilde{\epsilon}(\phi _m)\cdot \tilde{\epsilon}(\phi
_m)=\tilde{\epsilon}(\phi _m)$ for all $\phi _m\in A^\ast _m$, we
have
\[
\Pi(\epsilon (m))( \tilde{\epsilon}(\phi ^1_m),\ldots
,\tilde{\epsilon}(\phi ^k _m))=2\Pi(\epsilon (m))(
\tilde{\epsilon}(\phi ^1_m),\ldots ,\tilde{\epsilon}(\phi ^k _m)),
\]
$\forall \phi _m^1,\ldots ,\phi _m^k\in A^\ast _m$, $m\in M$.
Since we have $\tilde{\epsilon}(A^\ast )=N^\ast M$, we deduce that
$M$ is coisotropic with respect to $\Pi$.

Now we  proceed as in \cite{Weinstein:1990}. In the product
$\gm\times\bar{\gm}\times \bar{\gm}\times \gm\times \gm\times
\bar{\gm}$ (where $\bar{\gm}$ denotes $\gm$ endowed with the
$k$-vector field $(-1)^{k+1}\Pi$), we consider the coisotropic
submanifold $R=\{ (g,h,l,x,y,z)\, | \, gy=l \mbox{ and } hz=x\}$.
On the other hand, it is  simple  to see
 that the diagonal  $\Delta  \subset \gm\times \bar{\gm}$ is
a coisotropic submanifold.
Therefore, using Lemma \ref{coisotropic-lemma}, we have that
$R(\Delta)$ is coisotropic submanifold of $\gm\times \bar{\gm}\times
\bar{\gm}\times \gm$. It is simple to see that $R(\Delta)=\Omega$. Our
result thus follows.

Conversely, let $\Pi$ be a $k$-vector field on a Lie groupoid such
that $\Omega$ is coisotropic with respect to $\Pi\oplus
(-1)^{k+1}\Pi \oplus (-1)^{k+1}\Pi \oplus \Pi$ and $M$ is
coisotropic with respect to $\Pi$. Applying Lemma
\ref{coisotropic-lemma} to $R=\Omega$ and $C=M$, we obtain that
$\Lambda$ is coisotropic with respect to $\Pi \oplus \Pi \oplus
(-1)^{k+1}\Pi$. That is, $\Pi$ is multiplicative. \end{pf}

Next we recall the definition of an affine multivector field on a
Lie groupoid and  show that any multiplicative multivector field
is affine.
\begin{defn}
A multivector field $\Pi$ on $\gm$ is {\em affine} if for any $g,
h\in \gm$ such that $\bet (g)= \alp (h)=m$ and any bisections
${\cal X}, {\cal Y}$ through the points $g, h$, we have
\begin{equation}\label{affine}
\Pi (gh) = (R_{\cal Y })_\ast \Pi (g) + (L_{\cal X })_\ast \Pi (h)
- (R_{\cal Y }\smalcirc L_{\cal X })_\ast \Pi (\epsilon (m)).
\end{equation}
\end{defn}

A useful characterization of affine multivector fields is the
following

\begin{prop}\cite{MackenzieX:2000}\label{charact-affine}
Let $\Pi$ be a $k$-vector field on a Lie groupoid $\gm\gpd M$.
Then $\Pi$ is affine if and only if $[\Vec{X},\Pi ]$ is
right-invariant for all $X\in \gm (A)$.
\end{prop}

\begin{prop}\label{mult->affine}
If $\Pi$ is a multiplicative $k$-vector field on a Lie groupoid
$\gm$, then $\Pi$ is affine.
\end{prop}
\begin{pf} If $\Pi$ is multiplicative, according to  Proposition
\ref{charac}, we know that $\Omega$ is coisotropic with respect to
$\Pi\oplus (-1)^{k+1}\Pi \oplus (-1)^{k+1}\Pi \oplus \Pi$,
 and $M$ is coisotropic with respect to $\Pi$. For any $\mu \in T_{gh}\gm
$, it follows from Lemma 2.6 in \cite{Xu:1995} that $(-\mu,
L^{*}_{\cal X}\mu, R^{*}_{\cal Y}\mu, - L^{*}_{\cal X} R^{*}_{\cal
Y}\mu)$ is  conormal  to $\Omega$. Therefore, for any $\mu
^1,\ldots , \mu ^k \in T_{gh}\gm$, we have
\[
\begin{array}{rll}
-\Pi (gh) (\mu ^1, \ldots ,\mu ^k) +\Pi (h)( L^{*}_{\cal X }\mu
^1,\ldots, L^{*}_{\cal X}\mu ^k ) +\Pi (g) ( R^{*}_{\cal Y
}\mu ^1,\ldots, R^{*}_{\cal Y}\mu ^k ) && \\[5pt]-\Pi (\epsilon(m)) (
L^{*}_{\cal X} R^{*}_{\cal Y}\mu ^1,\ldots , L^{*}_{\cal X}
R^{*}_{\cal Y}\mu ^k ) &=& 0.
\end{array}
\]
Thus   Eq. (\ref{affine}) follows immediately. \end{pf}

Similar to the case of multiplicative bivector fields, we can
 give another useful characterization of
multiplicative $k$-vector fields.

\begin{them}\label{properties}
Let $\gm\gpd M$ be a Lie groupoid and $\Pi \in \mathfrak X ^k
(\gm)$  a $k$-vector field on $\gm$. Then $\Pi$ is
multiplicative if and only if the following conditions hold:
\begin{itemize}
\item[{\it i)}] $\Pi$ is affine, i.e., Eq.  (\ref{affine}) holds;
\item[{\it ii)}] $M$ is a coisotropic submanifold of $\gm$;
\item[{\it iii)}] both $\alp _\ast \Pi (g)$ and $\bet _\ast \Pi (g)$
only depend on $\alp (g)$ and $\bet (g)$, respectively;
\item[{\it iv)}]  $\forall \eta ^1,\eta^2\in \Omega ^1(M)$,
we have  $(\alp ^\ast \eta ^1\wedge \bet ^\ast \eta ^2)\per \Pi
=0$;
\item[{\it v)}]  $\forall \theta \in \Omega ^p (M)$,
 $1\leq p<k$, then $(\bet ^\ast\theta )\per \Pi$ is a
 left-invariant ($k-p$)-vector field on $\gm$.
\end{itemize}
\end{them}
\begin{pf} Let $\Pi$ be a multiplicative $k$-vector field on $\gm$. From
Propositions \ref{charac} and \ref{mult->affine}, we obtain {\it
i)} and {\it ii)}.

Next, since
\begin{equation}\label{alpha}
0_g\cdot (\bet ^\ast \eta )_{h}=(\bet ^\ast \eta )_{gh}
\end{equation}
for any $\eta \in \Omega ^1(M)$, from Eq. (\ref{multiplicative}),
 it follows that
\[
\begin{array}{lll}
\Pi( (\bet ^\ast \eta ^1)_{gh} ,\ldots , (\bet ^\ast \eta
^k)_{gh}) &=&\Pi _{}( 0_g\cdot (\bet ^\ast \eta ^1)_{h} ,\ldots
,0_g\cdot (\bet ^\ast \eta ^k)_{h})
\\&=&\Pi ( (\bet ^\ast \eta ^1)_{h} ,\ldots , (\bet^\ast \eta ^k)_{h})
\end{array}
\]
$\forall \eta ^1,\ldots ,\eta ^k \in \Omega ^1(M)$. Hence $\bet
_\ast \Pi _{g}$ only depends on $\bet (g)$. Similarly, from the
equation
\begin{equation}\label{beta}
(\alp ^\ast \eta )_{g}\cdot 0_h=(\alp ^\ast \eta )_{gh}
\end{equation}
$\forall \eta \in \Omega ^1 (M)$, we also deduce that $\alp _\ast
\Pi _{g}$ only depends on $\alp (g)$. Hence, {\it iii)} holds.

Moreover, from Eqs. (\ref{multiplicative}), (\ref{alpha}) and
(\ref{beta}), we have {\it iv)}.

Finally, let us prove {\it v)}. Let $\theta \in \Omega ^p(M)$ be a
$p$-form on $M$, $1\leq p<k$. Then, using Eqs.
(\ref{multiplicative}), (\ref{alpha}) and (\ref{beta}),
we see that  $(\bet^\ast \theta)\per\Pi $ is  tangent to  $\alp$-fibers.
Therefore, it  suffices to prove that
$$(L_{\cal X})_\ast ((\bet ^\ast \theta)\per \Pi )
_{(h)})=((\bet ^\ast \theta)\per \Pi )_{(gh)}, \ \ \forall
(g,h,gh)\in \Lambda, $$
 where $\cal X$ is an arbitrary bisection through $g$.
According to  \cite{Xu:1995},  for any $\mu \in
T^\ast_{gh}\gm$,
 there exists $\nu \in T^\ast _g\gm$, which is  characterized by
the equation:
\[ \langle \nu , v_g \rangle = \langle \mu, (R_h)_\ast
(v_g - (L_{\cal X} )_\ast \bet _\ast v_g)\rangle , \forall \, v_g
\in T_g\gm ,
\]
and $\nu \cdot (L_{\cal X})^\ast \mu =\mu$. Using this fact and
Eq. (\ref{beta}), it follows that
\[
((\bet ^\ast \theta)\per \Pi )_{(h)}((L_{\cal X})^\ast \mu
^1,\ldots ,(L_{\cal X})^\ast \mu ^{p-k})=((\bet ^\ast \theta)\per
\Pi )_{(gh)} (\mu ^1,\ldots , \mu ^{p-k})
\]
$\forall \mu^1,\ldots ,\mu^{p-k}\in \Omega^1(\gm)$. Thus, {\it v)}
follows.

To prove the converse, we first note that the following three
types of vectors span the whole conormal space of $\Omega$ at a
point $(g,h,l,w)$: $(-\mu, L^{*}_{\cal X}\mu, R^{*}_{\cal Y}\mu, -
L^{*}_{\cal X} R^{*}_{\cal Y}\mu)$ for any
 $\mu \in T^\ast _{l}\gm$,
$(-\bet ^\ast \eta ,\bet ^\ast \eta ,0,0)$ for any $\eta \in
T^\ast _{\alp (k)}\gm $, and $(-\alp ^\ast \zeta ,0,\alp ^\ast
\zeta ,0)$ for any $\zeta \in T^\ast _{\bet (z)}\gm $, where
${\cal X}$ and ${\cal Y}$ are any bisections through the points
$g$ and $h$, respectively (see \cite{Xu:1995}). From {\em i)},
{\em iii)}, {\em iv)} and {\em v)} we deduce that $\Omega$ is
coisotropic with respect to $\Pi\oplus (-1)^{k+1}\Pi \oplus
(-1)^{k+1}\Pi \oplus \Pi$ (see Theorem 2.8 in \cite{Xu:1995}).
Using this fact, {\em ii)} and Proposition \ref{charac}, the conclusion
follows.
\end{pf}

In a similar fashion that we prove {\it v)} in Theorem
\ref{properties}, we have  the following
\begin{cor}\label{right-inv}
If $\Pi$ is a multiplicative $k$-vector field on $\gm$ then for
all $\theta \in \Omega ^p (M)$, $1\leq p<k$, then $(\alp
^\ast\theta )\per \Pi$ is a right-invariant ($k-p$)-vector field
on $\gm$.
\end{cor}
Finally, let us show an interesting property that generalizes the
one obtained for multiplicative bivector fields in
\cite{Weinstein:1988}.

\begin{prop}\label{base}
Let $\gm \rightrightarrows M$ be a Lie groupoid and $\Pi \in
\mathfrak X ^k (\gm)$ be a multiplicative $k$-vector field on
$\gm$. Then there exists a unique $k$-vector field $\pi$  on $M$
such that
\[
\alp _\ast \Pi = \pi ,\qquad \bet _\ast \Pi = (-1)^{k+1}\pi .
\]
\end{prop}
\begin{pf} Since $\Pi$ is multiplicative, using property {\it iii)} in
Theorem \ref{properties}, we can define a $k$-vector field $\pi $
on $M$ by  setting $\pi =\alp _\ast \Pi$. Now, let us investigate
the relation between $\Pi$, $\pi$ and the map $\bet$.

First, we  show that if $i:\gm\to \gm$ is the groupoid
inversion then

\begin{equation}\label{inversion}
i_\ast \Pi  = (-1)^{k+1} \Pi .
\end{equation}
This  is an immediate consequence of property {\it ii)} in Theorem
\ref{properties} and the fact that the inverse $(\mu _g)^{-1}$ of
$\mu _g\in T^\ast _g\gm$ is given by $(\mu _g)^{-1}= -i ^\ast (\mu
_g)$. In fact,
\[
\begin{array}{lll}
0&=&\Pi ( \tilde{\epsilon}( \tilde{\bet}( \mu ^1_g)),\ldots ,
\tilde{\epsilon}(\tilde{\bet}( \mu
^k_g))) \\[5pt]
 &=&\Pi  (\mu ^1_g\cdot
(\mu^1_g)^{-1}, \ldots ,\mu ^k_g\cdot (\mu ^k_g)^{-1})\\[5pt]
&=& \Pi (\mu ^1 _g,\ldots ,\mu ^k _g) +(-1)^k (i_\ast\Pi ) (\mu ^1
_g,\ldots ,\mu ^k _g)
\end{array}
\]
for any $\mu ^1_g,\ldots ,\mu ^k_g\in T^\ast_g \gm $. Finally,
using Eq. (\ref{inversion}) and the relation $\alp \smalcirc i
=\bet$, we conclude that $\bet _\ast \Pi =(-1)^{k+1}\pi$. \end{pf}
\begin{numex}\label{coarse-example}
{\rm From Proposition \ref{base} and property {\it iv)} in Theorem
\ref{properties}, we obtain that the map
 $(\alp ,\bet ):\gm\to M\times M$ satisfies
\[
(\alp ,\bet )_\ast\Pi =\pi \oplus (-1)^{k+1}\pi .
\]
In  particular,  when $\gm$ is a  pair
 groupoid $M\times M\gpd M$,
since $(\alp ,\bet ):M\times M\to M\times M$ is a diffeomorphism,
we  conclude that the only multiplicative $k$-vector fields are of
the form $\pi \oplus (-1)^{k+1}\pi $, $\pi \in \mathfrak X ^k(M)$
(see \cite{MackenzieX:1998} and \cite{Weinstein:1988} for the case
$k=1$ and $k=2$, respectively).}
\end{numex}

\subsection{$k$-differentials on Lie algebroids}
We now turn to the study of the Lie algebroid counterpart of
multiplicative $k$-vector fields, namely, $k$-differentials.

\begin{defn}\label{k-differential}
Let $(A,\lcf \cdot ,\cdot \rcf, \rho )$ be  a  Lie algebroid over
$M$. An {\em almost $k$-differential} is a pair of linear maps
$\delta: C^{\infty}(M) \lon \gm (\wedge^{k-1}A)$ and $\delta: \gm
(A)\lon \gm (\wedge^{k}A)$ satisfying
\begin{itemize}
\item[{\em i)}]  $\delta (fg)=g(\delta f)+f(\delta g)$,
for all $f, g\in  C^{\infty}(M)$.
\item[{\em ii)}] $\delta (fX)=(\delta f)\wedge X+f\delta X$,
 for all $f \in C^{\infty}(M)$ and $X\in \gm (A)$.
\end{itemize}
An almost $k$-differential is said to be a {\em $k$-differential}
if it  satisfies the compatibility condition
\begin{equation}\label{compatibility}
\delta \lcf X, Y\rcf =\lcf \delta X, Y\rcf +\lcf X, \delta Y\rcf,
\end{equation}
for all $X, Y\in \gm (A)$.
\end{defn}

For a given Lie algebroid $A$, it is known that
 the anchor map together with the bracket on $\gm (A)$ extends to a graded Lie
bracket on $\oplus _k\gm (\wedge^k A)$, which makes it into a
Gerstenhaber algebra  $ (\oplus _k \gm (\wedge^k A), \lcf \cdot,
\cdot \rcf, \wedge )$ \cite{Xu:1999}.

A $k$-differential $\delta$ extends naturally to sections of
$\wedge A$ as follows:

\begin{equation}
\delta (X_{1}\wedge \cdots \wedge X_s )= \sum _{i=1}^s
(-1)^{(i+1)(k+1)}X_{1}\wedge \cdots \wedge (\delta X_i)\wedge
\cdots \wedge X_s
\end{equation}
$\forall X_1,\ldots X_s\in \gm (A)$. In this way,
 we obtain a linear operator $\delta :\gm (\wedge
^\bullet A)\to \gm (\wedge ^{\bullet +k-1}A)$.
This following proposition can be directly verified.

\begin{prop}
A $k$-differential  on a  given Lie algebroid $A$ is
equivalent to a derivation of the associated
Gerstenhaber algebra  $ (\oplus _k \gm (\wedge^k A), \lcf \cdot,
\cdot \rcf, \wedge )$, i.e., a linear operator $\delta :\gm (\wedge
^\bullet A)\to \gm (\wedge ^{\bullet +k-1}A)$  satisfying
\begin{equation}\label{k-dif-prop}
\begin{array}{l}
\delta (P\wedge Q)=(\delta P)\wedge Q+(-1)^{p(k-1)}P\wedge \delta
Q, \\[7pt] \delta \lcf P, Q\rcf =\lcf \delta P, Q\rcf +(-1)^{(p-1)(k-1)}
\lcf P, \delta Q\rcf ,
\end{array}
\end{equation}
for all $P\in \gm (\wedge ^pA)$ and $Q\in\gm (\wedge ^q A)$.
\end{prop}

 As we
see below,
 $k$-differentials reduce to  various  well-known notions
in special cases.

\begin{numex}
{\rm Let $\mathfrak g$ be a Lie algebra. A $k$-differential on
$\mathfrak g$ is just a linear map $\delta :\mathfrak g\to \wedge
^k\mathfrak g$ such that it is a 1-cocycle with respect to the
adjoint representation of $\mathfrak g$ on $\wedge ^k\mathfrak g$
\cite{LuW:1990}.}
\end{numex}
\begin{numex}
{\rm Let $\delta$ be a $0$-differential, that is, $\delta f=0$ for
$f\in C^\infty (M)$, and  $\delta X\in C^\infty (M)$ for $X\in \gm
(A)$. From {\em ii)} in Definition \ref{k-differential}, we deduce
that $\delta (fX)=f\delta X$. Therefore, there exists $\phi \in
\gm (A^\ast )$ such that
\begin{equation}
\delta X =\phi (X)\mbox{ for }X\in \gm (A).
\end{equation}
Moreover, using Eq. (\ref{compatibility}), we obtain that $\phi$
is a 1-cocycle in the Lie algebroid cohomology of $A$. Thus, we
can conclude that $0$-differentials are just Lie algebroid
1-cocycles with trivial coefficients. }
\end{numex}
\begin{numex}
{\rm  Let $\delta$ be an almost $1$-differential, that is, $\delta
f\in C^\infty (M)$ for $f\in C^\infty (M)$, and  $\delta X\in \gm
(A)$ for $X\in \gm (A)$. From {\em i)} in Definition
\ref{k-differential}, we deduce that there exists $X_0 \in
\mathfrak X (M)$ such that
\begin{equation}
\delta f =X_0(f) \mbox{ for } f\in C^\infty (M).
\end{equation}
Moreover, using {\em ii)}, we obtain that
\[
\delta (fX)=f\delta (X)+X_0 (f)X,
\]
that is, $\delta$ is a covariant differential operator on $A$,
with anchor $X_0$ (see \cite{Mackenzie:LGLADG,MackenzieX:1998}).
If, moreover, $\delta$ is a 1-differential, from Eq.  (\ref{compatibility}) we
see that $\delta \in \gm (CDO(A))$ is a derivation of the bracket
on $\gm (A)$. }
\end{numex}
\begin{numex}
{\rm It is well-known that a Lie algebroid structure on a vector
bundle $A\to M$ is equivalent to an almost 2-differential $\delta
:\gm (\wedge ^\bullet A^\ast )\to \gm (\wedge ^{\bullet +1}A^\ast
)$ of square $0$ (see, for instance, \cite{Kosmann,Xu:1999}).
Thus, we see that a Lie bialgebroid (a notion first introduced in
\cite{MackenzieX:1994}) corresponds to a 2-differential of square
$0$ on a Lie algebroid $A$. }
\end{numex}
\begin{numex}
{\rm  If $P \in \gm (\wedge^k A)$, then $ad(P) =\lcf P , \cdot
\rcf$ is clearly a $k$-differential, which is called the {\em
coboundary} $k$-differential associated to $P$.}
\end{numex}

The space of almost differentials can be endowed with a graded Lie
algebra structure as shown in the following:
\begin{prop}
Let $\hat \cala_k$ denote the space of almost $k$-differentials
and $\hat \cala = \oplus_k  \hat \cala_k$. If we define
\begin{equation}\label{graded-bracket}
[\delta_1 , \delta_2 ]=\delta_1 \smalcirc \delta_2
-(-1)^{(k+1)(l+1)} \delta_2 \smalcirc \delta_1
\end{equation}
for $\delta _1\in \hat \cala _k$ and $\delta _2\in \hat \cala _l$,
then
\begin{itemize}
\item[i)] $[\delta_1 , \delta_2 ]\in \hat\cala _{(k+l-1)}$;
\item[ii)] $(\hat\cala , [\cdot ,\cdot ])$ is a graded Lie algebra.
\end{itemize}
Moreover, the space $\cala$ of $k$-differentials, is a graded Lie
subalgebra.

\end{prop}
\begin{numrmk}
{\rm From Eqs. (\ref{k-dif-prop}) and (\ref{graded-bracket}), we
deduce that
\begin{equation}\label{property}
[\delta ,ad(P)]=ad(\delta P)
\end{equation}
for any $k$-differential $\delta \in \cala$ and any multisection
$P\in \gm (\wedge ^\bullet A)$.}
\end{numrmk}

To end this section, we note that one can introduce a graded Lie
algebra structure on $\hat {\cal A}\times \gm (\wedge A)$, where
$\gm (\wedge A)=\oplus _k \gm (\wedge ^k A)$. This is defined by
\[
[ (\delta _1,P_1),(\delta _2,P_2)] =([\delta _1,\delta _2], \delta
_1(P_2)-\delta _2 (P_1) ),
\]
for any $(\delta _1,P_1),(\delta _2,P_2)\in \hat {\cal A}\times
\gm (\wedge A)$.

Note that this is the semi-direct product Lie bracket when we
consider the natural representation of $\hat \cala$ on $\gm
(\wedge A)$.

\begin{lem}\label{diff->bracket}
If $\delta$ is a $k$-differential on a Lie algebroid $A$, then
there exists a $k$-vector field $\pi _M$ on $M$ given by
\[
\pi _M (df_1,\ldots ,df _k)=(-1)^{k+1}\langle \rho (\delta
f_1),df_2\wedge \ldots \wedge df_k\rangle ,\mbox{ for all
}f_1,\ldots ,f_k\in C^\infty (M),
\]
where $\rho :\gm (\wedge ^k A)\to \mathfrak X ^k(M)$ is the
natural extension of the anchor $\rho:\gm (A)\to \mathfrak X(M)$.
\end{lem}
\begin{pf}
Since $\delta f\in \gm (\wedge ^{k-1}A)$, we have for any $i\geq
2$
\[ \begin{array}{lll}
\{ f_1,f_2,\ldots ,f_i,f_{i+1},\ldots ,f_k\} &=& (-1)^{k+1}
\langle \rho (\delta f_1),df_2\wedge\ldots \wedge df_i\wedge
df_{i+1}\wedge
\ldots \wedge df_k \rangle \\[5pt] &=& -(-1)^{k+1}\langle \rho(\delta
f_1),df_2\wedge\ldots \wedge df_{i+1}\wedge df_i\wedge \ldots
\wedge df_k \rangle \\[5pt]&=& -\{ f_1,f_2,\ldots
,f_{i+1},f_i,\ldots ,f_k\}.
\end{array}
\]
On the other hand, since  $\delta$ is a $k$-differential, using
Eq. (\ref{k-dif-prop}) we get that, $\forall f,g\in C^\infty (M)$
\[ \begin{array}{lll}
0= \delta \lcf f,g \rcf &=&\lcf \delta f,g\rcf  +(-1)^{k+1}\lcf
f,
\delta g \rcf \\[5pt]
&=& (-1)^k d_Ag\per \delta f + (-1)^k d_Af\per \delta g,
\end{array}
\]
where $d_A$ is the differential of the Lie algebroid $A$. Thus we
can deduce that for any
 $f_1,\ldots ,f_k \in C^\infty (M)$,
\[ \begin{array}{lll}
0 &=&(-1)^{k+1} \langle df_2\per \rho (\delta
f_1),df_3\wedge\ldots\wedge df_k \rangle + (-1)^{k+1} \langle
df_1\per \rho (\delta f_2) ,df_3\wedge\ldots
\ldots \wedge df_k \rangle \\[5pt] &=& \{ f_1,f_2,\ldots ,f_k\}
+\{ f_2,f_1,\ldots ,f_k\}.
\end{array}
\]
Therefore $\{ \cdot ,\ldots ,\cdot \}$ is indeed
skew-sym\-metric. Moreover, from the fact that  both $\delta$ and
 $d$ are derivations, we can deduce that $\{\cdot ,\ldots ,\cdot \}$ is a
derivation with respect to each argument. That is, $\{\cdot
,\ldots ,\cdot \}$ induces a $k$-vector field
 $\pi _M\in \mathfrak X ^k(M)$.
\end{pf}
\begin{numex}
{\rm If $P\in \gm (\wedge ^k A)$ and $ad(P)$ is the coboundary
$k$-differential associated to $P$,  then the corresponding
$k$-vector field on $M$ is just $\rho (P)$. }
\end{numex}

\subsection{From multiplicative $k$-vector fields to $k$-differentials}
Assume that $\Pi$ is a multiplicative $k$-vector field on $\gm$.
For any $f \in C^{\infty} (M)$ and $X\in \gm (A)$, it is known
from Propositions \ref{charact-affine}, \ref{mult->affine} and
Corollary \ref{right-inv} that $[\Pi , \alp ^\ast f ] $ and $[\Pi,
\Vec{X} ]$ are right invariant, where $\Vec{X}$ denotes the right
invariant vector field on $\gm$ corresponding to $X\in \gm (A)$.
Therefore there exists $\delta_\Pi f \in \gm (\wedge^{k-1}A)$ and
$\delta_\Pi X \in \gm (\wedge^{k}A)$ such that
\begin{equation}\label{how-to-go-down}
\begin{array}{l}
\Vec{\delta _\Pi f} =[\Pi , \alp ^\ast f ] , \ \forall  f\in
C^\infty (M),
\\[7pt] \Vec{\delta _\Pi
X}=[\Pi, \Vec{X} ] , \ \forall X\in \gm (A).
\end{array}
\end{equation}

It is simple to see that $\delta _\Pi   $ is indeed a
$k$-differential. We are now ready to state the  main theorem of
the paper.

\begin{them}\label{derivation}
Assume that $\gm \toto M$ is an $\alpha$-simply connected and
$\alpha$-connected Lie groupoid with Lie algebroid $A$. Then the
map
\[
\begin{array}{ccc}
\delta : \oplus_k \mathfrak X^k_{mult} (\gm )&\to & \oplus _k \cala _k\\
\Pi &\mapsto &\delta _\Pi
\end{array}
\]
is a graded Lie algebra isomorphism.
\end{them}

We divide the proof into several steps. The surjectivity of $\delta$
will be postponed to Section 3. Here we prove the following
result.
\begin{prop}
Under the same hypothesis as in Theorem \ref{derivation}, $\delta$
is an injective graded Lie algebra homomorphism.
\end{prop}
\begin{pf}
Using the graded Jacobi identity of the Schouten brackets
 and Eq. (\ref{graded-bracket}), we deduce that if
$\Pi \in \mathfrak X ^k _{mult}(\gm)$ and  $\Pi '\in \mathfrak X
^l _{mult}(\gm)$ then
\[
\delta _{[\Pi ,\Pi ']} = [\delta _{\Pi },\delta _{\Pi '}].
\]
Therefore, $\delta$ is  a graded Lie algebra homomorphism.

Next, let us prove that $\delta$ is injective. We will use the
following lemma (see Theorem 2.6 in \cite{MackenzieX:2000}):

\begin{lem}\label{affine-uniqueness}
If $\Pi$ is an affine multivector field on an $\alpha$-connected
Lie groupoid $\gm \gpd M$ then $\Pi=0$ if and only if $\delta_\Pi X =0, \ \forall
X\in \gm (A)$,
and $\Pi$ vanishes on the unit space $M$.
\end{lem}
Suppose that $\Pi$ is  multiplicative $k$-vector field on $\gm$
such that $\delta _\Pi =0$.   It remains to show that
$\pi|_{M}=0$.
We know that $T^\ast _{\epsilon (m)}\Gamma$ is
spanned by the differential of functions of the type $\alp ^\ast
f$, with $f\in C^\infty (M)$,  and the  differential of functions
$F$ which are constant along  $M$, i.e., such that $dF \in N^* M$.
Since $M$ is coisotropic, we get that
\[
\Pi (\epsilon (m))(dF_1, \ldots ,dF_k)=0.
\]
Moreover,
\[
i_{\alp ^\ast f}\Pi = (-1)^{k+1}[\Pi ,\alp ^\ast f]
=(-1)^{k+1}\delta _\Pi (f)= 0,
\]
which implies that
\[
\Pi (\epsilon (m))(\alp ^\ast f_1, \ldots ,\alp ^\ast f_j, dF_1,
\ldots ,dF_l)=0,
\]
for $j+l=k$ and $j\geq 1$. Therefore, $\Pi_{|M}=0$.

From Lemma \ref{affine-uniqueness}, it follows that $\Pi =0$. Thus
$\delta$ is injective.
\end{pf}

Following \cite{LuW:1990}, we also call $\delta _{\Pi }$ the {\em
inner derivative} of $\Pi$.

Theorem \ref{derivation} has many interesting corollaries.
Below is a list of well known results which   are
in the literature.

\begin{numex}
{\rm Let $G$ be a simply connected and connected
 Lie group  with Lie algebra $\mathfrak g$ and
$\Pi \in \mathfrak X^k_{mult}(G)$. Then $\delta _\Pi :\mathfrak g
\to \wedge ^k\mathfrak g$ is the  1-cocycle obtained by taking the
inner derivative of $\Pi$  \cite{LuW:1990}. Thus
we have  a  one-to-one correspondence between 1-cocycles $\mathfrak g
\to \wedge ^*\mathfrak g$ and multiplicative
multivector fields on $G$ (see \cite{LuW:1990}).}
\end{numex}

\begin{numex}
{\rm If $\sigma$ is a multiplicative function on $\gm$, i.e.,
$\sigma: \gm \to \rr$ is a groupoid 1-cocycle, then
 $\delta _\sigma \in \gm (A^* )$ is exactly its corresponding Lie
 algebroid 1-cocycle. Thus we conclude that
for an $\alpha$-connected and $\alpha$-simply connected Lie
groupoid, there is one-to-one correspondence between
groupoid 1-cocycles $\sigma: \gm \to \rr$  and
 Lie algebroid 1-cocycles $\delta\in \gm (A^* )$.
This result was first proved in \cite{WX}.  }
\end{numex}

\begin{numex}
{\rm Multiplicative vector fields on a  Lie groupoid are
exactly infinitesimals  of Lie groupoid automorphisms
\cite{MackenzieX:1998}. One-differentials on a Lie algebroid,
on the other hand, are
 covariant differential operators on $A$
which are  derivations with respect to the bracket \cite{MackenzieX:1998}.
 These are exactly infinitesimals of the Lie algebroid automorphisms.
Thus  for an $\alpha$-connected and $\alpha$-simply connected Lie
groupoid, we have a one-to-one correspondence
between infinitesimals  of the Lie groupoid automorphisms  and
infinitesimals of the  corresponding Lie algebroid automorphisms
 \cite{MackenzieX:1998}.
}
\end{numex}

\begin{numex}\label{dif-cob}
{\rm Let $P \in \gm (\wedge ^k A)$ be a $k$-section of $A$, and
$\Pi =\Vec{P}-\ceV{P}$  the corresponding multiplicative
$k$-vector field on $\gm$. From the definition of $\delta _\Pi$,  we
see that $\delta _\Pi$ is just the coboundary $k$-differential
 $ad(P) =\lcf P,\cdot \rcf$. Thus we conclude that
for an $\alpha$-connected and $\alpha$-simply connected Lie
groupoid, there is a one-to-one correspondence between coboundary
multiplicative multivector fields  on the  Lie
groupoid and coboundary $k$-differentials on its Lie algebroid.}
\end{numex}

\begin{numex}
{\rm Let $(\gm \gpd M, \Pi )$ be a Poisson groupoid, i.e., $\Pi
\in \mathfrak X_{mult}^2 (\gm )$ such that $[\Pi ,\Pi ]=0$. From
Theorem \ref{derivation}, there exists a 2-differential $\delta
_\Pi$ on $A$. Moreover,
\[
\delta _\Pi ^2=\delta _\Pi \smalcirc \delta _\Pi =\frac{1}{2}
[\delta _\Pi , \delta _\Pi] =\frac{1}{2} \delta _{[\Pi ,\Pi ]}= 0.
\]
Thus, $\delta _\Pi$ defines a Lie algebroid structure on $A^\ast$.
Moreover, since $\delta _\Pi \lcf X, Y\rcf =\lcf \delta _\Pi X,
Y\rcf +\lcf X, \delta _\Pi Y\rcf $ for all $X, Y\in \gm (A)$, we
deduce that $(A,A^\ast )$ is a Lie bialgebroid. As a
consequence, we obtain the integration theorem of Mackenzie-Xu
\cite{MackenzieX:1994}:  there is a one-to-one correspondence between
$\alpha$-connected and $\alpha$-simply connected  Poisson
groupoids and Lie bialgebroids.}
\end{numex}

\section{Lifting of $k$-differentials}

This  section is devoted to the proof of  the surjectivity of
$\delta$ in Theorem \ref{derivation}.

\subsection{$A$-paths}

From now on, we use the
notation $I=[0,1]$.  Let $(A,\lcf \cdot ,\cdot \rcf ,\rho)$ be
 the  Lie algebroid of
an $\alpha$-connected and $\alpha$-simply connected Lie groupoid
$\gm$. Following \cite{Crainic}, by $\tilde{P}(A)$ we denote the
Banach manifold of all $C^1$-paths in $A$. A $C^1$-path $a :I \to
A$ is said to be an {\em $A$-path} if
 \begin{equation}\label{def:Apath}
\rho \big( a(t) \big) = \frac{\diff \gamma (t)}{\diff t},
\end{equation}
where $\gamma (t)=(p \smalcirc a)(t)$ is the base path ($p:A\to M$
is the bundle projection). The set of $A$-paths, denoted by
$P(A)$, is a Banach submanifold of $\tilde{P} (A)$.

It is well-known that integrating along $A$-paths yields a
$\Gamma$-path.
% (see Proposition 1.1 in \cite{Crainic}).
 We recall this construction in order to be self-contained. Following
\cite{Crainic}, we call a {\em $\Gamma$-path} a $C^2$-path $r(t)$
on the groupoid $\Gamma$ such that $r(0) \in M$ and
$\alpha(r(t))=r(0)$ for all $t\in I$.

There is a diffeomorphism ${\mathfrak I} $ from the space of
$A$-paths to the  space of  $ \Gamma$-paths. For any  $a \in P(A)
$, we define ${\mathfrak I}(a) $ to be the solution $r(t)$ of the initial
value problem
\begin{equation}\label{eq:calI}
\left\{
\begin{array}{ccc}
\frac{\diff r(t)}{\diff t}&=&\Vec{a(t)}_{r(t)}\\ r(0) &= &
\gamma(0)
\end{array}
\right.
\end{equation}
where $\gamma = p \circ a$ is the base path. The inverse of
${\mathfrak I}$ is, for any $\Gamma$-paths $r(t)$, given by $
{\mathfrak I}^{-1}(r(t))= L_{r^{-1}(t)\, *}    \frac{\diff
r(t)}{\diff t} $.
%Therefore we have a smooth map $\tau: P(A)\to \gm$,
%$a(t)\to r(1)$, where $r(t)={\mathfrak I}(a(t))$.

The purpose of this section is to study properties of the smooth
map $\tau$ from $P(A)$ to $\Gamma$
 given by
\begin{equation}
\label{eq:tau} \tau(a)=r(1),
\end{equation}
where $r={\mathfrak I}(a)$.

First, we study the covariance of $\tau$. We recall that the
(pseudo-) group of (local) bisections of $\Gamma \toto M$ acts on
$\Gamma$ by
\[
r \to  \tilde{g}^{-1}\cdot  r \cdot \tilde{g}
\]
for any bisection $\tilde{g}$ and $ r \in \Gamma$ (this last
expression makes sense provided that the local
bisection is chosen so that both products in the above
equation are defined). Differentiating this action with respect to
$r$, one obtains an automorphism of the
Lie algebroid  $A \to M$, which is
 denoted by $Ad_{\tilde{g}} : A_{m} \to A_{\tilde{g} \cdot
m}$. Differentiating this action with respect to
$\tilde{g}$, and using
the fact that the Lie algebra of the (pseudo-) group of bisections
is the Lie algebra of (local) sections of $A$, one constructs, for
all $\xi \in \Gamma(A) $ and $b \in A$, an element of $T_bA$,
which we denote by $ad_\xi b \in T_b A$.

The reader should not confuse $b \to ad_\xi b $,
which is a tangent vector  on $A$, with the  adjoint action
of the Lie algebra $\Gamma(A)$ on itself.

For any $\Gamma$-path $r(t)$ and any $C^2$-path $g(t)$ in
$\Gamma$, such that $r(t) $ and $g(t)$ are composable for all $t
\in I$ (i.e. $\beta \big(r(t)\big)=\alpha \big(g(t)\big)$), the
path $ g^{-1} (0) \cdot  r(t) \cdot  g(t) $ is a $\Gamma$-path
again. It is simple  to check that $ {\mathfrak
I}^{-1}\big(g(0)^{-1} \cdot r(t) \cdot g(t)\big) $ is equal to  $
Ad_{\widetilde{g_t}} {\mathfrak I}^{-1}\big(r(t)\big) +
(R_{\widetilde{g_t}^{-1}* } \frac{\diff \widetilde{g_t}}{\diff
t})_{|\tilde{g_t}  \gamma(t)}$, where $\widetilde{g_t} $ is, for
all $t \in I$, a (local) bisection of $\Gamma \toto M $  through
$g(t) $ (defined in a neighborhood of $g(t)$). In short, for any
$a \in P(A)$, and any time-dependent (local) bisection
$\tilde{g_t}$ through $g(t)$, with a $C^2$-dependence in $t$,
\begin{equation}
\label{eq:taucomp} \tau \Big (   Ad_{\widetilde{g_t}} a(t)  + (R_{
\widetilde{g_t}^{-1} *}
 \frac{\diff \widetilde{g_t}}{\diff t})_{|\tilde{g_t} \cdot \gamma(t)} \Big )=
g^{-1}(0)\cdot \tau(a) \cdot  g(1).
\end{equation}

For any $A$-paths $a_1 (t)$ and $a_2(t)$, if $\tau (a_1) =\tau(a_2) $,
then $ {\mathfrak I}(a_1) =  {\mathfrak I}(a_2) \cdot g(t) $ for
some path $g(t)$ on $\Gamma$ such that $g(0), g(1)
 \in M$. Hence, it follows from Eq. (\ref{eq:taucomp}) that
$\tau (a_1) =\tau(a_2) $ if and only if there exists a
time-dependent bisection $ \widetilde{g_t}$ such that
$$a_2  (t) = Ad_{\widetilde{g_t}} a_1 (t)  +R_{ \widetilde{g_t}^{-1}(t)\, *}
\frac{\diff \widetilde{g_t}}{\diff t}_{|\tilde{g_t}\cdot
\gamma(t)}$$
 and  $\widetilde{g_0} $ and $\widetilde{g_1}
$ are unital elements of the pseudo-group of local bisections.

Since $\Gamma$ is $\alpha$-simply connected, the pseudo-group of
time-dependent local bisections $\widetilde{g_t} $ is connected.
Therefore, $\tau(a_1)=\tau(a_2) $ if and only if $a_1$ and $a_2$
can be linked by a differentiable path in $P(A)$ which is tangent
to the differential of the action described by Eq.
(\ref{eq:taucomp}). But vectors in $T_aP(A) $ tangent to this
action are precisely the vectors of the form, at a given $A$-path
$a(t)$ with base path $\gamma(t)$,
\begin{equation}\label{eq:tauinfi}
G_{\xi\, |a}:  t \to  ad_{\xi(t)} a(t)  +
  \frac{\diff  \xi(t)}{\diff t}_{| \gamma(t)},
\end{equation}
where $\xi(t)$ is a $C^2$-time-dependent section of $A \to M$ with
$\xi (0)=\xi(1)=0$ and
 $ \frac{\diff  \xi(t)}{\diff t}_{ | \gamma(t)}    $,
an element of $A_{\gamma(t)}$,  is considered as  an element of
$T_{a(t)}A$.

For any  $C^2$-time dependent section $\xi(t)$ of $A \to M$
with $\xi (0)=\xi(1)=0$,  the vector field $G_{\xi} $  on $P(A)$
given by Eq. (\ref{eq:tauinfi}),  is called a {\em gauge vector field}.
 A smooth function $f:P(A) \to {\mathbb R}$ is   said to be {\em
invariant under the gauge transformation} if and only if
 $G_\xi ( f) =0 $ for any $G_\xi$ of the form described by Eq.
(\ref{eq:tauinfi}) with $\xi(0)=\xi(1)=0$. The following
proposition summarizes the above discussion.

\begin{prop}\label{prop:crainic}
Assume that $\gm$ is an $\alpha$-connected and  $\alpha$-simply
connected Lie groupoid. Then the map
  $\tau : P(A) \to \Gamma$ induces   an isomorphism between
${\mathcal C}^{\infty}(\Gamma)$ and the algebra of smooth
functions on $P(A)$ invariant under the gauge transformation.
\end{prop}

Note that it follows from Eqs. (\ref{eq:taucomp}-\ref{eq:tauinfi})
that for any time-dependent section $\xi(t)  $ of $\Gamma(A)$, we have

\begin{equation}
\label{eq:taucomp''} \tau_* \Big ( ad_{\xi(t)} a(t)  + \frac{\diff
\xi(t)}{\diff t}_{ | \gamma(t)}  \Big)= \Vec{\xi(1)}
-\ceV{\xi(0)}.
\end{equation}
This  relation  will be useful later on.

Next we need  to introduce regular extensions to $ \tilde{P}(A) $
 of the 1-form $\diff \tau^* f$
for  a smooth function  $f$ on $\Gamma$. First, we  give some
definitions related to the cotangent spaces of  $P(A)$ and
$\tilde{P}(A)$.

For any $a \in \tilde{P}(A)$ (resp. $P(A)$), the cotangent space
of $\tilde{P}(A)$  is denoted by $T^*_a\tilde{P}(A) $ (resp.
$T_a^*P(A)$). For any real-valued function $f$ on $\tilde{P}(A)$
(resp. on $P(A)$), the differential (if it exists) at the point
$a \in \tilde{P}(A)$ (resp. $a \in P(A)$) is an
 element of $T_a^*\tilde{P}(A)$  (resp. $T_a^*P(A)$),
and is  denoted by $ \diff f_{|_a}$.

For any $a \in \tilde{P}(A)$, denote by $P_a(T^*A)$, the space of
$C^1$-maps $\eta :I\to T^*A $
 such that $\forall t \in I$ the identity
 $    \pi \smalcirc \eta (t) =  a(t) $ holds,
where $\pi $  denotes the projection $T^*A \to A$. The space
$P_a(T^*A)$ can be considered as a sub-space of $T^*_a \tilde{P}(A)$,
 by associating to any $\eta(t)  \in P_a (T^*A)$
 the linear map
\begin{equation}\label{eq:pathdual} \begin{array}{ccc}  T_a\tilde{P}(A)
 & \to & {\mathbb R}\\  X(t)  &\to & \int_{0}^1  \langle\eta(t),
 X(t) \rangle dt , \end{array}\end{equation}
considered as an element of $T^*_a \tilde{P}(A)$,
 where  $\langle \ , \ \rangle$
denotes  the pairing between the cotangent and tangent
vectors. Throughout  this section, we will always consider $P_a(T^*A)$
as a subspace of $T^*_a \tilde{P}(A)$, and, therefore, the vector
bundle $P(T^*A) \to \tilde{P}(A)$ as a vector sub-bundle of $T^*
\tilde{P}(A) \to \tilde{P}(A)$. We denote by $ P(T^*A)_{|_{P(A)}}
\to P(A) $ and $T^* \tilde{P}(A)_{|_{P(A)}} \to P(A) $ the
restrictions of these vector bundles to $P(A)$.

\begin{defn}
Given a  1-form $\omega$ on the Banach submanifold $P(A)$ of
$A$-paths, by an {\em  extension (resp. regular extension)},
 we mean a smooth section  $\Phi_\omega
$ of $ T^* \tilde{P}(A)_{|_{P(A)}} \to P(A)  $
 (resp. $P T^* A_{|_{P(A)}} \to P(A)$) such that $\Phi_{\omega \, |_a} $,
 is for any $a \in P(A)$,  an  extension  of $\omega_{|_a}$ (i.e.,
  the restriction of $ \Phi_{\omega \, |_a} $
 to $T_a P(A)$ is $\omega_{|_a}$ for any  $a \in P(A)$).
\end{defn}

By  a {\em regular extension of a smooth function $g\in C^{\infty}(P(A))$},
we mean a regular extension of its
differential.
A regular extension of the zero function
on $P(A)$ will be called a {\em regular extension
of zero}. Also, we use the following notation:
for any regular extension $\Phi_{\omega}$, we denote by
$\Phi_{\omega}(t)$ the corresponding
path in $P_a
T^* A$.

Given a vector field $X$ on $\tilde{P}(A)$  tangent to
$P(A)$,  and a 1-form $\omega$ on $P(A)$, the Lie derivative of a
regular extension $\Phi_\omega $  of $\omega$ is defined by
\begin{equation}\label{eq:Lieder}
({\cal L}_X \Phi_{\omega}) (Y) = \Phi_{\omega} ([Y,X]) + X \big(
\Phi_\omega (Y)\big),
\end{equation}
where $Y $ is the restriction to $P(A)$ of a vector field on
$\tilde{P}(A)$. This definition needs to be justified. It is clear
that the right hand side of Eq. (\ref{eq:Lieder})  is
$C^{\infty}$-linear with respect to $Y$ and depends only on  its
restriction  to $P(A)$. It therefore
defines a section of $ T^* \tilde{P}(A)_{|_{P(A)}} \to P(A)  $. It
follows from Eq. (\ref{eq:Lieder})  that if $Y$  itself
is  tangent
to $P(A)$, then $({\cal L}_X \Phi_{\omega}) (Y)=({\cal L}_X
\omega) (Y)$. Therefore,  ${\cal L}_X \Phi_{\omega}$  is an
extension of ${\cal L}_X \omega$.

In the following subsections, we need to investigate  regular
extensions $ \Phi_{\diff \tau^* f}$ of $\diff \tau^* f$ for  a smooth function
$f\in C^{\infty}(\gm )$.
 The following technical lemma will be  very useful.

\begin{lem}
 \label{th:recap4.1}  $\forall a\in P(A) $,
\begin{enumerate}
\item[i)] for any smooth function $f:\Gamma \to {\mathbb R}$,
the pull-back  function $\tau^* f : P(A) \to {\mathbb R} $ admits  a regular
extension $\Phi_{\diff  \tau^* f }$;
\item[ii)] for any smooth function $f:\tilde{P}(A) \to {\mathbb R} $
whose restriction to $P(A) $ vanishes, there exists $g_t : I \to
{\mathcal C}^{\infty}(M) $ with $g_0=g_1=0 $ such that $\diff f_{|_a}
= \diff  {\mathcal F}_{\diff g \, |_a} $. Here for any time dependent
$1$-form $\omega: t \to \omega_t$ on $M$ and any $a \in \tilde{P}(A)$
with base path $\gamma(t)$,
\[
{\mathcal F}_{\omega}(a)= \int_{I} \langle \omega_{t
|_{\gamma(t)}},\frac{\diff \gamma(t)}{\diff t}- \rho(a)\rangle dt
.
\]
\end{enumerate}
\end{lem}
\begin{pf}
{\em i)} We divide the proof of {\em i)} into
 four steps. In what follows, $\gamma$ always denotes the base path
 of  $a \in P(A)$. We say that a   path $e(t)$ in a vector bundle $p:E \to
M$ is {\em over a base path $\gamma$} if $p \circ e = \gamma$.

{\em Step 1.} We  first describe  the differential $\tau_*$ of the
map $\tau$ defined in Eq. (\ref{eq:tau}).

The map $(r, a) \to \Vec{a}_{|_{r}}$ is a map  from the fibered
product $\Gamma \times_{\alpha ,M,p} A$ to $T \Gamma$ that maps a
pair $(r,a)$ (with $r \in \Gamma, a \in A$) to an element of $T_r
\Gamma$. Considering the fibered product $\Gamma \times_{\alpha
,M,p} A$ as a submanifold of $\Gamma \times A$, one can extend
this map to a smooth map $F(\alpha,a)$ from $ \Gamma \times A $ to
$ T\Gamma$ satisfying  the condition $F(r,a) \in T_{r}\Gamma $ for
all $r \in \Gamma$. This allows us to rewrite  $\tau$ more
conveniently. Consider the map $\tilde{\tau}: \tilde{P}(A)\times M
\to \Gamma$ given by
\begin{equation}
\label{eq:tildetau} \tilde{\tau}(a,m)=r(1) ,
\end{equation}
 where $r(t): I \to \Gamma$ is the solution of
the initial value problem
$$ \left\{ \begin{array}{ccc} \frac{\diff  r(t)}{ \diff t}&=&
F\big(r(t),a(t)\big) \\
r(0) &= & m   \end{array}\right. $$
By definition,  we have $\tau (a)=
\tilde{\tau}\big(a, (p\smalcirc a)(0)\big) $.

Linearizing the previous equation,
one obtains that the differential
of the map $ \tilde{\tau}$ can be written under the form
\begin{equation}
\label{eq:dtildetau}
 \tilde{\tau}_* (  \delta a, \delta m )=
L\big(\delta m (0)\big)
 + \int_{t \in I} L_t \big( \delta a (t) \big) dt,
 \end{equation}
where $\delta a , \delta m $ are elements of $T_a\tilde{P}(A)$ and
$ T_mM$ respectively, and $L_t,L$ are linear maps from $
T_{r(t)}%{\alpha(t)}
\Gamma$ and $T_mM$ to $ T_{\tilde{\tau}(a,m)} \Gamma$ respectively.
Thus  the differential of $\tau :
P(A) \to \Gamma$ can  be expressed under the form
 \begin{equation}
\label{eq:dtau}
 \tau_* (  \delta a )=
 \int_{t \in I} L_t \big( \delta a (t) \big) dt  + L(p_* \big(\delta a(0)\big) )
 \end{equation}
where $\delta a \in T_a P(A) $.

At this point, we need to  explain why some technical difficulties
arise in the construction of a regular extension of $\tau^* f $
and what remains to be done in order to avoid it. It follows from
Eq. (\ref{eq:dtau}) that the differential  $\diff \tau^* f =\diff
f \circ \tau_*$ is the sum of two $1$-forms, namely $\omega_1:
\delta a \to  \int_{t \in I} L_t \big( \delta a (t) \big) dt$ and
$ \omega_2:\delta a \to  L(p_* \big(\delta a(0)\big))$.  It is
easy to find a regular extension of $\omega_1$: we can just
choose, for any $\delta a \in T_a  \tilde{P}(A)$,
\[
\Phi_{\omega_1} (\delta a) =     \int_{t \in I} L_t \big( \delta a
(t) \big) dt .
\]
Unfortunately, it is not so easy to find a regular extension of
$\omega_2$, since this $1$-form is  ``concentrated" in $0$.

{\em Step 2.} We describe explicitly the tangent space $T_a P(A)$ of
$P(A)$.

   Let us choose a connection $\nabla^A$ on the vector bundle $A \to M$.
This connection allows us to decompose the tangent space of $A$ as
a direct sum $T_b A = T_{p(b)}M \oplus A_{p(b)}$ for any $b \in
A$. With this convention, for any $a \in \tilde{P}(A)$, an element
$\delta a $  of $T_a \tilde{P}(A)$ becomes a pair $\delta a=
(\epsilon,\beta) $ where $\epsilon: I \to TM$ and $\beta: I \to A$
are $C^1$-maps over the base path $\gamma = p \circ a$.

    We choose now a  connection $\nabla^M$ on
the tangent bundle $TM \to M$. By differentiating the relation
$\frac{\diff \gamma}{\diff t}= \rho(a)$, we obtain that $\delta a
= (\epsilon,\beta)\in T_a \tilde{P}(A)  $ is an element of the
tangent space of $T_a P(A) $ if and only if:

\begin{equation}\label{eq:tangentspace}
\nabla^M_{\frac{\diff \gamma}{\diff t}} \epsilon =
(\nabla_\epsilon \rho)  (a)+ \rho(\beta),
\end{equation}
where $(\nabla_{\epsilon} \rho )(a)= \nabla^M_\epsilon \rho(a)-
\rho (\nabla^A_\epsilon a)$ is a path in $TM$ (over $\gamma$ again).

{\em Step 3.} For any A-path $a$,
%we want to project an element $\delta a \in T_a\tilde{P}(A)$
%to the tangent space of the starting point of the base path
%$\gamma = p \circ a$ in a ``regular'' enough way. More precisely,
we want to construct
 a linear map $\Pi_a$ from $T_a\tilde{P}(A)$ to $T_{\gamma(0)}M$,
where $\gamma = p \smalcirc a$ is the base path,
depending smoothly on $a \in P(A)$, whose restriction to $
T_aP(A)$ is simply the differential $ a \to \gamma (0)$, and
which is ``given by an integral". That is,
the differential  is of the form
\begin{equation}\label{eq:Pa}
\Pi_a (\delta a) = \int_{t \in I} M_t \big( \delta a (t)\big) dt,
\end{equation}
where $M_t$ is,  for all $t \in I$,
 a linear map from $ T_{a(t)}A$ to $T_{\gamma(0)} M $. We  proceed as follows.
Set $\delta a = (\epsilon,\beta) \in T \tilde{P}(A)$ as in Step 2.

First, for any $ s \in I$, we define $\eta_{s}(t): I\to TM$
 to be the unique solution of the initial value problem:
\begin{equation} \left\{ \begin{array}{ccc}
\nabla^M_{\frac{\diff  \gamma (t)}{\diff  t}} \eta_s(t)  &=&
(\nabla_{\eta_s (t)} \rho)  (a)+ \rho(\beta(t))\\
 \eta_{s}(s)   &=&  \epsilon (s)  \\
\end{array} \right. \end{equation}
(this equation is a linear equation of order 1, which guarantees
the existence and uniqueness of the solution). Then we define
$\Pi_a(\epsilon,\beta)$ by
\[
\Pi_a(\epsilon,\beta) = \int_{s \in I} \eta_{s} (0)  ds  .
\]
It follows from a classical result of ordinary linear differential
equations (see, for instance, \cite{H}) that
$$\eta_{s} (0) = L(s) \big(\epsilon (s) \big) +
\int_{s}^0 M(s,u) (\beta(u))du$$
 for some smooth function $M(s,u)$
from $A_{\gamma(u)} $ to $T_{\gamma(0)}M$. It is  simple
to check  that for any $\delta a = (\epsilon,\beta)$
$$ \Pi_a(\delta a) =   \int_{s \in I} L(s) \big( \epsilon (s) \big) ds -
   \int_{s \in I}\int_{u=0}^s  M(s,u) \big(\beta(u)\big)   du ds . $$

The right-hand side of this  equation is of  the form
given  by Eq. (\ref{eq:Pa}). Now it remains to check  that the
restriction  of $\Pi_a$ to $TP(A)$  is equal to the map $ \delta a
\to p_* \big(\delta a (0)\big)$. For any  $\delta =
(\epsilon,\beta)$ tangent to $P(A)$, by the
uniqueness of the solution of an initial value problem,
 we have $ \eta_{s} = \epsilon$ for all $s \in I$. Therefore, the restriction of $\Pi_a$
to $T_a P(A)$ is equal to
\[
\Pi_a(\delta a) = \int_{I} \epsilon(0)   dt = \epsilon (0) =
p_*\big(\delta a (0)\big).
\]

{\em Step 4.} We now can define for any $f \in C^\infty (M)$
\[
\Phi_{\diff \tau^* f} ( \delta a) := \tilde{\tau}_* \big( \delta a
, \Pi_a (\delta a ) \big)
\]

It follows from Eqs. (\ref{eq:tildetau}-\ref{eq:Pa}) that $\Phi_{d
\tau^* f}$ is a regular extension of $\tau^* f$.
This achieves the proof of  {\em i)}.

{\em ii)}  We identify an element $\delta a \in \tilde{P}(A)$ with
 a pair  of paths $\epsilon (t)$
 in $ TM$,  and $\beta (t)$ in $A$ over a  base path $\gamma (t)$ as
previously.

Since $\Phi_{\diff  \tau^* f} $ is a regular extension, we have
$$  \Phi_{\diff  \tau^* f} (\epsilon,\beta) =\int_I \big< M (t), \epsilon(t) \big> dt
 + \int_I \big< A(t) ,\beta(t) \big> dt  $$
for some $C^1$-maps $ M(t) , A (t)$ from $I $ to $T^*M$ and
 $A^*$ respectively, which are   over the base path $\gamma (t)$.

Let $\omega (t): t \to T^*_{\gamma(t)} M $ be a path over $\gamma$
which is a solution of the initial value problem
\[
\left\{
\begin{array}{ccc}
\nabla^M_{\frac{d \gamma(t)}{dt}}  \omega (t)
+ \big( (\nabla
  \rho) a(t)\big)^* \big(\omega(t)\big) &=& M(t) \\ \omega(1)
& = & 0
\end{array}
\right.
\]
where, for any fixed $t\in I$,
 $\big( (\nabla \rho )a(t)\big)^* \in \mbox{End}(T^*_{\gamma(t)} M)$ is
the dual of the endomorphism $v  \to (\nabla_{v } \rho)
\big( a(t) \big) $ of $T_{\gamma (t)}M$.

Using integration  by part, we obtain
$$ \int_I \big< M(t),\epsilon(t) \big>dt=
- \int_I \big< \omega(t), \nabla^M_{\frac{d \gamma(t)}{dt}}
 \epsilon (t) \big>dt +  \int_I \big<
\omega(t),
 (\nabla_{\epsilon(t)} \rho)  (a(t)) \big> dt + \big<\omega (0),  \epsilon
 (0)\big>.$$
Since $ \Phi_{d\tau^* f} (\epsilon,\beta) $  vanishes as long
 as the conditions
$$ \left\{ \begin{array}{cccc} \beta (t)&=&0 & \mbox{ and } \\
 \nabla_{\frac{\diff \gamma(t)}{\diff
t}}^M \epsilon(t)-(
 \nabla_{\epsilon(t)} \rho)  a(t)&=&0 & \\ \end{array} \right.  $$
 are satisfied,
 we must have $\omega(0) =0$.

Now, since
$$  \Phi_{\diff \tau^*  f} (\epsilon,\beta)  =
  \int_I \big< \omega(t), \nabla^M_{\frac{\diff \gamma(t)}{\diff t}} \epsilon
  (t)+( \nabla_{\epsilon(t)} \rho)  a(t) \big> dt
+\int_I \big< A(t), \beta(t) \big> dt   .$$
must vanish whenever $\epsilon,\beta$ satisfy Eq.
(\ref{eq:tangentspace}), we must have  $ A (t)=\rho^* \omega(t)
$. As a consequence, we have
$$  \Phi_{\diff f} \big(\delta a \big)  = \int_I \big< \omega(t) ,
 \nabla^M_{\frac{\diff  \gamma(t)}{\diff   t}} \epsilon(t)   - ( \nabla_{\epsilon(t)} \rho)
  a(t)+ \rho \big( \beta (t) \big)   \big>  dt $$

According to  Lemma \ref{lem:annex} {\em i)}, there exists a
family of time dependent  functions $g_t$
from $I$ to $C^\infty (M)$ vanishing in $t=0,1$ such that $\diff
g_{t \, |_{\gamma(t)}} =\omega(t) $. The result now follows from
Lemma \ref{lem:annex} {\em ii)}.
\end{pf}

\begin{lem}\label{lem:annex}
\begin{itemize}
\item[i)] For any $C^1$-path $\omega(t): I \to T^*M$,
there exists a time-dependent function $g:t \to g_t$ vanishing  at
$t=0,1$ such that $\diff g_{t \, |_{\gamma(t)}} =\omega(t) $ for
any $t \in ]0,1[$.
\item[ii)]  The function $a \to \diff  {\mathcal F}_{\diff g \,  |_a}  $
is a regular extension of zero whose differential  is of the form
\begin{equation}
\label{eq:extension}
 \diff  {\mathcal F}_{\diff g \, |_a} (\epsilon,\beta)
  = \int_I \big< \omega(t) ,
\nabla^M_{\frac{\diff  \gamma(t)}{\diff t }} \epsilon(t)   - (
\nabla_{\epsilon(t)} \rho)
  a(t)+ \rho \big( \beta (t) \big)   \big>  dt,
\end{equation}
where $\epsilon, \beta ,\nabla^M, \nabla$ are as in Eq. (\ref{eq:tangentspace}).  \end{itemize}
\end{lem}
\begin{pf}
%We invite the reader to skip these technical proofs in a first reading.
{\em i)} We denote by $exp_m: T_mM \to M$ the exponential map
associated to the connection $\nabla^M$.
% (that is to say $exp_m(e)=
%g_e(1)$ where  for all $e \in T_mM$, $g_e(t)$ is the unique
%geodesic starting from  $m$ with initial velocity $\frac{\diff
%g_e}{dt}_{|_{t=0}} =u$). We
Recall that $exp_{m}$ is a local
diffeomorphism from a neighborhood of $0 \in T_mM$ to a
neighborhood of $m$. Let  $f(t)$ be a real-valued smooth function
that vanishes at $t=0$ and $t=1$.

Define $g_t(m) = f(t)  \psi(t,m) \langle \omega(t)
,exp_{\gamma(t)}^{-1}(m) \rangle$,  where $\psi(t,m) $ is a smooth
function on $I \times M$ satisfying   two conditions:
 (1) $\psi(t,m) $ is  identically equal to $1$ in a
neighborhood of the curve $(t,\gamma(t))$
 for all $t \in I$ (2) $\psi(t,m) $ vanishes outside
an open set on which $exp_{\gamma(t)}^{-1}$ is well-defined. We
leave it to the reader to check that these conditions can be
satisfied and  that the function $g_t$ has the requested
properties.

{\em ii)}  The function ${\mathcal F}_{dg}$ is identically zero on
$P(A)$ by definition. One checks that its differential
  is of the form (\ref{eq:extension}) by a direct computation.
\end{pf}

\subsection{Almost differentials and linear multivector fields}

In this subsection,  we study  a particular case of multivector
fields on a vector bundle.

Let $p:A\to M$ be a vector bundle. It is clear that there exists a
bijection between the space $\gm (A^\ast )$ of sections of the
dual bundle $p _\ast :A^\ast \to M$ and the set $C^\infty
_{lin}(A)$ of functions which are linear on each fiber. In fact,
for any section $\phi \in \gm (A^\ast )$ the corresponding linear
function $\ell _\phi$ is given by $\ell _\phi (X_m)=\langle \phi
(m),X_m\rangle \ \forall X_m\in A_m$. On the other hand, by  basic
functions, we mean functions on $A$
 which are  the pull-back  functions from  $M$.

\begin{defn}
Let $p:A\to M$ be a vector bundle and $\pi \in \mathfrak X ^k(A)$
a $k$-vector field on $A$. We say that $\pi$ is {\em linear} if
$\pi (df_1,\ldots ,df_k )$ is a linear function whenever all
$f_1,\ldots ,f_k\in C^\infty (A)$ are linear functions.
\end{defn}
Linear multivector fields were called homogeneous in \cite{GU1}.
They  satisfy the following properties.
\begin{prop}\label{prop:1cocycle}
Let $p:A\to M$ be a vector bundle and $\pi \in \mathfrak X ^k(A)$
a linear $k$-vector field on $A$ with $k \geq 2$. Then $\pi (
d\ell _{\phi _1}, \ldots ,d\ell _{\phi _{k-1}},d(p^\ast f_1 ))$ is
a basic function on $A$ and
\[
(d(p^\ast f_1)\wedge d(p^\ast f_2))\per \pi =0 \ \ \ \forall \phi
_1,\ldots ,\phi _{k-1} \in \gm (A^\ast ) \  \mbox{ and } f_1, \ f_2
\in C^\infty (M).
\]
\end{prop}
\begin{pf}
Proceed as in \cite{Co}, using Leibniz rule and the fact that
$(p^\ast f) \ell _\phi=\ell _{f\phi}$ for $f\in C^\infty (M)$ and
$\phi \in \gm (A^\ast )$.
\end{pf}

Let $\pi \in \mathfrak X _{lin}^k(A)$ be a linear $k$-vector field
on $A$.   Proposition \ref{prop:1cocycle} enables us to
introduce a  pair of operations
\[ \begin{array}{l}
\rho _\ast :\gm (A^\ast)\times \stackrel{(k-1)}{\ldots}\times \gm
(A^\ast )\to \mathfrak X (M) \ \ \mbox{ and} \\
[5pt] [ \cdot , \ldots , \cdot ] _\ast :\gm (A^\ast)\times
\stackrel{(k)}{\ldots}\times \gm (A^\ast )\to \gm (A^\ast)
\end{array}
\]
by
\begin{equation}\label{relations}
\begin{array}{l}
 \rho _\ast ( \phi _1,\ldots ,\phi _{k-1} )(f)\smalcirc p
 = \pi ( d\ell _{\phi_1}, \ldots ,d\ell _{\phi _{k-1}},d(f \smalcirc p)),\\[5pt]
\ell _{[ \phi_1,\ldots ,\phi _k ] _\ast}=\pi (d\ell _ {\phi _1}, \ldots ,d\ell
_{\phi _k}),
\end{array}
\end{equation}
$\forall \phi _1,\ldots ,\phi _k \in \gm (A^\ast )$ and $f\in
C^\infty (M)$.

It is simple to see that the following identities are satisfied:
\begin{equation}\label{deriv-prop}
\begin{array}{l}
\rho _\ast (\phi _1,\ldots ,f \phi _i, \ldots ,\phi _{k-1})= f\rho
_\ast (\phi _1,\ldots , \phi _i, \ldots ,\phi _{k-1}),\\[7pt]
\begin{array}{lcl} [\phi _1,\ldots ,f \phi _i, \ldots ,\phi _k ]_\ast
&=& (-1)^{i+1}\rho _\ast (\phi _1,\ldots , \hat{\phi} _i, \ldots
,\phi
_{k-1})(f)\phi _i\\[5pt]&&+
f[\phi _1,\ldots , \phi_i, \ldots ,\phi _k ] _\ast
\end{array}
\end{array}
\end{equation}
for any $\phi _1, \ldots ,\phi _k \in \gm (A^\ast )$ and $f\in
C^\infty (M)$. In particular, we see that $\rho_* $ induces a
bundle map $\wedge^{k-1} A^* \to TM$.

Conversely, if we have a pair $(\rho _\ast ,[ \cdot , \ldots ,
\cdot ]_\ast )$ satisfying Eq. (\ref{deriv-prop}), we can define a
linear $k$-vector field on $A$  using Eq. (\ref{relations}).

Now, assume that $\delta$ is an almost $k$-differential. We
construct a pair $([ \cdot , \ldots , \cdot ] _\ast ,\rho _\ast )$
satisfying Eq. (\ref{deriv-prop}) as follows. Let

\[
\begin{array}{lll}
\rho _\ast (\phi _1,\ldots ,\phi _{k-1})(f)&=&\langle \delta
f,\phi _1\wedge \ldots \wedge \phi _{k-1}\rangle , \ \ \mbox{ and}\\
[10pt] \langle [ \phi _1 , \ldots , \phi _k ]_\ast ,X\rangle &=&
\displaystyle \sum _{i=1}^k (-1)^{i+k} \rho _\ast (\phi _1 ,
\ldots , \hat{\phi}
_i,\ldots , \phi _k) (\phi _i(X))\\[5pt] && - \langle \delta X,
\phi _1\wedge \ldots \wedge \phi _k\rangle ,
\end{array}
\]
$\forall \phi _1,\ldots ,\phi _k \in \gm (A^\ast )$, $X\in \gm
(A)$ and $f\in C^\infty (M)$. Conversely, using these equations,
one  can construct an almost  $k$-differential from $([ \cdot ,
\ldots , \cdot ] _\ast ,\rho _\ast )$.

A combination of the above discussion leads to the following

\begin{prop}
For a given Lie algebroid $A$, there is a one-to-one
correspondence between linear multivector fields and almost
differentials on $A$.
\end{prop}

For an almost differential $\delta$ on $A$, we  denote its corresponding
linear multivector field by $\pi_\delta$. In local coordinates,
the correspondence between almost $k$-differentials and $k$-vector
fields $\pi_{\delta}$ on $A$ can be described as follows. Let
$(x_{1}, \cdots , x_n )$ be  local coordinates on $M$ and $\{ e_1,
\cdots , e_s\}$ a local basis of $\gm (A)$. Assume that \be
&&\delta x_i =\sum a_{i}^{i_{1}\ldots i_{k-1}}(x) e_{i_{1}}\wedge
\ldots \wedge e_{i_{k-1}},\\ &&\delta e_i =\sum c_{i}^{i_{1}\ldots
i_{k}}(x) e_{i_{1}}\wedge \ldots \wedge e_{i_{k}}, \ee
 Then
$$\pi_{\delta}=\sum  a_{i}^{i_{1}\ldots i_{k-1}}(x) \parrr{v_{i_1}}
\wedge \ldots \wedge
\parrr{v_{i_{k-1}}} \wedge \parrr{x_{i}}-c_{i}^{i_{1}\ldots i_{k}}(x)
v^{i}\parrr{v_{i_1}} \wedge \ldots \wedge
\parrr{v_{i_k}} ,$$
where $(v_{1}, \ldots , v_{s})$ are the corresponding linear
coordinates on the fibers.

The correspondence between almost differentials and linear
multivector fields preserves the graded Lie algebra structure, as
shown in the following

\begin{prop}\label{homomorphism}
Given a Lie algebroid $A$, the map $\delta \to \pi _\delta$
is a graded Lie algebra isomorphism from almost differentials on $A$ to linear multivector
fields on $A$. I.e., for any almost differentials $\delta_1$ and
$\delta_2$, we have,
$$[\pi_{\delta_1},
\pi_{\delta_2} ]=\pi_{[\delta_1, \delta_2 ]}.$$
\end{prop}
\begin{pf}
Note  that if $P\in \gm (\wedge ^p A)$ is  given locally  by
$\sum_{j_1, \cdots, j_p}
P^{j_1\ldots j_p}(x) e_{j_1}\wedge \ldots \wedge e_{j_p}$ under
 a given local basis $\{ e_1,\ldots ,e_s\}$ of $\gm (A)$, then
\[
\begin{array}{rcl}
\delta P &=& \displaystyle\sum _i \frac{\partial P^{j_1\ldots j_p}
(x)}{\partial x_i}
 a_{i}^{i_{1}\cdots i_{k-1}}(x) e_{i_{1}}\wedge
\cdots \wedge e_{i_{k-1}}\wedge
e_{j_1}\wedge \ldots \wedge e_{j_p}\\[8pt]&&+\displaystyle\sum _{i}(-1)^{(i+1)(k+1)}
P^{j_1\ldots j_p} (x)  e_{j_1}\wedge \ldots \wedge \delta e_{j_i}
\wedge \ldots \wedge e_{j_p} .
\end{array}
\]
The assertion follows from a tedious computation and is left to
the reader.
\end{pf}
\begin{numex}
{\rm If $\delta $ is an almost 0-differential, we know that
$\delta$ is just the contraction operator by an element
 $\phi\in \gm (A^\ast
)$. In this case, one shows that $\pi _\delta \in \mathfrak X
_{lin}^0 (A)=C^\infty _{lin}(A)$ is just $-\ell_\phi$. }
\end{numex}
\begin{numex}
{\rm When $k=1$, as  a consequence of Proposition
\ref{homomorphism}, one obtains a Lie algebra isomorphism between
$\gm(CDO(A))$, the space of covariant differential operators on
$A$, and $\gm ( T^{LIN}(A))$, the space of  linear vector fields
on $A$ (see \cite{MackenzieX:1998}). }
\end{numex}
\begin{numex}
{\rm Let $\delta$ be a 2-differential of square zero. We know that
$\delta$ induces a Lie algebroid structure on $A^\ast$ (see
\cite{Kosmann,Xu:1999}). On the other hand, from Proposition
\ref{homomorphism}, it follows  that $[\pi_\delta ,\pi _\delta
]=0$,  and therefore $\pi_\delta$ defines a Poisson structure on
$A$. Such a correspondence between Lie algebroid structures on
$A^\ast$ and linear Poisson structures on the dual bundle is
standard (see \cite{Co} for instance). A generalization to
arbitrary linear 2-vector fields and pre-Lie algebroids was
considered in \cite{GU}.}
\end{numex}

%\begin{numex}\label{complete-lift}
We end this subsection by recalling two kinds of liftings
from  $\gm (\wedge ^\bullet A)$ to multivector fields
on $A$, the complete and vertical lifts.
 Let  $P\in \gm (\wedge ^p A)$,
and $ad(P)$ its  corresponding coboundary differential. Then the
corresponding linear p-vector field $\pi _{ad(P)}$ is the
so-called {\em complete lift} of $P$, which is  denoted by $P^c$
(see \cite{GU1}).

On the other hand,  for any multisection $P \in \gm (\wedge ^p A)$,
 there exists another
kind of lift, called the {\em vertical lift}, denoted by $P^v$.
It  is obtained via   the natural  identification $A_m \simeq T_a^{Vert}(A)$,
 $\forall a\in A$. In
local coordinates, if $P=P^{j_1\ldots j_p}(x)e_{j_1}\wedge \ldots
\wedge e_{j_p}$, for a given local basis $\{ e_1,\ldots ,e_s\}$,
then $$P^v= P^{j_1\ldots j_p}(x)\frac{\partial }{\partial
v_{j_1}}\wedge \ldots \wedge \frac{\partial }{\partial v_{j_p}}.$$

Complete and vertical lifts satisfy the following properties
\cite{GU1}:

\begin{eqnarray}
&&\label{complete1}
f^v=p^\ast f;\\
&&\label{complete2}
f^c=\ell _{d_Af}; \\
&&\label{complete3}
(P\wedge Q)^c=P^c\wedge Q^v+P^v\wedge Q^c;\\
&&\label{complete35}
(P\wedge Q)^v=P^v\wedge Q^v;\\
&&\label{complete4}
[P^c,Q^c] =\lcf P,Q \rcf^c;\\
&&\label{complete5}
[P^c,Q^v] =\lcf P,Q\rcf ^v;\\
&&\label{complete6} [ P^v,Q^v] =0;
\end{eqnarray}
$\forall f\in C^{\infty}(M)$ and $P,Q\in \gm (\wedge^{\bullet} A)$.

\begin{prop}\label{dif+complete+vertical}
Let $\delta$ be a $k$-differential and $P\in \gm (\wedge^\bullet
A)$. Then
\[
\begin{array}{c}
[\pi _\delta ,P^c]=(\delta P)^c,\\[5pt] [\pi _\delta
,P^v]=(\delta P)^v.
\end{array}
\]
\end{prop}
\begin{pf}
Using Eq. (\ref{property}) and  Proposition \ref{homomorphism}, we have
\[
[\pi _\delta ,P^c]=[\pi _\delta ,\pi _{ad(P)}]=\pi _{[\delta
,ad(P)]}= \pi _{ad(\delta P)}=(\delta P)^c.
\]
The second identity can be checked directly using local coordinates.
\end{pf}

Let
$$P\Gamma (\wedge ^\bullet A)=\{ I\ni t\mapsto P(t) \in \gm ( \wedge ^\bullet A) :
P (t) \mbox{ is of class } C^2\mbox{ in }t\}$$ and
$$P_0\Gamma (\wedge ^\bullet A)=\{ I\ni t\mapsto P(t) \in \gm (\wedge ^\bullet A) :
P(0)=P(1)=0, \ P(t) \mbox{ is of class } C^2\mbox{ in }t\} .$$
 If $P\in P_0\gm (\wedge ^\bullet A)$, we define $G(P)$ as the
time-dependent multivector field  on $A$ given by:
\begin{equation}\label{c+v}
G(P)= P^c+ \Big ( \frac{dP}{dt}\Big )^v.
\end{equation}

The multivector fields $G(P)$ satisfy the following properties.
\begin{prop}\label{properties-G}
If $P,Q \in P\gm (\wedge ^\bullet A)$ and $\delta$ is
 a  multi-differential on $A$,
 then
\begin{eqnarray}
&&[G(P),G(Q)]=G([P,Q]), \label{eq:1}\\
&&G(P\wedge Q)= G(P)\wedge Q^{v}+P^v \wedge G(Q),  \label{eq:2}\\
&&[\pi_{\delta},G(P)]=G(\delta(P)) \label{eq:3}.
\end{eqnarray}
In particular, Eq.~(\ref{eq:2}) implies that
\begin{equation}\label{eq:fpsi}
G\big(f(t) P\big)= f(t) G(P)+\frac{df(t)}{dt}P^v .
 \end{equation}
% for any smooth function $f: I \to {\mathbb R} $ with $f(0)=f(1)=0 $.
\end{prop}
\begin{pf}
Using Eqs.~(\ref{complete4}), (\ref{complete5}) and
(\ref{complete6}) and that $\frac{d}{dt}$ is a derivation with
respect  to the Schouten bracket, we  deduce that Eq. (\ref{eq:1})
holds.

Eq. (\ref{eq:3}) is a direct consequence of Proposition
\ref{dif+complete+vertical} and  the fact that $\delta $ commutes
with  $\frac{d}{dt}$.
\end{pf}

\subsection{From $k$-vector fields on $A$ to $k$-vector fields
on $\tilde{P}(A)$}

We start this subsection with a construction that must be thought,
at least heuristically, as a lifting of  a time-dependent
multivector fields on $A$ to  a multivector field on
$\tilde{P}(A)$.

%Note that for any $ a\in \tilde{P}(A)$,
%$T^*_a(\tilde{P}(A)) $ is the set of
%smooth maps \comment{$C^1$??}
%$\eta :I \to T^* A$
%such that $\eta (t)\in T^\ast _{a(t)}A$ for any $t\in I$.

Let $\pi (t)$ be a  time-dependent $k$-vector field  on $A$ with
$k\geq 1$. For given $k$ one-forms
 $\eta_1,\dots,\eta_k \in \Omega^1 \big( P(A) \big)$
and their  regular extensions $\Phi_{\eta_1},\dots,\Phi_{\eta_k}$,
define a  smooth function
$\tilde{\pi}( \Phi_{\eta_1},\dots,\Phi_{\eta_k})$ on $P(A)$
 as follows. For any  $a\in P(A)$,
\begin{equation}\label{eq:tilde}
\tilde{\pi}( \Phi_{\eta_1},\dots,\Phi_{\eta_k})(a):= \int_I \pi
(t)_{|a(t)} \big(
\Phi_{\eta_1}(a(t)),\dots,\Phi_{\eta_k}(a(t))\big) dt
\end{equation}
\begin{numrmk}\label{rmk:onenotregular}
 Eq. (\ref{eq:tilde}) still makes sense
when one of the extension $\Phi_{\eta_i} $ is not necessary
regular, with the following modification of its definition. Assume
that $\Phi_{\eta_1} $ is not necessary regular, and that
$(\Phi_{\eta_2},\dots, \Phi_{\eta_k})$ are. Then $(     {\pi_t}  (
\Phi_{\eta_2}(a(t)),\dots,\Phi_{\eta_k}(a(t))) \big)$ is an
element of $T\tilde{P}(A)$ and we can therefore define
\[
\tilde{\pi}( \Phi_{\eta_1},\dots,\Phi_{\eta_k})(a) = \Phi_{\eta_1}
\Big ( {\pi_t}  ( \Phi_{\eta_2}(a(t)),\dots,\Phi_{\eta_k}(a(t)))
\Big ) .
\]
We recover of course the previous definition of $\tilde{\pi}$ if
$\Phi_{\eta_1}$ is regular.
\end{numrmk}

For instance, consider  the case $k=0$. If $f_t$ is a
time-dependent function on $A$, i.e., $f\in C^\infty (I\times A)$, then
Eq. (\ref{eq:tilde}) gives
\begin{equation}\label{eq:tildefunction}
\tilde{f}(a)=\int_I f\big(t,a(t)\big)dt .
\end{equation}
Eq. (\ref{eq:tildefunction}) still makes sense for any $a \in
\tilde{P}(A) $ and hence defines a function on $\tilde{P}(A)$ that
we denote by $\tilde{f}$ again. Note also that for any
time-dependent vector field $X: t \to X_t $ on $A$,
 $\tilde{X}$ is a  vector field on $\tilde{P}(A)$, which, at  any
$a(t)\in \tilde{P}(A)$, is given by
 $t \to X_{t \, |_{a(t)}}$.
%We have to explain a second abuse
%of notation for $k=1$.
%For any time-dependent vector field $X : t \mapsto X_t $
%on $A$, we had already defined (just before Propostion \ref{prop:crainic})
%a vector field on $\tilde{P}(A)$ denoted by $\tilde{X}$.
%Both definition agree
%in the sense that for any regular extension $\Phi_{\eta} $
%of some $1$-form $\eta$ on $P(A)$,
%both expressions of $ \langle \Phi_{\eta}, \tilde{X} \rangle $
%are equal.

For any time-dependent vector field $X : t \mapsto X_t$ on
$\tilde{P}(A)$, tangent to $P(A)$, we define the Lie derivative
${\cal L}_X \tilde{\pi}$ of $\tilde{\pi}$ as follows:
\begin{equation}\label{eq:Liederpi}
({\cal L}_{\tilde{X}} \tilde{\pi})( \Phi_{\eta_1},\dots,\Phi_{\eta_k}) =
 \tilde{X} (  \tilde{\pi}(\Phi_{\eta_1},
\dots,\Phi_{\eta_k}) )-\sum_{i=1}^k \tilde{\pi}(\Phi_{\eta_1} ,
\dots, {\cal L}_{\tilde{X}} \Phi_{\eta_i}, \dots, \Phi_{\eta_k}  )   .
\end{equation}
This definition needs to be justified. Of course, the Lie
derivative of a regular extension is not necessary a regular
extension, but it is  an extension and by Remark
\ref{rmk:onenotregular}, the left hand side of Eq.
(\ref{eq:Liederpi}) makes sense.

\begin{lem}\label{lem:subfol.2}
\begin{itemize}
\item[{\it i)}] For any  time-dependent function
$g: t \to  g_t $ from $I$ to $C^\infty (M)$ with $g_0=g_1=0$ ({\it
i.e., }  $g \in P \Gamma(\wedge^0 A)$),
 we have $ \tilde{G(g)}= {\mathcal F}_{\diff g}$.
\item[{\it ii)}]  If $\xi \in P_0\gm (A)$, then
 $\tilde{G(\xi)}$ is the gauge vector field
 given by  Eq.~(\ref{eq:tauinfi}).
\item[{\it iii)}] If $\xi \in P \gm (A)$, then
 $\tau_* \big( \widetilde{G(\xi)} \big)  =\Vec{\xi(0)}-\ceV{\xi(1)}$, where
$\tau: P(A) \to \gm$ is the map defined by Eq. (\ref{eq:tau}).
\item[{\it iv})]  For any  multi-differential $\delta$ on $A$, we have
\[
{\cal L}_{ \widetilde{ G(\xi) } } \tilde{\pi}_{\delta}=
  \widetilde{G(\delta(\xi))}.
\]
\item[{\it v)}]  For any  multi-differential $\delta$ on $A$
and time-dependent function $g: t \to  g_t$ on $M$, we have
$$  \tilde{\pi_{\delta}}(\widetilde{G(g)} ,\Phi_{\eta_2},\dots,\Phi_{\eta_k}   )
=   \widetilde{G(\delta (g))} ( \Phi_{\eta_2},\dots,\Phi_{\eta_k}).$$
\end{itemize}
\end{lem}
\begin{pf}
{\it i)} For any  time-dependent function $g: t \to  g_t $ from
$I$ to $C^\infty (M)$, it follows from Eq.~(\ref{complete2}) and
(\ref{eq:tildefunction}) that
\begin{equation}\label{eq:subfol.2.1}
\tilde{g^c}=\int_{0}^1 \langle \diff g_t,\rho \big(a(t)\big)
\rangle dt.
\end{equation}

Now using  Eqs.~(\ref{complete1}) and (\ref{eq:tildefunction}) we
have
\begin{equation}\label{eq:subfol.2.2}
\widetilde{\Big (\frac{\diff g}{\diff t}\Big )^v}=
 \int_{0}^1 \frac{\diff g_t}{\diff t}\big(p \smalcirc a(t)\big) dt
=-\int_{0}^1 \langle \diff g_t,\frac{\diff p \smalcirc a(t)}{\diff
t}\rangle dt,
\end{equation}
where the last equality is obtained by integration by parts using
the boundary condition $g_0=g_1=0$. Thus {\it i)} follows.

{\it ii)} follows  from the fact that for any $\xi \in \Gamma(A)$,
$ad_\xi$ is equal to $\xi^c$, a fact that can be easily checked in
local coordinates.

{\it iii)} is easily deduced from {\it (ii)} and  Eq.
(\ref{eq:taucomp''}).

{\it iv)} From the definition of  Lie derivatives
given by Eqs. (\ref{eq:Lieder}) and (\ref{eq:Liederpi}),
it follows  that the following identity holds
\begin{equation}\label{eq:Liedercomp}
{\cal L}_{\tilde{X_t}} \tilde{\pi_t} = \widetilde{[X_t,\pi_t]}
\end{equation}
for any time-dependent vector field $X: t \mapsto X_t$ on $A$ such
that $\tilde{X}$ is tangent to $P(A)$. Therefore ${\mathcal
L}_{\widetilde{G(\xi)}}\widetilde{\pi_\delta} =
\widetilde{[G(\xi),\pi_{\delta}]  }= G( \widetilde{\delta(\xi)} )$,
where the last identity follows  from Eq. (\ref{eq:3}).

{\it v)} For any time-dependent smooth function $f: I \to f_t$ on
$A$ and  any A-path $a(t)$, the differential at $a$ of
$\tilde{f_t}$ is a regular extension of  the restriction of
$\tilde{f_t} $ to $P(A) $ given by
\[
\diff  \tilde{f_{ t }}_{ |_{a} } (X)=  \int_I \langle \diff f_{ t
\, |_{a(t)} }   , X(t) \rangle dt
\]
for any $X \in T_a\tilde{P}(A)$, where $ \diff  f_{t \, |_{a(t)}}
$ is the differential of the smooth function $f_t$ at the point
$a(t)$.

By the  definition of $\widetilde{\pi_{\delta}}$,
we have
$$  \widetilde{\pi_{\delta}} (\diff \tilde{G(g)} ,\Phi_{\eta_2},
\cdots,\Phi_{\eta_k})|_{a}
 =
 \int_I   \pi_{\delta} \big( \diff G(g)_{|a(t)} ,
\Phi_{\eta_2}(t),\cdots,\Phi_{\eta_k}(t)\big)  dt     $$
 for any A-path $a(t)$.
By Eq. (\ref{eq:3}), we  have
 $$ \begin{array}{ccc}  \widetilde{\pi}_{\delta}(d\widetilde{G(g)}
  ,\Phi_{\eta_2},\cdots,\Phi_{\eta_k}   )_{|a}
%& = & \int_I G(\delta(  g))(\Phi_{\eta_2}(t),\cdots,\Phi_{\eta_k}(t)\big)  dt\\
& = & \int_I G( \delta(g))_{|_{a(t)}}
  (\Phi_{\eta_2}( t ) ,\dots,\Phi_{\eta_k}(t) \big) dt \\
 &  = &  \widetilde{G(\delta(g) )} \big( \Phi_{\eta_2},\dots,\Phi_{\eta_k}\big) .\\ \end{array}  $$
This proves {\it v)}.
\end{pf}

\begin{prop}\label{prop:lvectcanceled}
For any  $P \in P \gm (\wedge ^k  A)$, any  functions
  $f_1,\dots,f_k \in C^{\infty}(\gm  )$,
and any regular extensions $\Phi_{\diff \tau^*
f_1},\dots,\Phi_{\diff  \tau^* f_k}$ of $\tau^* f_1,\dots,\tau^*
f_k$,  we have
\begin{equation}
\label{eq:desc'}
 \tilde{G(P)}( \Phi_{\diff  \tau^*f_1},\dots, \Phi_{\diff  \tau^*f_k} ) =
 \Vec{P(1)} (\diff  f_1,\dots, \diff f_k )- \ceV{P(0)} (\diff  f_1,\dots, \diff  f_k )
 \end{equation}
\end{prop}

First we need the following

\begin{lem}
\label{lem:geeta} For any  $\xi \in P\Gamma(A)$, $a \in P(A)$, any
covector  $\eta  \in T^*_aP(A)$ conormal to the gauge orbit,
 and any its  regular extension  $\Phi_{\eta}$,
 the following identity holds
 \begin{equation}
\label{eq:geeta} \frac{\diff \langle  \Phi_{\eta} (a(t)),
\xi^v(t)\rangle}{\diff t}=\langle \Phi_\eta (t) ,{G(\xi)} (t)\rangle ,
\end{equation}
In particular,  for any   smooth function  $f : \Gamma \to {\mathbb
R}$ and any regular extension $\Phi_{d \tau^* f}$ of
 $ \tau^* f$, the following identity holds
\begin{equation} \label{eq:weak2}
\frac{\diff \langle\Phi_{\diff  \tau^*
f}(t),\xi^v(t)\rangle}{\diff t}= \langle \Phi_{\diff  \tau^*
f}(t),G( \xi ) (t) \rangle .
\end{equation}
  \end{lem}
\begin{pf}
Since  $\eta$  is conormal  to the gauge orbits, for any smooth
map $\chi(t): I\to \Gamma(A)$   with $\chi(0)=\chi(1)=0$, we have
$$  \langle\eta ,\widetilde{G\big( \chi \big)}_{|_a}\rangle =
 \langle\Phi_{\eta} ,\widetilde{G\big( \chi \big)}_{|_a}\rangle =   \int_I
\langle\Phi_{\eta}(t) ,G\big( \chi(t) \big)_{|_{a(t)}}\rangle dt
=0.$$ Applying this equality to $\chi(t):=\psi(t)\xi(t) $, where
$\psi(t)$ is any smooth function with $\psi(0)=\psi(1)=0$ and
using Eq.~(\ref{eq:fpsi}) we obtain

\begin{equation}\label{eq:geeta2}
0= \int_{0}^1 \psi(t) \langle\Phi_{\eta} (t),G\big(\xi (t)\big)
\rangle dt + \int_{0}^1 \frac{\diff \psi(t)}{\diff t}
\langle\Phi_{\eta} (t),\xi^v(t)\rangle dt .
\end{equation}
The result follows by integration by part.
\end{pf}

Now we are ready to prove Proposition \ref{prop:lvectcanceled}.
\begin{pf}
According to  Lemma \ref{lem:subfol.2}, we have
 $\widetilde{G(\xi)} (   \tau^* f)= \Vec{\xi(1)}( f)
 - \ceV{\xi(0)} ( f)$, $\forall f\in C^{\infty}(\gm )$. On the
other hand, we have, by definition,
  $\widetilde{G(\xi)}  (
\tau^* f ) = \int_{0}^1 \langle\Phi_{\diff \tau^* f}(t),G( \xi)(t) \rangle dt$.
 Therefore, we have
  $$  \int_{0}^1 \langle\Phi_{\diff  \tau^* f}(t),
G(\xi) (t)\rangle dt
 = \Vec{\xi(1)} (f) -\ceV{ \xi(0)} ( f).   $$
By Lemma \ref{lem:geeta}, it follows   that
\[
\int_{0}^1 \frac{\diff }{\diff t}(\langle   \Phi_{\diff  \tau^*
f}(t) ,\xi^v(t)  \rangle) dt = \Vec{\xi(1)}  ( f) -\ceV{ \xi(0)} (
f).
\]
Hence
$$   \langle\Phi_{\diff \tau^* f}(1),\xi^v(1)\rangle-
 \langle\Phi_{\diff \tau^* f}(0),\xi^v(0)\rangle =
  \Vec{\xi(1)}  (f) - \ceV{\xi(0)} ( f)  .$$
%Hence, we have, for any $P \in P\gm (\wedge ^\bullet A)$,
Since this identity holds for any time-dependent section $\xi$, we
have
\[
\langle\Phi_{\diff \tau^* f} (1),\xi^v(1) \rangle= \Vec{\xi(1)} (
f) \mbox{ and } \langle \Phi_{\diff \tau^* f}(0),\xi^v(0) \rangle
= \ceV{\xi(0)} ( f) .
\]
This implies that   for any $P\in P\gm (\wedge^k A)$, we have
\begin{equation}\label{eq:rightleft}
\left\{
\begin{array}{ccc}
P^v(1)( \Phi_{ \diff\tau^* f_1}(1), \cdots, \Phi_{\diff\tau^* f_k}(1))&=&\Vec{P(1)}(f_1,\dots,f_k)  \\   $$
   P^v(0)(\Phi_{\diff \tau^* f_1}(0), \dots,  \Phi_{\diff \tau^* f_k}(0))
&=&\ceV{P(0)}(f_1,\dots,f_k).  \\
\end{array}
\right.
\end{equation}

Now assume that $P(t) =\xi_1(t) \wedge \ldots \wedge \xi_k (t)$
for some time-dependent sections  $\xi_1,\dots,\xi_k$ of
$\Gamma(A)$. Then, according to Eq. (\ref{eq:2}),
\[
G(P)(t)= \sum_{i=1}^k \xi_1^v(t) \wedge \dots \wedge G\big( \xi_i
(t) \big)\wedge \dots \wedge  \xi_k ^v(t).
\]

Therefore, for any $a(t) \in P(A)$,
\be
&&  \tilde{G(P)}(\Phi_{\diff  \tau^* f_1},\dots,\Phi_{\diff  \tau^* f_k} ) (a(t)) \\
&=&\int_0^1 G(P)(t) (\Phi_{\diff  \tau^* f_1}(t),\dots, \Phi_{\diff  \tau^* f_k}(t))dt\\
&=& \sum_{i=1}^{k} \int_0^1 (\xi^v_1(t) \wedge \ldots \wedge
G\big(\xi_i(t)\big)_{|_{a(t)}} \wedge
 \ldots \wedge  {\xi_k ^v(t)})
\big( \Phi_{\diff  \tau^* f_1}(t),\ldots,\Phi_{\diff  \tau^* f_k} (t)\big) dt\\
 &=& \sum_{\sigma \in S_k} \sum_{i=1}^{k}(-1)^{|\sigma |} \int_{0}^1
\langle   \Phi_{\diff  \tau^* f_1}(t),\xi^v_{\sigma (1)}(t)\rangle
\ldots
  \langle \Phi_{\diff  \tau^* f_1}(t),G\big(\xi_{\sigma (i)}(t)\big)\rangle  \dots
  \langle  \Phi_{\diff  \tau^* f_k}(t),\xi ^v_{\sigma (k)}(t)\rangle  dt    \\
  &&\mbox{ (By Lemma  \ref{lem:geeta}) }\\
&=& \sum_{\sigma \in S_k} \sum_{i=1}^{k} (-1)^{|\sigma |}
\int_{0}^1 \langle \Phi_{\diff  \tau^* f_1}(t),\xi^v_{\sigma
(1)}(t)\rangle \dots \frac{\diff \langle  \Phi_{\diff  \tau^* f_i}
(t), {\xi_{\sigma (i)}(t)}\rangle}{\diff  t} \ldots
    \langle \Phi_{\diff  \tau^* f_k}(t),\xi^v_{\sigma (k)}(t) \rangle  dt \label{eq:bornes}\\
&=&  \sum_{\sigma \in S_k} (-1)^{\sigma} \int_{0}^1
  \frac{\diff }{ \diff  t}\big(\langle  \Phi_{\diff  \tau^* f_1}(t),\xi^v_{\sigma (1)}(t)\rangle \ldots
  \langle \Phi_{\diff \tau^* f_k}(t),\xi^v_{\sigma (k) }(t)\rangle \big)dt \\
&=&  P^v(t)( \Phi_{\diff  \tau^* f_1}(t),\dots, \Phi_{\diff \tau^*
f_k}(t))|_{0}^{t=1} \ee The result now follows from Eq.
(\ref{eq:rightleft}).
\end{pf}

\subsection{From $k$-vector  fields
on $\tilde{P}(A)$ to $k$-vector fields on $\Gamma$}

Assume that $\gm$ is an $\alpha$-simply connected and
$\alpha$-connected Lie groupoid with Lie algebroid $A$. Let
$\delta$ be a $k$-differential on $A$ and $\tilde{\pi}_{\delta}$ the
 corresponding  $k$-vector field  on $\tilde{P}(A)$.
The goal of this section is to construct  a $k$-vector field
$\Pi_{\delta}$ on $\Gamma$ from $\tilde{\pi}_{\delta}$.

\begin{prop}\label{prop:down}
Let $f_1, \ldots , f_k\in C^{\infty}(\gm )$ be a family of smooth
functions on $\Gamma$ and $\Phi_{\diff \tau^*
f_1},\ldots,\Phi_{\diff \tau^* f_k}$ regular extensions of
$\tau^*f_1, \ldots , \tau^* f_k$. Then $\tilde{\pi}_\delta
(\Phi_{\diff \tau^* f_1},\ldots,\Phi_{\diff \tau^* f_k})$ is a
smooth function on $P(A)$,  which is
\begin{itemize}
\item [{\it i)}] independent of the choice of the chosen regular
extensions $ \Phi_{\diff  \tau^* f_1},\dots,\Phi_{\diff \tau^*
f_k}$ of $\diff  \tau^* f_1, \dots,\diff  \tau^* f_k$, and
\item[{\it ii)}] invariant under gauge transformations.
\end{itemize}
\end{prop}
\begin{pf}{\it i)}
It suffices to prove that if $\Phi_{\diff \tau^* {f_1}}$ is a
regular extension of zero,
 the function $\tilde{\pi}_\delta
(\Phi_{\diff \tau^*f_1},\dots,\Phi_{\diff \tau^*{f_k}})$ vanishes.
By Lemma \ref{th:recap4.1} {\em ii)}, for any $a \in P(A)$, there
exists $g :I\to C^\infty (M)$ with $g_0=g_1=0$ such that
${\Phi_{\diff \tau^* {f_1}}}_{|_{a}}=\diff {\mathcal F}_{\diff g
\, |_{a}} $. Then \be && \tilde{\pi}_\delta (\Phi_{\diff
\tau^*f_1},\dots,\Phi_{\diff  \tau^*{f_k}})_{ |_a}  =
\tilde{\pi}_{\delta}(\diff {\mathcal F}_{\diff g}, \Phi_{\diff
\tau^*f_2},\dots,\Phi_{\diff \tau^*{f_k}}). \ee
 According to Lemma
\ref{lem:subfol.2} {\it i)}, we have ${\mathcal F}_{\diff
g}=\tilde{G(g)}$. Lemma \ref{lem:subfol.2} {\it v)} implies that
\[
 \tilde{\pi}_\delta (\Phi_{\diff \tau^*f_1},\dots,\Phi_{\diff \tau^*{f_k}})_{|_a}=
 \tilde{G(\delta(g))}(   \Phi_{\diff \tau^*f_2},\dots,\Phi_{\diff \tau^*{f_k}}  )_{|_a}.
\]
The latter vanishes according to  Proposition
\ref{prop:lvectcanceled},  since by assumption $\xi(0)=\xi(1)=0$.
This proves {\it i)}.

Before starting the proof of {\it ii)}, we have to add a comment.
We have proven in {\it (i)}  that for any regular extension $\Phi$
of $0$, and any regular extensions $\Phi_{\diff
\tau^*f_2},\dots,\Phi_{\diff  \tau^*{f_k}}$ of $\tau^* f_2,\dots,
\tau^* f_n$, we have
$$ \tilde{\pi}_\delta (\Phi, \Phi_{\diff \tau^*f_2},\dots,\Phi_{\diff \tau^*{f_k}})
=    0$$

Regular extensions are dense with respect to the induced topology
of $T^*_a\tilde{P}(A)$. Thus,
\begin{equation}\label{eq:onenon}
\tilde{\pi}_{\delta} (\Phi, \Phi_{\diff
\tau^*f_2},\dots,\Phi_{\diff \tau^*{f_k}}) =    0,
\end{equation}
for any extension $\Phi$ of 0.

{\it ii)}  Since the functions $\tau^* f_1, \cdots,\tau^* f_k$ are
invariant under the gauge transformation, $\forall \xi \in P_0 \gm
(A)$,  the Lie derivative ${\cal L}_{\widetilde{G(\xi)}} \Phi_{
\diff \tau^* f_i   }  $ is,  for all $i \in \{ 1,\dots,k \}$, an
extension of zero. Therefore, by Eq. (\ref{eq:onenon}), for all $i
\in \{1,\dots,k  \}$, $  \tilde{\pi}_\delta (\Phi_{\diff
\tau^*f_1},\dots,{\cal L}_{\widetilde{G(\xi)} } \Phi_{\diff
\tau^*{f_i}}, \dots \Phi_{\diff \tau^*{f_k}})= 0    $. By
% definition of the Lie derivative of $\tilde{\pi}_{\delta}$ with respect to a vector field given in
Eq. (\ref{eq:Liederpi}), we have
$$ \tilde{G(\xi)}   (\tilde{\pi}_\delta ( \Phi_{\diff \tau^*f_1},\dots, \Phi_{\diff \tau^*{f_k}} ))
 = ({\cal L}_{\tilde{G(\xi)}}  \tilde{\pi}_\delta) (  \Phi_{\diff \tau^*f_1},\dots,
  \Phi _{\diff \tau^* f_k }   ).$$
By Lemma \ref{lem:subfol.2} {\it v)}, we have
\[
\tilde{G(\xi)}   (\tilde{\pi}_\delta ( \Phi_{\diff
\tau^*f_1},\dots, \Phi_{\diff \tau^*{f_k}} ))    =
\tilde{G(\delta(\xi))} (   \Phi_{\diff
\tau^*f_1},\dots,\Phi_{\diff \tau^*{f_k}}  ),
\]
which is identically equal to $0$ according to Proposition
\ref{prop:lvectcanceled}.
 This proves {\it ii)}.
\end{pf}

By Proposition \ref{prop:crainic} and \ref{prop:down} {\em ii)}
above, the function $\tilde{\pi}_\delta (  \Phi_{\diff
\tau^*f_1},\dots,   \Phi_{\diff \tau^*f_k})$ descends to a smooth
function  on $\gm$, which will be denoted by
$\{f_1,\cdots,f_k\}$.
 I.e.,
\begin{equation}\label{eq:godoeq} \tilde{\pi}_{\delta}
\big( \Phi_{\diff \tau^*f_1},\dots,\Phi_{\diff \tau^*{f_k}} \big)
           =
\tau ^* \{f_1,\dots,f_k\} .
\end{equation}
It is straightforward to check that the map $f_1,\dots,f_k \to
\{f_1,\dots,f_k\}$ indeed defines a $k$-vector
field $\Pi_{\delta} $ on $\Gamma$, i.e.,
\begin{equation}
\label{eq:58A}
\{f_1,\dots,f_k\}=\Pi_\delta (df_1, \cdots , df_k ).
\end{equation}

\begin{prop}
$\Pi_{\delta}$ is a multiplicative $k$-vector field on $\gm$.
\end{prop}
\begin{pf}
By definition, $\Pi_{\delta} $ is given, for any
$\eta_1,\dots,\eta_k \in T^*_g\Gamma$ by
\begin{equation}\label{eq:deffpi}
\Pi_\delta(\eta_1,\dots,\eta_k)=\int_I (\pi_{\delta})_{|_{a(t)}}
\big( \Phi_{\tau^* \eta_1 }(t), \dots , \Phi_{\tau^* \eta_k }
(t)\big) dt \end{equation} where $\Phi_{\tau^* \eta_1 }(t), \dots
, \Phi_{\tau^* \eta_k } (t)$ are any regular extensions of $
\tau^* \eta_1 , \dots, \tau^* \eta_k$, with a smooth dependence on
the variable $t$ and $a(t)$ is any $A$-path with $\tau(a)=g$.
Since smooth functions are dense in the space of piecewise
continuous functions with finitely many discontinuities, we obtain
that Eq. (\ref{eq:deffpi}) remains valid when the extensions
$\Phi_{\tau^* \eta_1 }(t), \dots , \Phi_{\tau^* \eta_k } (t) $ are
just assumed to be piecewise continuous in $t$ (with finitely many
discontinuities).

It is straightforward to check that for any $g \in \Gamma$ with
$\alpha (g)=m$ and $\beta(g)=n$, there exists an $A$-path $a(t)$
with $\tau(a)=g$ such that $a(t)$ is constantly equal to $m$ in a
neighborhood of $t=0$ and constantly equal to $n$ in a
neighborhood of $t=1$. Consider two composable elements $g_1,g_2
\in \Gamma$ and two composable elements $\eta_1\in T^*_{g_1}
\Gamma$ and $\eta_2 \in T^*_{g_2} \Gamma$. Choose now $A$-paths
$a_1(t) $ and $a_2(t)$ satisfying the previous condition. Define
$a(t)$ by $a(t)=a_1(2t)$ for $t\in [0,\frac{1}{2}] $ and by $a(t)
= a_2(2t -1) $  for $t\in [\frac{1}{2},1] $. By construction,
$a(t) $ is an $A$-path and $\tau(a) =g_1 g_2$.

Consider  $ \Phi_{\tau^* \eta_1}(t)$ and $ \Phi_{\tau^*
\eta_1}(t)$ two regular extensions of $\tau^* \eta_1 $ and $\tau^*
\eta_2 $. Then the map defined by $  \Phi(t) =  \Phi_{\tau^*
\eta_1}(2t) $ for $t \in [0,\frac{1}{2}] $  and $\Phi(t)=
\Phi_{\tau^* \eta_2}(2t-1) $ for $t \in [\frac{1}{2},1] $ is a
regular extension of $ \Phi_{\eta_1 \eta_2} \in T^*_{g_1 g_2}
\Gamma$, and it may have a point of discontinuity in $t =
\frac{1}{2}$.

We choose now $k$ compatible pairs $\eta_1^i, \eta_2^i $ for $
i=1,\dots, k$ of elements of $ T^*_{g_1} \Gamma$ and  $ T^*_{g_2}
\Gamma$. For all $i=1,\dots, k$ we consider two regular extensions
$\Phi_{\tau^* \eta^i_1} $ and $\Phi_{\tau^* \eta^i_2} $ of $\tau^*
\eta_1^i $ and $\tau^* \eta_2^i $ respectively. And, for all $i
=1,\dots, k$ again, we form $\Phi_{\eta_1^i \eta_2^i}(t) $ as
above.

Eq. (\ref{eq:deffpi}) being valid even for piecewise regular
extension, the first of the identities below is valid; the other
ones are routine.
\begin{equation}
\begin{array}{ccc} \Pi_\delta(\eta_1^1 \eta_2^1 ,\dots,\eta_1^k \eta_2^k )
&=&\int_I (\pi_{\delta})_{|_{a(t)}} \big( \Phi_{\eta_1^1 \eta_2^1}(t), \dots ,
\Phi_{\eta_1^k \eta_2^k } (t)\big) dt     \\
& =&  \int_{t=0}^{\frac{1}{2}}   (\pi_{\delta})_{|_{a_1(t)}} \big(
\Phi_{ \eta_1^1 \eta_2^1}(\frac{t}{2}), \dots ,
\Phi_{ \eta_1^k \eta_2^k} (\frac{t}{2})\big)  dt \\
 &+& \int_{t=\frac{1}{2}}^1
(\pi_{\delta})_{|_{a_2(t)}}  \big( \Phi_{\eta_1^1 \eta_2^1 }(\frac{1+t}{2}), \dots ,
\Phi_{ \eta_1^k \eta_2^k} (\frac{1+t}{2})\big)     dt  \\
 & =&  \int_I (\pi_{\delta})_{|_{a_1(t)}}    \big( \Phi_{\tau^* \eta_1^1 }(t), \dots ,
\Phi_{\tau^* \eta_1^k } (t)\big)  dt \\
 &+& \int_I
(\pi_{\delta})_{|_{a_2(t)}}  \big( \Phi_{\tau^* \eta_2^k  }(t), \dots ,
\Phi_{\tau^* \eta_2^k } (t)\big)     dt  \\
 &= &   \Pi_\delta(\eta_1^1  ,\dots,\eta_1^k  )   + \Pi_\delta(\eta_2^1  ,\dots,\eta_2^k  )  \\
\end{array}
\end{equation}

By Proposition \ref{prop:1cocycle}{\it (ii)}, $ \Pi_\delta$ is multiplicative.
\end{pf}
%Moreover, since  our construction preserves the graded bracket in each step,  the
%following result is immediate.
%\begin{prop}
%The map $\delta \to \Pi_{\delta}$ from ${\mathcal A}$ to $\oplus _k\mathfrak
%X ^k _{mult }(\gm )$ is a homomorphism of graded Lie algebras.
%\end{prop}

Finally to complete  the proof of Theorem \ref{derivation},
 we need to show that the map
$\delta \to \Pi_{\delta}$  constructed above  is indeed the
inverse of $\Pi \to \delta_{\Pi}$. This is due to the following

\begin{prop}\label{prop:overway}
For any $k$-differential $\delta$, the identities \be
&& [\Pi_{\delta}, \Vec{X}]=\Vec{\delta X}\mbox{ and }\\
 &&[\Pi_{\delta}, \alp ^\ast f]=\Vec{\delta f}
\ee
hold.
\end{prop}
\begin{pf}
For any functions $f_1,\dots,f_k \in C^{\infty}(\gm)$ and any $ X
\in \Gamma(A)$, one has
\begin{equation} \label{eq:Pidelta}  [\Pi_{\delta}, \Vec{X}](f_1,\cdots,
 f_k)=
\Vec{X}  \big(\Pi_{\delta}(f_1,\cdots,f_k)\big)- \sum_{i=1}^k
\Pi_{\delta}(f_1,\cdots, \Vec{X}( f_i), \cdots,f_k)
.\end{equation}

Let $\xi \in P\gm (A)$ be an element such that  $\xi(1)=X$ and
$\xi(0)=0$. According to Lemma \ref{lem:subfol.2}, we have  $\tau_*
(\widetilde{G(\xi)} ) = \Vec{X} $. Therefore, for any regular
extension $\Phi_{\diff \tau^* f} $ of $\tau^* f$, where
$f \in C^{\infty}(\Gamma)$,
${\cal L}_{\tilde{G(\xi)}}   \Phi_{\diff \tau^* f} $ is a regular
extension of $  \tau^* (\Vec{X} (f))$.

Let us choose some regular extensions $\Phi_{\diff \tau^* f_1
},\dots, \Phi_{\diff \tau^* f_k }$ of $\tau^*f_1,\cdots, \tau^* f_k$.
Applying $\tau^*$ to both sides of Eq. (\ref{eq:Pidelta}) and using Eq.
(\ref{eq:58A}), we obtain
\[
\begin{array}{rcl}
\tau^* ([\Pi_{\delta}, \Vec{X}](f_1,\cdots , f_k))   &=&
\widetilde{G(\xi)} \cdot \tilde{\pi}_{\delta}(\Phi_{\diff \tau^*
f_1},\ldots,\Phi_{\diff \tau^* f_k}) \\ & &\displaystyle -
\sum_{i=1}^k \tilde{\pi}_{\delta}(\Phi_{\diff \tau^* f_1},\dots,
{\cal L}_{\tilde{G(\xi)}} \Phi_{\diff \tau^* f_i} ,
\ldots,\Phi_{\diff \tau^* f_k}).
\end{array}
\]

Using the  definition as given by  Eq. (\ref{eq:Liederpi}),
%of the Lie derivative of $\tilde{\pi_{\delta}}$ with
%respect to a vector field (as defined in Eq. (\ref{eq:Liederpi}))
the right hand side of the  equation above is
$$ ({\cal L}_{\widetilde{G(\xi)}} \tilde{\pi}_{\xi}) (\Phi_{\diff \tau^* f_1},\dots,  \Phi_{\diff
\tau^* f_k})    $$
 which is again, by Lemma \ref{lem:subfol.2} {\em (iv)},
equal to $ \tilde{ G \big( \delta(\xi) \big)} (  \Phi_{\diff
\tau^* f_1},\dots,  \Phi_{\diff \tau^* f_k}) $. By Proposition
\ref{prop:lvectcanceled}, the latter  is equal to $\tau^*
\Vec{\delta(X)} (f_1,\dots,f_k)  $. This proves the first
equality. The proof of the other equality is similar.
\end{pf}

%\subsection{Examples and applications}
%\begin{numex}
%{\rm For the particular case of Lie groups, we recover the
%one-to-one correspondence between ad-1-cocyles and multiplicative
%multivector fields (see \cite{LuW:1990}). }
%\end{numex}
%\begin{numex}
%{\rm If $\phi \in \gm (A^\ast )$ is a 1-cocycle on the Lie
%algebroid cohomology of $A$ and $\delta _\phi$ is the
%corresponding 0-differential then, applying our theorem, we get a
%multiplicative 0-vector field, i.e., a multiplicative function.
%Conversely, if $\sigma $ is a multiplicative function on $\Gamma$
%then we obtain that $\diff \sigma _{|A}$ is a 1-cocycle on the Lie
%algebroid $A$. }
%\end{numex}
%\begin{numex}
%{\rm For the case of $k=1$, we recover the one-to-one
%correspondence between covariant differential operators on $A$
%which are a derivation with respect to the bracket and
%multiplicative vector fields on the Lie groupoid (see
%\cite{MackenzieX:1998,Moerdijk-Mrcun}).}
%\end{numex}
%\begin{numex}
%{\rm If $P\in \gm (\wedge ^kA)$ is a $k$-section and $ad(P)=\lcf
%P,\cdot \rcf$ is the corresponding coboundary $k$-differential
%%then we get that the corresponding multiplicative $k$-vector field
%is  $\Pi =\Vec{P}-\ceV{P}$, recovering the coboundary
%multiplicative multivector fields in the particular case when
%$\gm$ is a Lie group $G$. }
%\end{numex}
%

\section{Quasi-Poisson groupoids}
In this section, we will introduce the notion of  quasi-Poisson
groupoids. We  describe their  infinitesimal invariants, namely
quasi-Lie bialgebroids, and study the corresponding momentum map
theory.
\subsection{Definition and properties}
\begin{defn}
A {\em quasi-Poisson groupoid} is a triple $(\gm\gpd M ,\Pi
,\Omega )$, where $\gm\gpd M$ is a Lie groupoid, $\Pi$ is a
multiplicative bivector field on $\gm$ and $\Omega \in \gm (\wedge
^3 A)$ such that the following compatibility conditions hold
\begin{equation}\label{quasiPoisson1}
\half [\Pi ,\Pi ]=\Vec{\Omega}-\ceV{\Omega}, \ \ \mbox{ and}
\end{equation}
\begin{equation}\label{quasiPoisson2}
[ \Pi ,\Vec{\Omega} ]=0 .
\end{equation}
\end{defn}
Using Proposition \ref{base} and Eqs. (\ref{quasiPoisson1}) and
(\ref{quasiPoisson2}), we obtain the following

\begin{prop}\label{bivector-onthebase}
Given a quasi-Poisson groupoid $(\gm \gpd M, \Pi ,\Omega )$ there
exists a bivector field $\Pi _M$ on $M$ such that
\[
\Pi _M =\alp _\ast \Pi =-\bet _\ast \Pi.
\]
$\Pi _M$ satisfies the relations
\begin{equation}\label{compat-base}
\half [\Pi _M,\Pi _M ] = \Omega _M, \ \ \mbox{ and}
\end{equation}
\begin{equation}\label{base-quasi-Poisson}
[\Pi _M ,\Omega _M]=0,
\end{equation}
where $\Omega _M$ is the 3-vector field $\Omega _M=\rho (\Omega )$,
and $\rho :\gm (\wedge ^3A)\to \mathfrak X ^3(M)$ is the extension
of the anchor map.
\end{prop}
Some interesting examples of quasi-Poisson groupoids are listed
below.
\begin{numex}
If $\Omega =0$ then we just have a multiplicative Poisson
structure $\Pi$ on a Lie groupoid $\gm \gpd M$. I.e., $(\gm\gpd M,
\Pi )$ is a Poisson groupoid.
\end{numex}
\begin{numex}
{\rm If $G$ is a Lie group, then we recover the notion of
quasi-Poisson structures on a Lie group of Kosmann-Schwarzbach
\cite{Kosmann:1991}. That is, a multiplicative bivector field
$\Pi$ on $G$ and an element $\Omega \in \wedge ^3 \mathfrak g$
such that $ \half [\Pi ,\Pi ]=\Vec{\Omega}-\ceV{\Omega}$ and $[\Pi
,\Vec{\Omega}]=0$. }
\end{numex}
%\begin{numex}
%{\rm Let $(M, \pi , \phi )$ be a twisted-Poisson manifold, that
%is, a bivector $\pi \in \mathfrak X^2(M)$ and a closed 3-form
%$\phi\in \Omega ^3(M)$ such that $\half [\pi,\pi]=(\wedge ^3\pi
%^\sharp )(\phi)$. It is straightforward to see, from Example
%\ref{coarse-example}, that $\Pi = \pi \oplus -\pi\in \mathfrak X
%^2(M\times M)$ is a multiplicative bivector field on the pair
%groupoid $\Gamma =M\times M$. Moreover, we directly deduce that
%\[
%\half [\Pi ,\Pi ]=\Vec{\Omega}-\ceV{\Omega}
%\]
%where $\Omega=(\wedge ^3\pi ^\sharp)(\phi )\in \mathfrak X ^3
%(M)\cong \gm (\wedge ^3 TM)$. On the other hand, the graded Jacobi
%identity and Eq. (\ref{Eq-twisted}) imply that
%\[
%0=[\pi,[\pi ,\pi ]]=[\pi ,(\wedge ^3 \pi ^\sharp )(\phi )].
%\]
%As a consequence, $[\Pi ,\Vec{\Omega}]=0$. Thus, $(\Gamma \gpd
%M,\Pi ,\Omega )$ is a quasi-Poisson groupoid. }
%\end{numex}
\begin{prop}\label{non-degenerate-quasi}
Let $(\gm \gpd M, \Pi, \Omega )$ be a quasi-Poisson groupoid such
that $\Pi\in \mathfrak X^2 (\gm)$ is non-degenerate. Let $\omega
\in \Omega ^2(\gm)$ be its corresponding non-degenerate 2-form and
$\phi \in \Omega ^3 (M)$ be the 3-form on $M$ defined by $(\wedge
^3 \omega ^\flat )(\Vec{\Omega})=\alpha ^\ast \phi $. Then $(\gm
\gpd M,\omega ,\phi )$ is a non-degenerate twisted symplectic
groupoid in the sense of \cite{CX}. That is,
\begin{enumerate}
\item $d\phi =0$;
\item $d\omega =\alp ^\ast \phi -\bet ^\ast \phi$; and
\item $\omega$ is multiplicative, i.e., the 2-form $(\omega,
\omega ,-\omega )$ vanishes when being restricted to the graph of
the groupoid multiplication $\Lambda \subset \gm\times\gm\times
\gm$.
\end{enumerate}
\end{prop}
\begin{pf}
Since $\Pi$ is multiplicative, it follows from Remark
\ref{multiplicative-bivectors} that $\omega$ satisfies the formula
$m^\ast \omega =pr_1^\ast \omega +pr_2^\ast \omega$. That is,
$\omega$ is multiplicative.

Since $\omega$ is multiplicative, we know (see \cite{CX}) that the
Lie algebroid $A$ of $\gm$ is isomorphic to $T^\ast M$ as a vector
bundle and the isomorphism $\lambda: A\to T^\ast M$ is
characterized by $\omega ^\flat (\Vec{X})=\alp ^\ast \eta$, for
$X\in \Gamma (A)$ and $\eta \in \Omega ^1 (M)$. In general, we
have
\[
(\wedge ^k \omega ^\flat )(\Vec{P})=\alp ^\ast \varphi, \qquad
(\wedge ^k \omega ^\flat )(\ceV{P})=\bet ^\ast \varphi,
\]
$\forall P\in \Gamma (\wedge ^k A)$ and $\varphi \in \Omega ^k (M)$.

Define $\phi \in \Omega ^3(M)$ as the 3-form on $M$ such that
$(\wedge ^3 \omega ^\flat)(\Vec{\Omega})=\alp ^\ast \phi$ or,
equivalently, $\Vec{\Omega}=(\wedge ^3 \Pi  ^\sharp )(\alp ^\ast
\phi)$. Using that $\Pi ^\sharp$ is the inverse of $\omega
^\flat$,  $(\wedge ^3 \Pi ^\sharp )(d\omega )=\half [\Pi ,\Pi
]$ and Eq. (\ref{quasiPoisson1}),
 we deduce that
\[
\begin{array}{rcl}
(\wedge ^3 \Pi ^\sharp )(d\omega )&=&\half [\Pi ,\Pi ] =
\Vec{\Omega}-\ceV{\Omega} \\&=& (\wedge ^3 \Pi  ^\sharp )(\alp
^\ast \phi -\bet ^\ast \phi).
\end{array}
\]
As a consequence, we have  $d\omega =\alp ^\ast \phi -\bet ^\ast \phi$.

Let $\Pi _M=\alp _\ast \Pi$. We will prove that
\begin{equation}\label{relation-quasi}
\lambda (-\delta _\Pi P)=\delta (\lambda (P)),\mbox{ for }P\in \gm
(\wedge ^\bullet A),
\end{equation}
where $\delta _\Pi$ is the 2-differential corresponding
 to $\Pi$ (see Theorem \ref{derivation}),
 and $\delta :\Omega ^\bullet (M)\to
\Omega ^{\bullet +1}(M)$ is the map characterized by
\[
\begin{array}{l}
\delta f=df, \ \forall f\in C^\infty (M),\\ \delta \eta=d\eta
-\Pi ^\sharp _M(\eta) \per \phi, \forall \eta\in \Omega ^1 (M).
\end{array}
\]
From Lemma 2.4 in \cite{CX}, we know that $\delta f=df=-\lambda
(\delta _\Pi f)$. If $\eta \in \Omega ^1 (M)$ and $X\in \Gamma
(A)$ such that $\lambda (X)=\eta$ (that is, $\Vec{X}=\Pi ^\sharp
(\alp ^\ast \eta)$),  then using the relation

$$(\wedge ^2\Pi ^\sharp)d\gamma  =[\Pi ^\sharp (\gamma), \Pi ]
-\half (\gamma \per [\Pi ,\Pi]),  \ \ \ \forall \gamma\in \Omega^1 (M)$$
  and Eqs. (\ref{how-to-go-down})
and (\ref{quasiPoisson1}), we obtain
\[
\begin{array}{rcl}
(\wedge ^2\Pi ^\sharp )d\alp ^\ast \eta &=& [\Pi ^\sharp (\alp
^\ast \eta ), \Pi ]- (\alp ^\ast \eta \per
(\Vec{\Omega}-\ceV{\Omega}) )\\&=& -\Vec{\delta _\Pi X} - (\alp
^\ast \eta \per  ((\wedge ^3 \Pi  ^\sharp )(\alp ^\ast \phi - \bet
^\ast \phi))) .
\end{array}
\]
Applying $\wedge ^2\omega ^\flat$ to both sides of this equation
and using that
\[
(\wedge ^2\omega ^\flat )(\alp ^\ast \eta \per ((\wedge ^3 \Pi
^\sharp )(\alp ^\ast \phi - \bet ^\ast \phi)))=-\alp ^\ast (\Pi
^\sharp _M(\eta )\per \phi ),
\]
we conclude that
\[
- (\wedge ^2\omega ^\flat )(\Vec{\delta _\Pi X})=\alp ^\ast (d\eta
-\Pi ^\sharp _M(\eta) \per \phi ).
\]
I.e., $- \lambda (\delta _\Pi X )=\delta \eta$. As a consequence,
by  Eqs. (\ref{quasiPoisson2}), (\ref{relation-quasi}) and that
$\lambda(\Omega)=\phi$, we  have  $\delta \phi=0$. Thus
$d\phi=0$.
Thus, we conclude that $(\gm \gpd M, \omega ,\phi )$ is a
non-degenerate twisted symplectic groupoid.
\end{pf}

\subsection{Quasi-Poisson groupoids and quasi-Lie bialgebroids}

In this subsection, we will describe the infinitesimal invariants of
quasi-Poisson groupoids,  i.e., quasi-Lie bialgebroids.
The notion of quasi-Lie bialgebroids was first introduced by Roytenberg
\cite{Roy}. Here, we give an alternative definition using
2-differentials.
\begin{defn}\label{quasi-Lie-bialgebroid}
{\rm A {\em quasi-Lie bialgebroid} corresponds to a 2-differential
whose square is a coboundary, i.e., $\delta :\gm (\wedge ^\bullet
A)\to \gm (\wedge ^{\bullet +1}A)$ such that $\delta \smalcirc
\delta =\lcf \Omega , \cdot \rcf $ for some $\Omega \in \gm
(\wedge^3 A )$ satisfying $\delta \Omega =0$.}
\end{defn}

An interesting example of a quasi-Lie bialgebroid is the
following.
\begin{numex}\label{twisted-Poisson}
{\rm Recall that a {\em twisted Poisson structure} $(M, \pi , \phi
)$ is a bivector field $\pi\in \mathfrak X ^2(M)$ and a closed
3-form $\phi \in \Omega ^3(M)$ such that
\begin{equation}\label{Eq-twisted}
\half [\pi ,\pi ]=(\wedge ^3\pi ^\sharp)(\phi ).
\end{equation}
In this situation, let $A=T^* M$ be the corresponding twisted Lie
algebroid (see Eq. (\ref{tLie}) and, for more details, \cite{SW}).
Define $C^{\infty}(M)\stackrel{\delta _{\pi,\phi}}{\lon}\Omega^1
(M) \stackrel{\delta _{\pi,\phi}}{\lon}\Omega^2 (M)$ as follows.
On $C^{\infty}(M)$, $\delta _{\pi,\phi}$ is  the usual de~Rham
differential, while on $\Omega^1 (M) $, $\delta _{\pi,\phi} \eta
=d\eta  - \pi^\sharp (\eta) \per \,\phi$,  for all $\eta \in
\Omega^1 (M) $. Then $(T^\ast M,\delta _{\pi,\phi},\phi )$ is a
quasi-Lie bialgebroid \cite{CX}. }
\end{numex}

A direct consequence of Lemma \ref{diff->bracket} is the following

\begin{prop}\label{4.9}
Let $(A,\delta ,\Omega )$ be a quasi-Lie bialgebroid. Then,
$\delta$ induces a bivector field $\pi _M$ on $M$ such that
\[
\begin{array}{l}
\half [\pi _M, \pi _M] =\Omega_M ,\\[5pt]
[\pi _M , \Omega_M]=0,
\end{array}
\]
where $\Omega _M$ is the 3-vector field $\Omega _M=\rho (\Omega )$,
and $\rho :\gm (\wedge ^3A)\to \mathfrak X ^3(M)$ is the extension
of the anchor map.
\end{prop}

Now we are ready to state the main theorem of this section:
quasi-Lie bialgebroids are indeed the infinitesimal invariants of
quasi-Poisson groupoids.

\begin{them}\label{inf-inv-quasi}
If $(\gm\gpd M, \Pi ,\Omega )$ is a quasi-Poisson groupoid, then
there exists a natural quasi-Lie bialgebroid structure $(\delta
,\Omega )$ on the Lie algebroid $A$ of $\gm$.

Conversely, if $(A, \delta ,\Omega )$ is a quasi-Lie bialgebroid,
where $A$ is the Lie algebroid of an $\alp$-connected and
$\alp$-simply connected Lie groupoid $\gm$, there exists a
quasi-Poisson groupoid structure $(\Pi ,\Omega )$ on $\gm$ such
that the corresponding quasi-Lie bialgebroid is $(\delta ,\Omega
)$.
\end{them}
\begin{pf}
Let $(\gm \gpd M, \Pi ,\Omega )$ be a quasi-Poisson groupoid.
Since $\Pi $ is a multiplicative bivector field, it induces a
2-differential $\delta _\Pi$ on $A$ according to Theorem
\ref{derivation}. In addition since  $\Pi \mapsto \delta
_\Pi$ preserves the graded Lie algebra structures,
from   Eqs.
(\ref{graded-bracket}) and (\ref{quasiPoisson1}),  it follows
that
\[
\delta _\Pi \smalcirc \delta _\Pi = \half [\delta _\Pi ,\delta
_\Pi ]=ad (\Omega ).
\]
Moreover Eq. (\ref{quasiPoisson2})  implies that $\delta
_\Pi (\Omega )=0$. As a consequence, $(A,\delta _\Pi ,\Omega )$ is
a quasi-Lie bialgebroid.

Conversely, let $(A,\delta ,\Omega)$ be a quasi-Lie bialgebroid.
Since $\delta$ is a 2-differential, according to Theorem
\ref{derivation}, there exists a multiplicative bivector field
$\Pi$ on $\gm$ such that $\delta _\Pi =\delta$. On the other hand,
the 3-differential $ad(\Omega )$ can be integrated to the
multiplicative 3-vector field $\Vec{\Omega}-\ceV{\Omega}$. In
addition, using the identity $\delta \smalcirc \delta =\lcf \Omega
,\cdot \rcf $,  we have
\[
\delta _{[\Pi ,\Pi ]}=[\delta  ,\delta  ]=2\delta \smalcirc \delta
=2 ad(\Omega )=2\delta _{ \Vec{\Omega}-\ceV{\Omega}}.
\]
Thus, by Theorem \ref{derivation} again, we have
$$
\half [\Pi ,\Pi ]=\Vec{\Omega}-\ceV{\Omega}.
$$
Moreover, $[ \Pi ,\Vec{\Omega} ]=\Vec{\delta _\Pi
\Omega}=\Vec{\delta \Omega}=0$. Thus, we conclude that $(\gm\gpd
M, \Pi ,\Omega )$ is a quasi-Poisson groupoid integrating the
quasi-Lie bialgebroid $(A,\delta ,\Omega )$.
\end{pf}

\begin{rmk}
{\rm Theorem \ref{inf-inv-quasi} generalizes a result of
Kosmann-Schwarzbach regarding quasi-Poisson Lie groups and
quasi-Lie bialgebras \cite{Kosmann:1991}. On the other hand, when
$\Omega =0$, that is, $(A,\delta )$ is a Lie bialgebroid, we
recover the classical results in
\cite{MackenzieX:1994,MackenzieX:2000}: there exists a one-to-one
correspondence between Lie bialgebroids and Poisson groupoids. }
\end{rmk}

As another consequence of  Theorem \ref{inf-inv-quasi},
 we recover the construction of
the non-degenerate twisted symplectic groupoid associated with  a
twisted Poisson manifold \cite{CX}.

\begin{cor}
Let $(M,\pi , \phi )$ be a twisted Poisson structure. If the Lie
algebroid  ${T^*M_{(\pi,\phi)}}$ can be integrated to an
$\alp$-simply connected and $\alp$-connected Lie groupoid $\Gamma$,
then $\Gamma$ is a non-degenerate twisted symplectic groupoid.
\end{cor}
\begin{pf}
Let  $\lambda :T^\ast M\to A$ be the Lie algebroid isomorphism
between the Lie algebroid $A$ of $\gm$ and $T^\ast M$.
Denote   by $\delta$  and $\Omega$  the almost 2-differential
and the  3-section of $A$ respectively such that
\begin{equation}\label{isomorphism}
\begin{array}{l}
\Omega = \lambda ( \phi );\\
-\lambda (\delta _{\pi ,\phi} \varphi )=\delta (\lambda
 (\varphi )), \forall \varphi \in \Omega ^\bullet (M),
\end{array}
\end{equation}
where $\delta _{\pi,\phi }$ is defined as  in Example
\ref{twisted-Poisson}. $(A,\delta ,\Omega )$ is clearly a
quasi-Lie bialgebroid induced by the quasi-Lie bialgebroid
structure $(T^\ast M,\delta _{\pi ,\phi},\phi )$. Therefore,
according to Theorem \ref{inf-inv-quasi}, there exists a bivector
field $\Pi$ such that $(\gm \gpd M,\Pi ,\Omega )$ is a
quasi-Poisson groupoid satisfying $\delta _\Pi =\delta$. Using
Eqs. (\ref{how-to-go-down}) and (\ref{isomorphism}),  we  have
\[
\lambda (df )=\lambda (\delta _{\pi ,\phi} f )=-\delta f=\Pi
^\sharp (\alp ^\ast df), \forall f\in C^\infty (M).
\]
Thus, $\lambda (\eta )=\Pi ^\sharp (\alp ^\ast \eta)$, $
\forall \eta \in T^\ast M$. As a consequence, $\lambda ^\ast :A^\ast \to TM$ is
given by $\lambda (\xi )=-\alpha _\ast \Pi ^\sharp (\xi )$
$ \forall \xi \in A^\ast$, where $A^\ast$ is identified with the conormal
bundle of $M$. Hence, $\lambda ^\ast (\xi )=-\Pi ^\sharp (\xi )$,
because $\Pi ^\sharp (\xi )$ is tangent to $M$. Therefore,
\[
\Pi ^\sharp (\xi +\alp ^\ast df )=-\lambda ^\ast (\xi )+\lambda
(df), \forall \xi \in A^\ast \textrm{ and }f\in C^\infty (M).
\]
Since any element of the cotangent space
 $T^\ast _{\epsilon (m)}\gm$
($m\in M$) can be written as $\xi +\alp ^\ast df$, with $\xi \in
A^\ast$ and $f\in C^\infty (M)$, and  $\lambda$ and
$\lambda ^\ast$ are injective ($\lambda$ is a vector bundle
isomorphism), it follows that $\Pi$ is non-degenerate along $M$.
Following Theorem 5.3 in \cite{MackenzieX:2000}, one can extend
this non-degeneracy for every point in $\Gamma$. From Proposition
\ref{non-degenerate-quasi}, we conclude that $\gm$ is the twisted
symplectic groupoid integrating $(M,\pi,\phi)$
\end{pf}
%\begin{rmk}
%{\rm Theorem \ref{inf-inv-quasi} generalizes the result of
%Y.~Kosmann-Schwarzbach regarding quasi-Poisson Lie groups and
%quasi-Lie bialgebras \cite{Kosmann:1991}. On the other hand, when
%$\Omega =0$, that is, $(A,\delta )$ is a Lie bialgebroid, we
%recover the classical results in
%\cite{MackenzieX:1994,MackenzieX:2000}: there exists a one-to-one
%correspondence between Lie bialgebroids and Poisson groupoids. }
%\end{rmk}

We end this section with the following proposition,  which reveals
the relation between the bivector fields obtained by Propositions
\ref{bivector-onthebase} and \ref{4.9}.

\begin{prop}\label{bivectors-relation}
Let $(\gm\gpd M,\Pi, \Omega )$ be a quasi-Poisson groupoid,
 and $(A, \delta _\Pi, \Omega )$ its corresponding quasi-Lie bialgebroid.
By $\Pi_M$ and $\pi _M$ we denote the bivector fields on $M$
induced from the quasi-Poisson groupoid structure on $\gm$
and the quasi-Lie bialgebroid structure on  $A$ as in  Proposition
\ref{bivector-onthebase}  and Proposition \ref{4.9} respectively.
Then
\[
\Pi_M =\pi _M .
\]
\end{prop}
\begin{pf}
Let $f,g\in C^\infty (M)$. Then, using Eq. (\ref{how-to-go-down}),
\[
\begin{array}{rcl}
-\langle \rho (\delta f),dg \rangle &=&-\langle \alp _\ast ( [\Pi
,\alp ^\ast f] ),dg \rangle\\ &=&  \Pi (\alp ^\ast df ,\alp ^\ast
dg ) = (\alp _\ast \Pi )(df ,dg) .
\end{array}
\]
The conclusion thus  follows.
\end{pf}
\subsection{Hamiltonian $\gm$-spaces of quasi-Poisson groupoids}
Let $\Gamma \gpd M$ be a Lie groupoid. Recall that a {\em
$\gm$-space} is a smooth manifold $X$ with a map $J:X\to M$,
called the {\em momentum map},
and an action
$$\Gamma \times _M X=\{ (g,x)\in \Gamma \times X \, |\, \bet (g)=J(x)\}\to X,
\ \ \
(g,x)\mapsto g\cdot x$$ satisfying
\begin{enumerate}
\item $J(g\cdot x)=\alp (g)$, for $(g,x)\in \Gamma \times _M X$;
\item $(g h)\cdot x=g\cdot (h\cdot x)$, for $g,h\in \Gamma $
and $x\in X$ such that $\bet (g)=\alp (h)$ and $J(x)=\bet
(h)$;
\item $\epsilon (J(x))\cdot x =x$, for $x\in X$.
\end{enumerate}
Hamiltonian $\Gamma$-spaces for Poisson groupoids were studied in
\cite{LWX}. For quasi-Poisson groupoids, one can introduce
Hamiltonian $\gm$-spaces in a similar fashion.
\begin{defn}
Let $(\Gamma \gpd M, \Pi , \Omega )$ be a quasi-Poisson groupoid.
A {\em Hamiltonian $\gm$-space} is a $\gm$-space $X$ with momentum
map $J:X\to M$ and a bivector field $\Pi _X\in \mathfrak X ^2(X)$
such that:
\begin{enumerate}
\item the graph of the action $\{ (g,x,g\cdot x) \,|\, J(x)=\bet (g)
\}$ is a coisotropic submanifold of $(\Gamma \times X \times X,\Pi
\oplus \Pi _X\oplus -\Pi _X)$;
\item $\half[\Pi _X ,\Pi _X]=\hat{\Omega}$, where the hat denotes the map
$\Gamma (\wedge ^3A)\to \mathfrak X ^3(X)$, induced by the
infinitesimal action of the Lie algebroid on $X$: $\Gamma (A)\to
\mathfrak X (X)$, $Y\mapsto \hat{Y}$ .
\end{enumerate}
\end{defn}
\begin{prop}
Let $(\Gamma \gpd M, \Pi , \Omega )$ be a quasi-Poisson groupoid.
If $(X,\Pi _X)$ is a Hamiltonian $\Gamma$-space with momentum map
$J,$ then $J$ maps $\Pi _X$ to $\Pi _M$, where $\Pi _M$ is given
by Proposition \ref{bivector-onthebase}.
\end{prop}
\begin{pf}
The result follows using the coisotropy condition and the fact
that $(-\bet ^\ast \eta, J^\ast \eta, 0)$, with $\eta \in T^\ast
M$, is conormal to the graph of the groupoid action. Indeed,
\[
0=\Pi (-\bet ^\ast \eta _1,-\bet ^\ast \eta _2) +\Pi (J^\ast \eta
_1,J ^\ast \eta _2)=\Big ( \bet _\ast \Pi + J_\ast \Pi _X\Big
)(\eta _1,\eta _2)
\]
for any $\eta _1,\eta _2\in T_m^\ast M$. Thus, $J_\ast \Pi _X =
-\bet _\ast \Pi =\Pi _M$.
\end{pf}
The following theorem gives an equivalent description of
Hamiltonian $\Gamma$-spaces of quasi-Poisson groupoids in terms of
their infinitesimal objects.
\begin{them}\label{inf-action-charact}
Let $(\Gamma \gpd M, \Pi , \Omega )$ be a quasi-Poisson groupoid
with corresponding quasi-Lie bialgebroid $(A,\delta _\Pi ,\Omega
)$. Then $(X,\Pi _X)$ is a Hamiltonian $\Gamma$-space with
momentum map $J:X\to M$ if and only if $\half[\Pi _X ,\Pi
_X]=\hat{\Omega}$ and
\begin{equation}\label{inf-eqs}
\begin{array}{l}
[\Pi _X,J^\ast f]= \widehat{\delta _\Pi f}, \ \forall f\in C^\infty (M),\\[5pt]
[\Pi _X,\hat{Y}]=\widehat{\delta _\Pi (Y)}, \ \forall  Y\in \Gamma
(A) .
\end{array}
\end{equation}
\end{them}
\begin{pf}
Using Theorem 7.1 in \cite{LWX} we know that the coisotropy
condition is equivalent to the following conditions:
\begin{enumerate}
\item For any $f\in C^\infty (M)$, $X_{J^\ast f}(x)=
(r_x)_\ast X_{\alp ^\ast f}(u)$, where $x\in X$, $u=J(x)$ and
$r_x$ denotes the map $g\mapsto g\cdot x$ from $\bet ^{-1}(u)$ to
$X$.
\item For any compatible $(g,x)\in \Gamma \times _M X$,
\[
\Pi _X(g\cdot x)=(L_{\cal X})_\ast \Pi _X(x)+(R_{\cal Y})_\ast \Pi
(g)- (R_{\cal Y})_\ast (L_{\cal X})_\ast \Pi (u),
\]
where $u=\bet (g)=J(x)$, $\cal X$ is any local bisection
through $g$, and $\cal Y$  is any local section of $J$ through the
point $x$.
\end{enumerate}
From the first condition and the equation
 $\Vec{\delta _\Pi f}= [\Pi
,\alp ^\ast f]$, it follows that $[\Pi _X,J^\ast f]=
\widehat{\delta _\Pi f}$.

On the other hand, let $x\in X$ and $u=J(x)$. Applying the second
condition to
 the family of (local) bisections ${\cal X}_t=\exp
tY$ associated to $Y\in \Gamma (A)$ and $g_t=(\exp tY)(u)$, one
obtains that
\[
(L^{-1}_{{\cal X}_t})_\ast \Pi _X(g_t\cdot x)= \Pi _X(x)+(R_{\cal
Y})_\ast (L^{-1}_{{\cal X}_t})_\ast \Pi (g_t)- (R_{\cal Y})_\ast
\Pi (u),
\]
for any section ${\cal Y}$ of $J$. Then, taking derivatives, we
have
\[
[\Pi _X,\hat{Y}]=\widehat{\delta _\Pi (Y)}.
\]
The other direction follows just by going backwards.
\end{pf}
\begin{rmk}
{\rm Note that Eq. (\ref{inf-eqs}) is equivalent to
that  $[\Pi _X,\hat{P}]= \widehat{\delta _\Pi P}$
 for any $P\in \Gamma (\wedge
^k A)$, $k\geq 0$. I.e., the following diagram is commutative,

\begin{picture}(375,60)(40,50)
\put(200,20){\makebox(0,0){$\gm (\wedge ^{k+1}A)$}}
\put(305,20){\makebox(0,0){$\mathfrak X^{k+1}(X)$}}
\put(225,20){\vector(1,0){55}} \put(200,80){\makebox(0,0){$\gm
(\wedge ^{k}A)$}} \put(200,70){\vector(0,-1){40}}
\put(300,70){\vector(0,-1){40}} \put(225,80){\vector(1,0){55}}
\put(300,80){\makebox(0,0){$\mathfrak X^k(X)$}}

\put(250,85){$\widehat{}$} \put(250,25){$\widehat{}$}

\put(180,50){\makebox(0,0){$\delta _\Pi$}}
\put(320,50){\makebox(0,0){$[\Pi _X,\cdot ]$}}

\put(375,50){\makebox(0,0){$\forall \mbox{ k }\geq 0.$}}
\end{picture}

}
\end{rmk}

\vspace{1cm}

\subsection{Twists of quasi-Lie bialgebroids}
Just like quasi-Lie bialgebras, one can talk about twists of a
quasi-Lie bialgebroid. Given a quasi-Lie bialgebroid $(A, \delta,
\Omega )$ and a section $t\in \gm (\wedge^2 A)$, let
$\delta^t=\delta +\lcf t, \cdot \rcf$, and $ \Omega^t = \Omega
+\delta t+\half \lcf t, t\rcf$. Using Eq. (\ref{k-dif-prop}) and
the properties of the Schouten bracket, it is simple to see that
$(A, \delta^t, \Omega^t )$ is also a quasi-Lie bialgebroid, which
will be called the {\em twist of $(A, \delta, \Omega )$ by $t \in
\gm (\wedge^2 A)$}.

The following proposition describes how a twist affects the
integration (see \cite{AK} for the case of quasi-Lie bialgebras).
\begin{prop}
Assume that $(\gm \gpd M, \Pi ,\Omega )$ is a
 quasi-Poisson groupoid with corresponding
quasi-Lie bialgebroid $(A, \delta, \Omega )$, and $t\in \gm
(\wedge^2 A)$. Then $(\gm \gpd M, \Pi^t ,\Omega^t )$ is a
quasi-Poisson groupoid with quasi-Lie bialgebroid $(A, \delta^t,
\Omega^t )$, where $\Pi^t=\Pi+\Vec{t}- \ceV{t}$.
\end{prop}
\begin{pf}
It is a direct consequence of Theorem \ref{derivation}.
Note  that  if $\Pi$ and $\Pi '$ are multiplicative bivector fields then
$\delta _{\Pi +\Pi '}=\delta _{\Pi} +\delta_{\Pi '}$.
Also we have
$\delta_{\Vec{t}-\ceV{t}}=ad(t)$ (see Example \ref{dif-cob}).
\end{pf}

Next, we will show how a Hamiltonian $\gm$-space of a
quasi-Poisson groupoid $(\gm\gpd M, \Pi ,\Omega)$ must be modified
in order to obtain a Hamiltonian $\gm$-space for the twisted
quasi-Poisson structure $(\Pi ^t,\Omega ^t)$.
\begin{prop}
\label{pro:twist-ham} Assume that $(\gm \gpd M, \Pi ,\Omega )$ is
a quasi-Poisson groupoid and $t\in \gm (\wedge^2 A)$. There is a
bijection between Hamiltonian $\gm$-spaces of $(\gm \gpd M, \Pi
,\Omega )$ and those of $(\gm \gpd M, \Pi^t ,\Omega^t )$.

More precisely, if $(X\to M, \Pi_X)$ is a Hamiltonian $\gm$-space
of $(\gm \gpd M, \Pi ,\Omega )$, then $(X\to M,  \Pi_X+\hat{t})$
is a Hamiltonian $\gm$-space of $(\gm \gpd M, \Pi^t ,\Omega^t )$.
\end{prop}
\begin{pf}
Let $(X\to M, \Pi_X)$ be a Hamiltonian $\gm$-space of the
quasi-Poisson groupoid $(\gm \gpd M, \Pi ,\Omega )$. Using Theorem
\ref{inf-action-charact}, we deduce that
\[
[\Pi _X,\hat{P}]= \widehat{\delta _\Pi P}\mbox{ for any }P\in
\Gamma (\wedge ^k A),\quad k\geq 0.
\]
Therefore, from the fact that $\oplus_k \Gamma (\wedge^k A)\to
\oplus_k {\mathfrak X}^k (X)$, $Y\mapsto \hat{Y}$ is a graded Lie algebra morphism, one
gets that
\[
[\Pi _X+\hat{t},\hat{P}]=\widehat{\delta _\Pi P}+\widehat{\lcf
t,P\rcf}=\widehat{\delta _{\Pi ^t}P}.
\]
On the other hand, it is trivial to see that
\[
\half [\Pi _X+\hat{t},\Pi _X+\hat{t}]=\widehat{\Omega ^t}.
\]
As a consequence of Theorem \ref{inf-action-charact}, $(X\to M,
\Pi _X+\hat{t})$ is a Hamiltonian $\gm$-space of the quasi-Poisson
groupoid $(\gm\gpd M, \Pi ^t,\Omega ^t)$.
\end{pf}

\subsection{Quasi-Poisson groupoids associated to  Manin pairs}
\label{Manin-pair}
In this subsection we will  describe an
example of quasi-Poisson groupoid associated to a Manin
quasi-triple.

Let $(\mathfrak d, \mathfrak g)$ be a {\em Manin pair}, that is,
${\mathfrak d}$ is an even dimensional Lie algebra with an
invariant, nondegenerate symmetric bilinear form, and $\mathfrak
g$ is a maximal isotropic subalgebra of $\mathfrak d$. In this
case, one can integrate the Manin pair $(\mathfrak d, \mathfrak
g)$ to the so-called {\em group pair} $(D, G)$, where $D$ and $G$
are connected and simply connected Lie groups with Lie algebra
${\mathfrak d}$ and $\mathfrak g$ respectively. Furthermore, the
action of the Lie group $D$ on itself by left multiplication
induces an action of $D$ on $S=D/G$,  and in particular a $G$-action
 on $S$,  which is called the {\em dressing action}. As in \cite{AK},
the infinitesimal dressing action is
 denoted by $v\mapsto v_S$ for any $v\in \mathfrak d$.

If $\mathfrak h$ is an isotropic complement of $\mathfrak g$ in
$\mathfrak d$, by identifying $\mathfrak h$ with $\mathfrak
g^\ast$, we obtain a quasi-Lie bialgebra structure on $\mathfrak
g$, with cobracket $F:\mathfrak g\to \wedge ^2 \mathfrak g$ and
$\Omega \in \wedge ^3 \mathfrak g$. If $\{ e_i \}$ is a basis of
$\mathfrak g$ and $\{ \epsilon ^i\}$ the dual basis of $\mathfrak
g^\ast \cong \mathfrak h$, then $F(e_i)=\half \sum _{j,k}F^{jk}_i
e_j\wedge e_k$ and $\Omega =\frac{1}{6}\sum _{i,j,k} \Omega
^{ijk}e_i\wedge e_j\wedge e_k$. Moreover, the bracket on
$\mathfrak d \cong \mathfrak g\oplus \mathfrak h$ can be written
as
\begin{equation}\label{d-bracket}
[e_i,e_j]_{\mathfrak d}=\displaystyle\sum _{k=1}^n c_{ij}^k e_k,
\quad [e_i,\epsilon ^j]_{\mathfrak d}=\displaystyle\sum _{k=1}^n
-c^j_{ik}\epsilon ^k +F_i^{jk}e_k, \quad [\epsilon ^i,\epsilon
^j]_{\mathfrak d}=\displaystyle\sum _{k=1}^n F_k^{ij}\epsilon ^k+
\Omega ^{ijk}e_k,
\end{equation}
where $c_{ij}^k$ are the structure constants of the Lie algebra
$\mathfrak g$ with respect to the basis $\{ e_i \}$.
\begin{numex}\label{double}
Let $\mathfrak g$ be a Lie algebra endowed with a nondegenerate
symmetric bilinear form $K$. On the direct sum $\mathfrak
d=\mathfrak g\oplus \mathfrak g$ one can construct a scalar
product $(\cdot | \cdot )$ by
\[
((u_1,u_2)|(v_1,v_2))=K(u_1,v_1)-K(u_2,v_2),
\]
for $(u_1,u_2),(v_1,v_2)\in \mathfrak d$. Then, $(\mathfrak
d,\Delta (\mathfrak g),\half \Delta _- (\mathfrak g))$ is a Manin
quasi-triple, where $\Delta (v)=(v,v)$ and $\Delta _-(v)=(v,-v)$,
 $\forall v\in\mathfrak g$ (see \cite{AK}). In this case, as far as the
corresponding quasi-Lie bialgebra is concerned, the cobracket $F$
vanishes and $\Omega$ can be identified with the trilinear form on
$\mathfrak g$ given by $(u,v,w)\mapsto
\frac{1}{4}K(w,[u,v]_\mathfrak g)$.
\end{numex}
Let $\lambda :T^\ast _sS\to \mathfrak g$ be the dual map of the
infinitesimal dressing action $\mathfrak g^\ast \cong \mathfrak h
\to T_sS$. That is,
\[
\langle \lambda (\theta _s),\eta \rangle = \langle \theta _s,\eta
_S(s) \rangle, \forall \theta _s \in T^\ast _sS\mbox{ and }\eta
\in \mathfrak h.
\]
A direct consequence is that
\begin{equation}\label{def-lambda}
\lambda (df)=\displaystyle \sum _{i=1}^n (\epsilon
^i)_S(f)e_i,\mbox{ for }f \in C^\infty (S).
\end{equation}
\begin{rmk}
Recall   that an isotropic complement $\mathfrak h$
is said to be {\em admissible} at a point $s\in S=D/G$ if the
infinitesimal dressing action restricted to $\mathfrak h$ defines
an isomorphism from $\mathfrak h$ onto $T_sS$  \cite{AK}. In this case, we
can define an isomorphism from $\mathfrak g$ to $T^\ast _s S$,
$\xi\mapsto \xi _{\mathfrak h} (s)$, as follows,
\[
\langle \xi_{\mathfrak h}{}{(s)},\eta _S {}{(s)} \rangle = - ( \xi
|\, \eta ), \ \forall \eta \in \mathfrak h,
\]
where in the right-hand side $(\cdot | \cdot )$ is the bilinear
form on $\mathfrak d$.  $(\mathfrak d, \mathfrak g,
\mathfrak h)$ is said to be an {\em admissible} Manin quasi-triple if
$\mathfrak h$ is admissible at every point of $S$. In this case,
if $\lambda (\theta _s)=\xi$ then $\xi_{\mathfrak
h}{}{(s)}=-\theta _s$.
\end{rmk}
Next, we will show that on the transformation Lie algebroid
$\mathfrak g \times S\to S$ there exists a natural quasi-Lie
bialgebroid structure.
\begin{prop}
Assume that $(\mathfrak d, \mathfrak g, \mathfrak h)$ is a Manin
quasi-triple with associated quasi-Lie bialgebra $(\mathfrak
g,F,\Omega )$. On the transformation Lie algebroid $\mathfrak
g\times S\to S$ (where $\mathfrak g$ acts on $S$ by the
infinitesimal dressing action) define an almost 2-differential
\begin{equation}\label{2-diff-trans}
\begin{array}{l}
\delta  (f)=\lambda (df),\mbox{ for }f \in
C^\infty (S),\\[10pt]
\displaystyle \delta  \xi = - F(\xi),\mbox{ for }\xi \in \mathfrak
g,
\end{array}
\end{equation}
where $\lambda$ is defined by Eq. (\ref{def-lambda}) and $\xi \in
\mathfrak g$ is considered as a constant section of the Lie
algebroid $\mathfrak g\times S\to S$, extending this operation to
arbitrary section using the derivation law. Then, $(\mathfrak
g\times S, \delta ,\Omega )$ is a quasi-Lie bialgebroid.
\end{prop}
\begin{pf}
We remark that it suffices to check the axioms of a quasi-Lie
bialgebroid for functions on $S$ and constant sections of
$\mathfrak g\times S$, since the general case follows from the
derivation law.

If $f\in C^\infty (S)$,  then using Eq. (\ref{d-bracket}),
\[
\begin{array}{rcl}
\delta ^2 f &=&\displaystyle\sum _{i=1}^n \delta ( (\epsilon ^i)_S
(f) e_i )\\
& =& \displaystyle\sum _{i,j=1}^n(\epsilon ^j)_S (\epsilon
^i)_S (f) e_j\wedge e_i +\sum _{i=1}^n (\epsilon ^i)_S(f) \delta
(e_i)
\\&=& \displaystyle\sum _{i<j} \Big ( (\epsilon
^i)_S (\epsilon ^j)_S (f) -  (\epsilon ^j)_S (\epsilon ^i)_S (f)
\Big ) e_i\wedge e_j +\sum _{i=1}^n (\epsilon ^i)_S(f) \delta
(e_i)\\&=& \displaystyle\sum _{i<j}  ( [\epsilon ^i, \epsilon ^j
]_\mathfrak d)_S (f)  e_i\wedge e_j +\sum _{k=1}^n (\epsilon
^k)_S(f) \delta (e_ k)\\ &=&\displaystyle \sum _{i<j}\sum _{k=1}^n
( F^{ij}_k (\epsilon ^k )_S(f) +\Omega ^{ijk} (e_k)_S(f) )
e_i\wedge e_j -\sum _{k=1}^n (\epsilon ^k)_S(f) \sum _{i<j}
F^{ij}_k e_i\wedge e_j\\&=& \displaystyle\sum _{i<j} \sum _{k=1}^n
\Omega ^{ijk} (e_k)_S(f) e_i\wedge e_j\\
& =& \displaystyle \half \sum
_{i,j,k=1}^n \Omega ^{ijk} (e_k)_S(f) e_i\wedge e_j
\\&=& \lcf \Omega ,f\rcf .
\end{array}
\]
Moreover, since $(\mathfrak g, F,\Omega )$ is a quasi-Lie
bialgebra, we have
\[
\delta ^2 \xi = \lcf \Omega ,\xi \rcf ,  \forall \xi \in
\mathfrak g .
\]
Next, let us show that $\delta \lcf \xi ,f \rcf = \lcf \delta \xi
,f \rcf +  \lcf \xi ,\delta f \rcf$, for $\xi \in \mathfrak g$ and
$f\in C^\infty (S)$. Taking $\xi =e_i$ and using Eq.
(\ref{d-bracket}),
\[
\begin{array}{rcl}
\lcf \delta e_i ,f \rcf \kern-2pt+\kern-2pt\lcf e_i,\delta f\rcf
\kern-2pt-\kern-2pt\delta \lcf e_i,f\rcf % &=& \lcf \displaystyle
%\half \sum _{j,k=1}^n - F_i^{jk} e_j\wedge e_k ,f\rcf +\lcf
%e_i,\sum _{j=1}^n(\epsilon ^j)_S(f)e_j \rcf
%\\ && - \displaystyle
%\sum _{j=1}^n (\epsilon ^j)_S(e_i)_S(f)e_j\\
&=& \displaystyle\sum _{j,k =1}^n (e_j)_S(f)F_i^{jk} e_k + \sum
_{j=1}^n(e_i)_S(\epsilon ^j )_S(f) e_j \\ &&+\displaystyle\sum
_{j,k =1}^n (\epsilon ^j)_S(f) c_{ij}^k e_k -\displaystyle\sum _{j
=1}^n (\epsilon ^j)_S(e_i)_S(f)e_j
\\&=& \displaystyle\sum _{j =1}^n \Big ( (e_i)_S(\epsilon ^j)_S(f)e_j -
(\epsilon ^j )_S(e_i)_S(f) e_j \Big )
\\ && +\displaystyle\sum _{j,k =1}^n \Big ( (e_j)_S(f)F_i^{jk}  +
(\epsilon ^j)_S(f) c_{ij}^k \Big ) e_k
\\&=& \displaystyle\sum _{j =1}^n [ (e_i)_S,(\epsilon ^j)_S](f) e_j
- \displaystyle\sum _{j,k =1}^n \Big (  F_i^{jk}e_k - c_{ik}^j
\epsilon ^k\Big )_S(f) e_j
\\&=& \displaystyle\sum _{j =1}^n [ (e_i)_S,(\epsilon ^j)_S](f) e_j
- \displaystyle\sum _{j =1}^n (  [ e_i,\epsilon ^j]_\mathfrak d
)_S(f) e_j
%\\&=& \displaystyle\sum _{j =1}^n \Big ( [ (e_i)_S,(\epsilon ^j)_S] -
%( [ e_i,\epsilon ^j]_\mathfrak d )_S \Big )(f) e_j
\\ &=& 0.
\end{array}
\]
Since we also know that $\delta [\xi _1,\xi _2]=[\delta \xi _1
,\xi _2 ]+[\xi _1,\delta \xi _2]$ for any $\xi _1,\xi _2\in
\mathfrak g$, the conclusion thus follows.
\end{pf}
Now, consider the transformation groupoid $\Gamma : G\times S \gpd
S$ associated to the dressing action. Theorem \ref{inf-inv-quasi}
implies that $\Gamma$ is a quasi-Poisson groupoid. In what
follows, we will explicitly describe the  multiplicative bivector
field $\Pi$ on $\Gamma$.

The quasi-Lie bialgebra $(\mathfrak g,F,\Omega )$ implies that $G$
is a quasi-Poisson Lie group with multiplicative bivector field
denoted by $\Pi _G$. Moreover, there exists a bivector field $\Pi
_S$ on $S$ given by $\Pi _S =- \sum _{i=1}^n (e_i)_S\otimes
(\epsilon ^i)_S$, i.e.,
\begin{equation}\label{bivector-S}
\Pi _S{} (df, dg)=- \sum _{i=1}^n (\epsilon ^i)_S(f) (e_i)_S(g),
\forall f,g\in C^\infty (S).
\end{equation}
Using Proposition \ref{4.9}, Eqs. (\ref{2-diff-trans}) and
(\ref{bivector-S}), we directly deduce the following

\begin{prop}\label{relation-bivector-S}
Let $(\mathfrak d, \mathfrak g, \mathfrak h)$ be a Manin
quasi-triple and $(\mathfrak g\times S \to S,\delta ,\Omega )$ the
corresponding quasi-Lie bialgebroid. Then, the bivector field on
$S$ induced by $\delta$ as in Proposition \ref{4.9} coincides with
$\Pi _S=- \sum _{i=1}^n (e_i)_S\otimes (\epsilon ^i)_S$.
\end{prop}
\begin{pf}
If $f,g\in C^\infty (S)$, then
\[
- (\delta f)_S (g)= -\displaystyle \sum _{i=1}^n (\epsilon
^i)_S(f) (e_ i)_S(g)= \Pi _S (df, dg).
\]
\end{pf}
\begin{numex}\label{double2}
Let $\mathfrak g$ be a Lie algebra endowed with a nondegenerate
symmetric bilinear form $K$ and consider the corresponding Manin
quasi-triple $(\mathfrak d,\Delta (\mathfrak g),\half \Delta _-
(\mathfrak g))$ (see Example \ref{double}). Then, $\mathfrak g$
acts on $S=G$ by the adjoint action and the 2-differential is
given by
\[
\begin{array}{l}
\delta  (f)=\displaystyle \half \sum _{i=1}^n
(\Vec{e_i}+\ceV{e_i})(f)e_i,  \ \ \forall f \in
C^\infty (G),\\[10pt]
\displaystyle \delta  \xi = 0, \ \ \forall  \xi \in \mathfrak g,
\end{array}
\]
where $\{ e_i\}$ is an orthonormal basis of $\mathfrak g$.
Moreover, the bivector field on $G$ induced by $\delta$ is
\[
\Pi _G=\half \sum _{i=1}^n \ceV{e_i}\wedge \Vec{e_i},
\]
which was first obtained in \cite{AK} (see also \cite{AKM}).
\end{numex}

Now we are ready to   describe the quasi-Poisson groupoid structure on the
transformation groupoid $G\times S \gpd S$.
\begin{them}\label{Manin-groupoid}
Assume that $(\mathfrak d, \mathfrak g, \mathfrak h)$ is a Manin
quasi-triple. Define a bivector field $\Pi $ on $G\times S$ by
\begin{equation}\label{bivector-Manin}
\begin{array}{rcl}
\Pi \Big ( (\theta _g,\theta _s),(\theta'_g,\theta ' _s)\Big )
 &=& \Pi _G{} (\theta _g, \theta
'_g)-\Pi _S (\theta _s,\theta '_s) \\ && +\langle \theta '_s,
(L_g^\ast \theta _g)_S \rangle -\langle \theta _s,(L_g^\ast \theta
' _g)_S\rangle ,
\end{array}
\end{equation}
for any $(g,s)\in G\times S$, $\theta _g,\theta '_g \in T^\ast
_gG$, $\theta _s,\theta '_s \in T^\ast _sS$, where $(L_g^\ast
\theta _g)_S$ denotes the  vector field on $S$
 corresponding to the dressing action  of $L_g^\ast \theta
_g\in \mathfrak g^\ast \cong \mathfrak h\subset \mathfrak d$,
similarly for $(L_g^\ast \theta' _g)_S$. Then $(G\times S\toto S,
\Pi, \Omega )$ is a quasi-Poisson groupoid integrating the
2-differential $\delta$ given by Eq. (\ref{2-diff-trans}).
\end{them}
\begin{pf}
Let $\Pi$ be the multiplicative bivector field on $G\times S$
integrating $\delta$. Then, using Propositions
\ref{bivectors-relation} and \ref{relation-bivector-S} and the
fact that $\beta (g,s)=s$ for any $(g,s)\in G\times S$,
\[
\begin{array}{rcl}
\Pi ( (0,df ),(0,dg))&=&\Pi(\beta ^\ast (df),\bet ^\ast (dg))=
(\bet _\ast \Pi )(df, dg)\\&=&-\Pi _S(df,dg).
\end{array}
\]
Therefore,
\[
\Pi ( (0,\theta _s ),(0,\theta _s'))=-\Pi _S(\theta _s,\theta_s'),
\forall \theta _s,\theta _s' \in T^\ast _sS.
\]
Next, if $f\in C^\infty (S)$ and $\theta _g\in T^\ast _gG$, then
from Eqs. (\ref{how-to-go-down}), (\ref{def-lambda}) and
(\ref{2-diff-trans}), it follows that
\[
\begin{array}{rcl}
\Pi ( (0,df ),(\theta _g,0))&=&\langle \Pi ^\sharp (\beta ^\ast
(df)),(\theta _g,0)\rangle  = - \langle L_{(g,s)}{}_\ast ( \delta
f ),(\theta _g,0)\rangle \\&=& - \langle (L_{g})_\ast ( \delta f
),\theta _g \rangle= - \langle  \delta f ,L^\ast _{g} \theta _g
\rangle \\&=&-\langle  df ,(L^\ast _{g}\theta _g)_S \rangle .
\end{array}
\]
Thus,
\[
\Pi ( (0,\theta _s ),(\theta _g,0))= - \langle  \theta _s ,(L^\ast
_{g}\theta _g)_S \rangle ,\mbox{ for }\theta _g\in T^\ast _gG, \
\forall  \theta _s\in T^\ast _sS.
\]
Finally, fixing $s\in S$, we write $\Pi _s$ the bivector field on
$G$ defined by $\Pi _s(g)(\theta _g,\theta '_g)=\Pi (g,s)((\theta
_g,0),(\theta '_g,0))$. Then, it is clear that $\Pi _s$ is a
multiplicative bivector field on $G$. Moreover, if $\xi \in
\mathfrak g$  is considered as a constant section of $\mathfrak
g\times S\to S$, then $\Vec{\xi}=(\Vec{\xi}_G,0)$, where $\Vec{\xi}_G$
is  the right-invariant vector field on $G$. Using Eqs.
(\ref{how-to-go-down}) and (\ref{2-diff-trans}), we  get
\[
([\Pi _s,\Vec{\xi}_G],0)=[\Pi, \Vec{\xi}]=\Vec{\delta
\xi}=(-\Vec{F(\xi)}_G,0)=(-[\Vec{\xi}_G,\Pi _G],0)=([\Pi
_G,\Vec{\xi}]_G,0).
\]
As a consequence, $\Pi _s$ coincides with $\Pi _G$, that is,
\[
\Pi ( (\theta _g,0),(\theta'_g,0) ) = \Pi _G{} (\theta _g, \theta
'_g), \ \forall \theta _g,\theta '_g \in T^\ast
 _gG.
\]
Summing up, we can conclude that the multiplicative bivector field
$\Pi$ integrating $\delta$ is the one given by Eq.
(\ref{bivector-Manin}) and that $(\gm \gpd S, \Pi ,\Omega )$ is a
quasi-Poisson groupoid.
\end{pf}
\begin{numex}
{\rm In the particular case when $\mathfrak h$ is also a Lie
subalgebra of $\mathfrak d$, that is, $(\mathfrak d, \mathfrak g,
\mathfrak h)$ is a Manin triple, then $\Omega =0$ and $(G,\Pi _G)$
is a Poisson Lie group. Moreover, if the dressing action is
complete, we have $D/G\cong G^\ast ,$ the dual Poisson Lie group.
Thus, we recover a Poisson groupoid structure on $G\times G^\ast
\gpd G^\ast ,$ whose Poisson bivector field is described in
\cite{Lu:1990}.

On the other hand, when the dressing action is not complete, we
can still have a Poisson groupoid $G\times (D/G)\gpd (D/G)$, while
the groupoid $G\times G^\ast\gpd G^\ast$ does not exist any more.
It would be interesting  to study the relation between this  Poisson
groupoid  and the symplectic groupoid by Lu-Weinstein in \cite{LW}. }
\end{numex}
From Examples \ref{double} and \ref{double2} and Theorem
\ref{Manin-groupoid}, one can deduce the following

\begin{cor}\label{GxG}
Assume that $\mathfrak g$ is a Lie algebra  endowed
with a nondegenerate symmetric bilinear form $K$ and $G$ is its
corresponding connected and simply connected Lie group. Then the
transformation groupoid $G\times G\gpd G$, where $G$ acts on $G$ by
conjugation, together with
the
multiplicative bivector  field $\Pi$ on $G\times G$:
\[
\Pi (g,s) = \half \sum _{i=1}^n \ceV{e^2_i}\wedge \Vec{e^2_i}-
\ceV{e^2_i}\wedge \ceV{e^1_i}-\Vec{(Ad_{g ^{-1}}e_i)^2}\wedge
\Vec{e^1_i},
\]
the bi-invariant 3-form $\Omega: =\frac{1}{4}K(\cdot ,[\cdot ,\cdot ]_\mathfrak g)\in \wedge^3 {\mathfrak g}^*\cong \Omega^3 (G)^G$ on $G$,
is  a quasi-Poisson groupoid. Here $\{ e_i\}$ is an orthonormal
basis of ${\mathfrak g}$ and the superscripts
refer to the respective $G$-component.
%\[
%\Pi = \half \sum _{i=1}^n \ceV{e_i}\wedge \Vec{\epsilon _i}+
%\Vec{e_i}\wedge \ceV{\epsilon _i}.
%\]
\end{cor}
\begin{rmk}
We remark that, under the change of coordinates $(g,s)\mapsto
(a,b)=(s^{-1}g^{-1},g)$, $\Pi$ becomes the bivector field
 on $G\times G$ obtained in Example 5.3 of \cite{AKM}, i.e.,
\[
\Pi = \half \sum _{i=1}^n \ceV{e^1_i}\wedge \Vec{e^2 _i}+
\Vec{e^1_i}\wedge \ceV{e^2_i}.
\]
\end{rmk}

\subsection{$D/G$ momentum maps}

Next, we will investigate  the relation between quasi-hamiltonian
spaces with $D/G$-momentum map in the sense of \cite{AK} and
Hamiltonian $\gm$-spaces, where $\gm$ is the quasi-Poisson
groupoid $G\times S\gpd S$ associated to a Manin quasi-triple
$(\mathfrak d, \mathfrak g,  \mathfrak h )$ as described in
Section \ref{Manin-pair}.

First, we recall the notion of a quasi-hamiltonian space with
$D/G$-momentum map.
\begin{defn}\cite{AK}
Let $(G,\Pi_G,\Omega )$ be a connected quasi-Poisson Lie group
acting on a manifold $X$ with a bivector field $\Pi_X$. The action
$\Phi :G\times X \to X$ of $G$ on $X$ is said to be a {\em
quasi-Poisson action} if and only if
\begin{equation} \label{quasiPoissonAction1}
\Phi _\ast (\Pi _G\oplus \Pi _X) = \Pi _X,
\end{equation}
\begin{equation} \label{quasiPoissonAction2}
\half [\Pi _X,\Pi _X]=\Omega _X
\end{equation}
where $\Omega _X\in \mathfrak X ^3(X)$ is defined using the map
$\wedge ^3\mathfrak g \to \mathfrak X ^3(X)$ induced by the
infinitesimal action.
\end{defn}
Now, we recall the definition of a quasi-hamiltonian action with a
momentum map $J$ \cite{AK}. In fact our definition does not
require the assumption of $(\mathfrak d, \mathfrak g, \mathfrak
h)$ being admissible, so is more general.
\begin{defn}
Let $\Phi$ be a quasi-Poisson action of a quasi-Poisson Lie group
$(G,\Pi_G,\Omega )$ on $(X,\Pi_X)$. A $G$-equivariant map $J: X\to
S$ is called a {\em momentum map} if
\begin{equation}\label{quasiPoissonAction3}
\Pi _X^\sharp (J^\ast \theta _s)=-(\lambda (\theta _s))_X,\mbox{
for }\theta _s\in T^\ast _sS,
\end{equation}
where $G$ acts on $S$ by dressing action. The action is called
{\em quasi-hamiltonian} if it admits a momentum map and $X$ is
called a {\em quasi-hamiltonian space}.
\end{defn}
\begin{rmk}
If we consider an admissible isotropic complement $\mathfrak h$ of
$(\mathfrak d, \mathfrak g)$,
 then our definition coincides with the one given in \cite{AK}.
\end{rmk}
\begin{prop}
Let $(G,\Pi _G,\Omega)$ be a quasi-Poisson Lie group. If $X$ is a
quasi-hamiltonian space then the momentum map is a bivector map from
$(X,\Pi _X)$ to $(S,\Pi _S)$, i.e.,
\begin{equation}\label{mapping}
J_\ast \Pi _X=\Pi _S.
\end{equation}
\end{prop}
\begin{pf}
Let $f,g\in C^\infty (S)$. Using Eqs. (\ref{def-lambda}),
(\ref{quasiPoissonAction3}) and the $G$-equivariance, we have that
\[
\begin{array}{rcl}
\Pi _X(J^\ast df, J^\ast dg)&=&\langle \Pi _X^\sharp (J^\ast df ),
J^\ast dg \rangle = -J_\ast (\lambda (df )_X)(g)\\&=&- (\lambda
(df))_S(g)=\Pi _S(df, dg).
\end{array}
\]
\end{pf}
\begin{them}\label{hamiltonian-equivalence}
Let $(\mathfrak d ,\mathfrak g, \mathfrak h)$ be a Manin
quasi-triple. If  $(X,\Pi _X)$ is  a quasi-hamiltonian space with
momentum map $J:X\to S$,
 then $X$ is a Hamiltonian $\gm$-space, where $\gm$ is the
quasi-Poisson groupoid $(G\times S\gpd S, \Pi, \Omega)$. Here the
$\gm$-action is given by
\begin{equation}\label{extendedAction}
(g,s)\cdot x =\Phi (g,x), \ \forall  g\in G, s\in S \mbox{ and }
x\in X \mbox{ such that } J(x)=s,
\end{equation}
and $\Phi (g,x)$ denotes the $G$-action on $X$.

Conversely, if $(X,\Pi _X)$ is a Hamiltonian $\gm$-space with
momentum map $J: X\to S$, then $(X,\Pi _X)$ is a quasi-hamiltonian
space with momentum map $J:X\to S$, where the $G$-action  on $X$
is  given by
\[
\Phi (g,x)= (g,J(x))\cdot x, \ \forall g\in G\mbox{ and }x\in X.
\]
\end{them}
\begin{pf}
First of all, using the fact that $J:X\to S$ is $G$-equivariant,
we know that $\Phi$ can be extended to an action of $\Gamma$ on
$X$ by Eq. (\ref{extendedAction}). In addition, the conormal space
of the graph of the $\gm$-action at
 a point $((g,x),s,\Phi(g,x))$
is spanned by vectors of the form $(-(\Phi _x)^\ast \theta ,0
,-(\Phi _g)^\ast \theta ,\theta)$ for $\theta \in T^\ast _{\Phi
(g,x)}X$ and $(0 ,\gamma ,-J^\ast \gamma ,0)$ for $\gamma \in
T^\ast _sS$. We accordingly divide our proof into three different
cases:\\\\
{\bf Case 1.} The 2 covectors are of the first type. From the
definition of $\Pi$ (see Eq. (\ref{bivector-Manin})) we see that
the coisotropy condition is equivalent to
\[
\Phi_* (\Pi _G\oplus \Pi _X)=\Pi _X.
\]
{\bf Case 2.} The 2 covectors are of the second type. From Eq.
(\ref{bivector-Manin}),
 we deduce that $\Pi\oplus \Pi _X\oplus -\Pi
_X$ is coisotropic is equivalent to
\[
J_\ast \Pi _X=\Pi _S .
\]
{\bf Case 3.} One covector  is
of the first type and the other  of
the second. In this case, the coisotropy condition is just the
momentum map condition, i.e.,
\[
\Pi _X^\sharp (J^\ast \gamma )+(\lambda (\gamma ))_X=0.
\]
Finally, we note that Eq.~(\ref{quasiPoissonAction2}) is
equivalent to $\half[\Pi _X ,\Pi _X]=\hat{\Omega}$. We conclude
our proof.
\end{pf}

In \cite{AKM}, the authors study quasi-hamiltonian spaces for the
particular case when the Manin pair is the one associated with a
Lie algebra $\mathfrak g$ endowed with a nondegenerate symmetric
bilinear form $K$ (see Example \ref{double}). We now recall their
definitions and discuss the relation with
 Hamiltonian $\gm$-spaces of quasi-Poisson groupoids.

\begin{defn}\cite{AK,AKM}
A {\em quasi-Poisson manifold} is a $G$-manifold $X$, equipped
with an invariant bivector field $\Pi _X$ such that
\[
\half [\Pi _X,\Pi _X]= \Omega _X .
\]
\end{defn}
Moreover, one  can also  introduce the following specific
definition of momentum  maps.

\begin{defn}\cite{AK,AKM}
An Ad-equivariant map $J:X\to G$ is called a {\em momentum map}
for the quasi-Poisson manifold $(X,\Pi _X)$ if
\[
\Pi _X^\sharp (d(J^\ast f))=(J^\ast ({\cal D}f))_X, \
\forall f\in C^\infty (G),
\]
where ${\cal D}:C^\infty (G)\to C^\infty (G,\mathfrak g)$ is
defined by ${\cal D}f = \half \sum _{i=1}^n
(\Vec{e_i}+\ceV{e_i})(f)e_i$. The triple $(X,\Pi _X,J)$ is then
called a {\em Hamiltonian quasi-Poisson manifold}.
\end{defn}
\begin{rmk}
Note that the operator ${\cal D}$ is  exactly
 the 2-differential $\delta$ on the
Lie algebroid $\mathfrak g\times S \to S$
applying to functions (see Example \ref{double2}).
\end{rmk}

As a consequence of Theorem \ref{hamiltonian-equivalence} we
deduce the following
\begin{cor}
Let $G$ be a Lie group with Lie algebra $\mathfrak g$ endowed with
a nondegenerate symmetric bilinear form $K$. If $(X,\Pi _X)$ is a
hamiltonian quasi-Poisson manifold with momentum map $J:X\to G$,
where the action is denoted by $\Phi (g,x)$, then $X$ is a
Hamiltonian $\gm$-space, where $\gm$ is the quasi-Poisson groupoid
$(G\times G\gpd G, \Pi, \Omega)$ obtained in Corollary \ref{GxG}.
Here the $\gm$-action is given by
\[
(g,s)\cdot x =\Phi (g,x),  \ \ \forall  g\in G, s\in S \mbox{ and }
x\in X \mbox{ such that } J(x)=s.
\]
Conversely, if $(X,\Pi _X)$ is a Hamiltonian $\gm$-space with
momentum map $J: X\to G$, then $(X,\Pi _X)$ is a quasi-Poisson
manifold with momentum map $J:X\to G$, where the $G$-action on $X$
is given by
\[
\Phi (g,x)= (g,J(x))\cdot x,\mbox{ for }g\in G\mbox{ and }x\in X.
\]
\end{cor}


\begin{thebibliography}{99}
\bibitem{AK} A. Alekseev and Y. Kosmann-Schwarzbach, Manin Pairs and
moment maps, {\em J. Differential Geom.} {\bf 56} (2000) 133--165.

\bibitem{AKM} A. Alekseev, Y. Kosmann-Schwarzbach and E. Meinrenken,
Quasi-Poisson manifolds, {\em Canad. J. Math.} {\bf 54} (2002)
3--29.

\bibitem{BhaskaraV}
K.H. Bhaskara and K. Viswanath, {\em Poisson algebras and Poisson
manifolds}. Research Notes in Mathema\-tics, 174, Pitman, London,
1988.

\bibitem{CF} A.S. Cattaneo, G. Felder, Poisson sigma models and
symplectic groupoids, in {\it Quantization of Singular Symplectic
Quotients} (N.P. Landsman, M. Pflaum, and M. Schlichenmeier,
eds.), {\it Progr. in Math.} \textbf{198} (2001), 41--73.

\bibitem{CX}
A.~S. Cattaneo and P. Xu, Integration of twisted Poisson
structures, {\em J. Geom. Phys.} {\bf 49} (2004) 187--196.

%\bibitem{CDW} A. Coste, P. Dazord and A. Weinstein,
%Groupo{\"\i}des symplectiques, Publications du D{\'e}partement de
%Math{\'e}matiques de l'Universit{\'e} de Lyon, {I}, {\bf 2/A}
%(1987) 1--65.

\bibitem{Co}
T. Courant, Dirac manifolds, {\em Trans. A.M.S.} {\bf 319} (1990)
631--661.

\bibitem{Crainic} M. Crainic and R.L. Fernandes, Integrability of
Lie brackets, {\em Ann. of Math.} {\bf 157} (2003) 575--620.

%\bibitem{CF2} M. Crainic and R.L. Fernandes, Integrability of Poisson
%brackets, {\em to appear in J. Differential Geom.}.

\bibitem{Drinfeld} V. Drinfeld, Quasi-Hopf algebras, {\it
Leningrad Math. J.} \textbf{1} (1990) 1419--1457.

\bibitem{GU1}
J. Grabowski and P. Urbanski, Lie algebroids and Poisson-Nijenhuis
structures,  {\sl Rep. Math. Phys.} {\bf 40} (1997) 195--208.

\bibitem{GU}
J. Grabowski and P. Urbanski, Tangent and cotangent lifts and
graded Lie algebras associated with Lie algebroids, {\sl Ann.
Global Anal. Geom.} {\bf 15} (1997) 447--486.

\bibitem{H} H. Reinhard, \'Equations diff\'erentielles, {\em Fondements 
et applications}, Proceedings of Symposia in Pure Mathematics, 
38. Dunod, Paris, 1982.

\bibitem{Kosmann:1991}
Y.~Kosmann-Schwarzbach, Quasi-big\`ebres de Lie et groupes de Lie
quasi-Poisson, {\em C.R. Acad. Sci. Paris} {\bf 312} S\'erie I
(1991) 391--394.

\bibitem{Kosmann}
Y.~Kosmann-Schwarzbach, {E}xact Gerstenhaber algebras and Lie
bialgebroids, {\em Acta Appl. Math.} {\bf 41} (1995) 153--165.

\bibitem{LWX} Z.-J. Liu, A. Weinstein and P. Xu, Dirac structures and
Poisson homogeneous spaces, {\em Comm. Math. Phys.} {\bf 192}
(1998) 121--144.

\bibitem{Lu:1990} J.-H. Lu, Multiplicative and affine Poisson structures on
Lie groups, {\sl Ph.D. Thesis}, University of California at
Berkeley, (1990).

\bibitem{LuW:1990} J.-H. Lu and A. Weinstein,   Poisson Lie groups,
dressing transformations, and Bruhat decompositions, {\em J. Diff.
Geom.} {\bf 31} (1990) 501--526.

\bibitem{LW} J.-H. Lu and A. Weinstein, Groupo\"{\i}des symplectiques
doubles des groupes de Lie-Poisson, {\it C.R. Acad. Sci. Paris
S\'er. I Math.} \textbf{309} (1989) 951--954.

\bibitem{Mackenzie:LGLADG}
K.~Mackenzie, {\em Lie groupoids and {L}ie algebroids in
differential geometry}, London Mathematical Society Lecture Note
Series, no.~124. Cambridge University Press, 1987.

\bibitem{MackenzieX:1994}
K.~C.~H. Mackenzie and P. Xu, Lie bialgebroids and {P}oisson
groupoids, {\em Duke Math.~J.} {\bf 73} (1994) 415--452.

\bibitem{MackenzieX:1998} K.~C.~H. Mackenzie and P. Xu,
 Classical lifting processes and multiplicative vector fields,
 {\em Quarterly J. Math.} {\bf 49} (1998) 59--85.

\bibitem{MackenzieX:2000}
K.~C.~H. Mackenzie and P. Xu,
 Integration of Lie bialgebroids,
 {\em Topology} {\bf 39} (2000) 445--467.

\bibitem{Moerdijk-Mrcun} I. Moerdijk and J. Mr\v cun, On
integrability of infinitesimal actions,
 {\em Amer. J. Math.} {\bf 124} (2002) 567--593.

\bibitem{Roy}
D. Roytenberg, Quasi-Lie bialgebroids and twisted Poisson
manifolds,
 {\em Lett. Math. Phys.} {\bf 61} (2002) 123--137.

\bibitem{SW}
P. Severa and  A. Weinstein, Poisson geometry with a 3-form
background, {\em Prog. Theor. Phys. Suppl.} {\bf  144} (2001)
145--154.

\bibitem{Weinstein:1988}
A.~Weinstein,
 Coisotropic calculus and {P}oisson groupoids,
 {\em J. Math. Soc. Japan} {\bf 40} (1988) 705--727.

\bibitem{Weinstein:1990}
A. Weinstein, Affine {P}oisson structures, {\em Internat. J.
Math.} {\bf  1} (1990) 343--360.

\bibitem{WX}
A. Weinstein and P. Xu,
Extensions of symplectic groupoids and quantization,
{\em J. Reine. Angew. Math.} {\bf 417} (1991), 159-189.


\bibitem{Xu:1995}
P. Xu, On {P}oisson groupoids, {\em Internat. J.~Math.} {\bf 6}
(1995) 101--124.

\bibitem{Xu:1999}
P. Xu, Gerstenhaber algebras and BV-algebras in Poisson geometry,
{\em Comm. Math. Phys.} {\bf 200} (1999) 545--560.

\end{thebibliography}
\end{document}